\documentclass[1p]{elsarticle}
\journal{Stochastic Processes and their Applications}
\usepackage[colorlinks=true,linkcolor=blue,pdfborder={0 0 0}]{hyperref}

\usepackage{xcolor}
\usepackage{ifthen}
\usepackage{cmap}
\usepackage[utf8x]{inputenc}
\usepackage[T1]{fontenc}
\usepackage[english]{babel}
\usepackage{latexsym}
\usepackage{bbm}
\usepackage{amssymb,amsmath}
\usepackage{mathtools}
\usepackage{stmaryrd}
\usepackage{amsthm}
\usepackage{comment}
\newcommand{\markupdraft}[2]{
    \ifthenelse{\equal{#1}{display}}{#2}{}
    \ifthenelse{\equal{#1}{color}}{\color{#2}}{}
}

\newcommand{\newcolored}[3][]{{\markupdraft{color}{#2}#3}
    \ifthenelse{\equal{#1}{}}{}{\markupdraft{display}{{\color{yellow!70!black}[#1]}}}} 

\newcommand{\del}[2][]{{\markupdraft{display}{{\color{orange}[removed: "#2"[#1]]}}}} 
\newcommand{\new}[2][]{\newcolored[#1]{blue!90!green!90!black}{#2}}
    
\renewcommand{\del}[2][]{{}} 
\renewcommand{\markupdraft}[2]{}

\newtheorem{theorem}{Theorem}
\newtheorem{proposition}[theorem]{Proposition}
\newtheorem{corollary}[theorem]{Corollary} 
\newtheorem{lemma}[theorem]{Lemma} 
 
\theoremstyle{definition} 
 
\newtheorem{remark}[theorem]{Remark}

\DeclareMathOperator{\diag}{diag}
\DeclareMathOperator{\Tr}{Tr}

\DeclareMathOperator*{\argmin}{argmin}

\DeclareMathOperator{\E}{\mathbb{E}} 

\DeclarePairedDelimiterX{\inner}[2]{\langle}{\rangle}{#1, #2}
\DeclarePairedDelimiter{\norm}{\lVert}{\rVert}
\DeclarePairedDelimiter{\abs}{\lvert}{\rvert}

\providecommand{\T}{\mathrm{T}} 
\providecommand{\N}{\mathbb{N}} 
\providecommand{\X}{\mathbb{X}} 
\providecommand{\R}{\mathbb{R}} 
\providecommand{\rmd}{\mathrm{d}} 
\providecommand{\indicator}[1]{\mathbbm{I}\{#1\}} 
\providecommand{\1}[1]{\indicator{#1}} 

\renewcommand{\geq}{\geqslant}
\renewcommand{\leq}{\leqslant}

\renewcommand{\epsilon}{\varepsilon}

\renewcommand{\phi}{\varphi}
\let\oldsqrt\sqrt
\renewcommand{\sqrt}[2][]{\oldsqrt[\leftroot{-3}\uproot{3}#1]{\mathstrut #2}}

\providecommand{\deltat}{\alpha}

\newcommand{\flow}{\varphi}
\providecommand{\mm}{m}
\providecommand{\vv}{v}

\providecommand{\stheta}{\theta^{*}}

\providecommand{\ww}{w}

\providecommand{\xstar}{x^*}

\providecommand{\tildeY}{\tilde Y}

\providecommand{\ind}[1]{\indicator{#1}}

\providecommand{\asm}[1]{\textnormal{#1}} 

\newcommand{\Ratio}{c_\mathrm{ratio}}%

\newcommand{\wrt}{w.r.t.\ }
\newcommand{\dimtheta}{\ensuremath{{\mathrm{dim}(\theta)}}}

\newcommand{\sfield}{\mathcal{F}}

\newcommand{\qftle}{q_{\theta}^{<}}
\newcommand{\qfteq}{q_{\theta}^{=}}
\newcommand{\qftnle}[1]{q_{\theta_{#1}}^{<}}
\newcommand{\qftneq}[1]{q_{\theta_{#1}}^{=}}

\newcommand{\trinomial}{P_T}
\newcommand{\potential}{\Psi}
\newcommand{\Aone}{\asm{A1}}
\newcommand{\Atwo}{\asm{A2}}
\newcommand{\Athr}{\asm{A3}}
\newcommand{\Afou}{\asm{A4}}
\newcommand{\Bone}{\asm{B1}}
\newcommand{\Btwo}{\asm{B2}}
\newcommand{\DeltaBone}{\Delta_\Bone}
\newcommand{\DeltaBtwo}{\Delta_\Btwo}
\newcommand{\DeltaAone}{\Delta_\Aone}
\newcommand{\DeltaAtwo}{\Delta_\Atwo}
\newcommand{\DeltaAthr}{\Delta_\Athr}
\newcommand{\DeltaAfou}{\Delta_\Afou}
\newcommand{\CR}{\gamma_{\deltat}} 
\newcommand{\optN}{N_{\deltat}} 
\newcommand{\remCR}{{\bar \gamma_{\deltat}}}
\newcommand{\deltatset}{\Lambda} 

\newcommand{\EventU}{\Omega}
\newcommand{\supU}{\zeta}

\newcommand{\bigO}{O}

\newcommand{\eye}{I}
\newcommand{\cl}{\text{cl}}

\renewcommand{\S}{Section~}

\begin{document}

\providecommand{\correct}[1]{\textcolor{blue}{#1}}
\renewcommand{\correct}[1]{#1}

\sloppy

\begin{frontmatter}
\title{An ODE Method to Prove the Geometric Convergence of Adaptive Stochastic Algorithms\tnoteref{mytitlenote}}
\tnotetext[mytitlenote]{This is the author version of the paper accepted to Stochastic Processes and their Applications. DOI:   \url{https://doi.org/10.1016/j.spa.2021.12.005}.}

\author[tsukuba,riken]{Youhei Akimoto\corref{mycorrespondingauthor}}\ead{akimoto@cs.tsukuba.ac.jp}
\author[inria]{Anne Auger}\ead{anne.auger@inria.fr}
\author[inria]{Nikolaus Hansen}\ead{nikolaus.hansen@inria.fr}
\address[tsukuba]{University of Tsukuba, 305-8573 Ibaraki, Japan}
\address[riken]{RIKEN Center for Advanced Intelligence Project, 103-0027 Tokyo, Japan}
\address[inria]{Inria and CMAP, Ecole Polytechnique, IP Paris, CNRS, 91128 Palaiseau, France}
\cortext[mycorrespondingauthor]{Corresponding author}

\begin{abstract}
We consider stochastic algorithms \correct{derived} from methods for solving deterministic optimization problems, especially comparison-based algorithms derived from stochastic approximation algorithms with a constant step-size. We develop a methodology for proving geometric convergence of the parameter sequence $\{\theta_n\}_{n\geq 0}$ of such algorithms. We employ the ordinary differential equation (ODE) method, which relates a stochastic algorithm to its mean ODE, along with a Lyapunov-like function $\potential$ such that the geometric convergence of $\potential(\theta_n)$ implies---in the case of an optimization algorithm---the geometric convergence of the expected distance between the optimum and the search point generated by the algorithm.
  We provide two sufficient conditions for $\potential(\theta_n)$ to decrease at a geometric rate: $\potential$ should decrease ``exponentially'' along the solution to the mean ODE, and the deviation between the stochastic algorithm and the ODE solution (measured by $\potential$) should be bounded by $\potential(\theta_n)$ times a constant.  We also provide practical conditions under which the two sufficient conditions may be verified easily without knowing the solution of the mean ODE.
  Our results are any-time bounds on $\potential(\theta_n)$, so we can deduce not only the asymptotic upper bound on the convergence rate, but also the first hitting time of the algorithm.
  The main results are applied to a comparison-based stochastic algorithm with a constant step-size for optimization on continuous domains.
\end{abstract}

\begin{keyword}
  adaptive stochastic algorithm \sep
  comparison-based algorithm \sep
  geometric convergence \sep
  ordinary differential equation method \sep
  Lyapunov stability \sep
  optimization
\end{keyword}
\end{frontmatter}

\section{Introduction}
\label{sec:intro}

The ordinary differential equation (ODE) method is a standard technique to prove the convergence of stochastic algorithms of the form 
\begin{equation}
	\theta_{n+1} = \theta_{n} + \deltat F_n \enspace, \label{eq:algo-original}
\end{equation}
where $\{\theta_n\}_{n\geq0}$ is a sequence of parameters taking values in $\Theta \subseteq \R^\dimtheta$;  $\dimtheta$ is a positive integer;
$\{F_n\}_{n\geq0}$ is a sequence of random vectors in $\R^\dimtheta$; and
$\deltat > 0$ is the step size (or learning rate), which can depend on the time index $n$ \cite{Kushner2003book,Ljung1977tac,Borkar2008book,benaim2006stochastic2,benaim2005stochastic}.
The method connects the convergence of the stochastic algorithm to that of the solutions of the underlying ``mean'' ODE
\begin{equation}
\frac{\rmd \theta}{\rmd t} = F(\theta), \quad \theta(0) = \theta_0 \enspace,
\label{eq:associated-ode}
\end{equation}
where $F(\theta_n) = \E[F_n \mid \mathcal{F}_n ]=\E[F_n \mid \theta_n ]$ is the conditional expectation given the natural filtration $\{\mathcal{F}_n\}_{n \geq 0}$ associated to $\{\theta_n\}_{n \geq 0}$.
Here, we assume that $F$ is well defined and that the dependency of $F_n$ on the past is only through $\theta_n$.
The stochastic algorithm \eqref{eq:algo-original} provides a stochastic approximation of the solution of \eqref{eq:associated-ode}.
Often, the error between the stochastic algorithm \eqref{eq:algo-original} and the solution of the mean ODE is controlled by making $\alpha$ time dependent and taking a sequence of learning rates that decreases to zero, but not too quickly \cite{Kushner2003book,Ljung1977tac,Borkar2008book,Borkar:2000vd}.

In this paper, we explore the use of the ODE method to prove the convergence of some specific algorithms that, at the most abstract level, are of the form \eqref{eq:algo-original} but where \emph{geometric convergence} of a function of $\theta_n$ occurs.
For concreteness, consider the following example, an algorithm arising in the context of the optimization of a (black-box) function $f: \R^d \to \R$.
The algorithm state is $\theta_n = (m_n,\correct{\sigma_n}) \in \R^d \times \R_+ = \Theta$, where $m_n$ encodes the mean of a multivariate normal sampling distribution and $\correct{\sigma_n}$ its standard deviation.
The updates for $\theta_n$ read
\begin{equation}
  \begin{split}
    \mm_{n+1} &=  \mm_{n} + \deltat \sum_{i=1}^{\lambda} W(i;x_{n,1},\dots,x_{n,\lambda}) \Ratio (x_{n,i} - \mm_n) \enspace, \\
    \correct{\sigma_{n+1}} &=  \correct{\sigma_{n} + \deltat \sum_{i=1}^{\lambda} W(i;x_{n,1},\dots,x_{n,\lambda}) \left(\frac{1}{2d \sigma_n }\bigl( \norm{x_{n,i} - \mm_n}^2 - d \sigma_n^2\bigr)\right)} \enspace,
  \end{split}
\label{eq:es-igo-intro}
\end{equation}
where $(x_{n,i})_{1 \leq i \leq \lambda}$ are candidate solutions sampled independently from a multivariate normal distribution $\mathcal{N}(\mm_n, \correct{\sigma_n^2}\eye)$ with mean $\mm_n$ and overall \correct{standard deviation $\sigma_n$}; $W(i;x_{n,1},\dots,x_{n,\lambda})$ are weights assigned to the candidate solutions that are decreasing if the candidate solutions are ordered according to their $f$-values (i.e., better solutions have larger weights); and $\Ratio$ is a constant ensuring different learning rates for the mean and \correct{standard deviation} updates.

Algorithm \eqref{eq:es-igo-intro} is termed comparison-based because it uses objective-function values only through \emph{comparisons} that rank the candidate solutions.
Comparison-based algorithms are invariant to strictly increasing transformations of the objective function.
Therefore, they do not assume convexity or even continuity of the objective function, whereas most gradient-based algorithms require strong convexity of the function to guarantee geometric convergence.
Algorithm \eqref{eq:es-igo-intro} is (mainly empirically) known to converge \emph{geometrically} to a local minimum \new{$x^* \in \R^d$} of $f$ on wide classes of functions for a fixed learning rate $\alpha$.
Both the mean and the \correct{standard deviation} converge geometrically (towards a local optimum $x^* \in \R^d$ of $f$ and zero, respectively).
Algorithm~\eqref{eq:es-igo-intro} is a simplified version of a state-of-the art stochastic search algorithm, the Covariance Matrix Adaptation Evolution Strategy (CMA-ES), which adapts a full covariance matrix instead of \correct{$\sigma_n$} \cite{Hansen2001ec,Hansen2003ec,Hansen:2013vt,AH2019ecj}.
Despite successful applications of such randomized search algorithms \cite{Salimans2017}, proofs of their convergence are difficult to achieve. 

Note that is essential not only to prove that the algorithms converge, but also to characterize the rates at which they do so.

In this context, we extend the ODE method to be able to prove the geometric convergence of algorithms of the form \eqref{eq:es-igo-intro}. 
Our setting differs from those where the ODE method is usually employed, which makes extending the ODE method non-trivial.
Typically, we must cope with the following issues (some of which have been addressed individually in the literature) simultaneously:
\begin{enumerate}
\item
  The step-size $\deltat$ cannot decrease to zero, or the geometric convergence is jeopardized.
  Indeed, in the example algorithm described above, we empirically observe that the distance $\norm{X_n - x^*}$ from a candidate solution to the optimum decreases geometrically to zero with the factor proportional to $\exp(- \deltat)$ for sufficiently small $\deltat$. If $\deltat \to 0$, geometric convergence can no longer occur.
  Previous works using the ODE method generally considered a sequence of step-sizes decreasing to zero, because most stochastic approximation methods with constant step-sizes are not guaranteed to converge. (Some studies that did consider a constant step-size will be discussed later.)
  
\item
  All the (uncountably many) boundary points of the domain $\Theta$ are typically equilibrium points of the mean ODE. We want convergence towards one of them only: convergence to only one specific equilibrium point means that the underlying optimization algorithm has found the optimum while convergence to the others would mean convergence towards non-critical point.
  In the above example, all the boundary points $\theta = (\mm, 0)$ for any $\mm \in \R^d$ are equilibrium points of the associated ODE, so any neighborhood of the optimal parameter $\theta^*$ (equal to $(x^*,0)$ in the example) contains uncountably many equilibrium points.
\item
  Geometric convergence of the algorithm does not usually mean geometric convergence of $\theta_n$ with the Euclidean distance.
  We therefore introduce a function $\potential: \Theta \to \R^+$, the geometric convergence of which implies the geometric convergence of the quantities in which we are interested.
  Thus, in the above example, we would like to show the geometric convergence of the distance between the generated solutions $x_{n,i}$ and the optimal solution $x^*$ in some stochastic sense, so $\potential$ must be chosen so that its geometric convergence implies that of the distance.
\end{enumerate}

Having equilibrium points on the boundary is not a critical issue.
Previous works (e.g., \cite{Buche2001sicon}) addressed such cases.
However, having multiple equilibrium points that are \emph{all connected} complicates the analysis.
We want to prove the convergence only towards a specific equilibrium point; convergence towards the other points mean a failure of the optimization algorithm.
Previous studies \cite{Ljung1978,Pemantle1990,Brandiere1998sicon,Fang2000ieeetac} have addressed multiple equilibria. However, in those studies, the objective was to show that the algorithm converges towards a set of equilibria, whereas we want it to converge to a single point.
Bena\"{i}m \cite{Benaim1996} showed that the limit set of the stochastic process \eqref{eq:algo-original} is internally chain transitive for the mean dynamics \eqref{eq:associated-ode}, but step-sizes were assumed to decrease to zero.
  Therefore, the combination of the first and second difficulties makes it hard to apply existing results.
The third point of the list is critical, especially since any neighborhood of the optimal parameter typically has uncountably many equilibrium points.
Hence we cannot choose (for instance) the Euclidean distance to the optimal parameter as the potential function $\potential$. 

Our approach consists in finding a function $\potential: \Theta \to \R_{\geq 0}$  that satisfies conditions \Aone{} and \Atwo{} presented below so that $\potential(\theta_n)$ converges geometrically to zero for an appropriate choice of the step-size $\alpha$.
The function $\potential$ can be thought as the distance to a desired parameter set such that, if $\potential(\theta_n)$ goes to zero geometrically, we can conclude that the algorithm we investigate also converges geometrically, i.e.,~$\potential(\theta_n)$ is an upper bound on that quantity, the geometric convergence of which we seek to prove.
The first condition on $\potential$ states that $\potential$ should decrease along the trajectories of the solution of the ODE \eqref{eq:associated-ode} (behaving similarly to the Lyapunov functions used for proving the stability of ODEs), and the decrease should be geometric.
More precisely:
\begin{itemize}
\item[\Aone{}] $\exists \DeltaAone: \R_{\geq 0} \to \R_{\geq 0}$ nonincreasing such that $\DeltaAone(t) \downarrow 0$ as $t \uparrow \infty$, and for any $\theta \in \Theta$ and any $t \in \R_{\geq 0}$, 
\begin{equation}\label{eq:A1intro}
\potential(\flow(t; \theta)) \leq \DeltaAone(t) \potential(\theta) \enspace,
\end{equation}
where $\flow(t; \theta)$  is the solution of \eqref{eq:associated-ode} at time $t$ and $\flow(0; \theta) = \theta$.
\end{itemize}
The second condition is to control the deviation between the stochastic algorithm and the solution of the ODE: 
\begin{itemize}
\item[\Atwo{}] $\exists \DeltaAtwo: \R_{+} \times \R_{+} \to \R_{+}$ nondecreasing \wrt each argument such that $\DeltaAtwo(\deltat, T) \downarrow 0$ as $\deltat \downarrow 0$ for any fixed $T > 0$, and for any $N \in \N_{+}$ and $\theta_{0} \in \Theta$,
\begin{equation}\label{eq:A2intro}
\E_0[\potential(\theta_{N})] \leq \potential(\flow(N \deltat; \theta_{0}))  + \DeltaAtwo(\deltat, N\deltat) \potential(\theta_{0}) \enspace ,
\end{equation}
where $\E_{0}$ denotes the conditional expectation given $\theta_{0}$.
\end{itemize}
Note that these conditions are stronger than the assumptions typically made in standard stochastic-approximation settings.
In particular, the deviation term (second term on the right-hand side of \eqref{eq:A2intro}) must be proportional to $\potential(\theta_0)$; by contrast, standard assumptions lead to \eqref{eq:ode-lemma} below, where the corresponding term is a constant over $\theta_0 \in \Theta$.
Our goal is to obtain the convergence rate of the algorithm, not just to prove its convergence. 

The flow solution $\flow$ of the ODE comes into play in both conditions \eqref{eq:A1intro} and \eqref{eq:A2intro}.
However, only rarely is the solution of the ODE known explicitly.
Hence, we provide some practical conditions under which it is possible to verify the sufficient conditions without knowing $\varphi$.
The condition developed in Theorem~\ref{lem:prac-one} (which implies \eqref{eq:A1intro}) is similar to the conditions for obtaining exponential stability of equilibrium points of ODEs in Lyapunov's theory.
It states that the upper Dini directional derivative of $\potential$ in the direction of $F(\theta)$ has to be smaller than or equal to a negative constant times $\potential(\theta)$.
The additional conditions developed in Theorem~\ref{lem:ode-igo} (implying \eqref{eq:A2intro}) are based on the conditions for the expected Euclidean distance between the stochastic algorithm and the solution of the underlying ODE being bounded by a constant (see Lemma~1 in Chapter~9 of \cite{Borkar2008book}).
We replace the constant with a constant times $\potential(\theta)$.
We illustrate how to use the different results (and in particular the practical conditions) on \correct{the step-size adaptive evolution strategy described in \eqref{eq:es-igo-intro}}.

This paper is organized as follows:
After summarizing the notations employed throughout, we describe in \S\ref{sec:ex} some rank-based or comparison-based stochastic search algorithms to which the methodology of the paper can be applied, including a running example (a step-size adaptive evolution strategy).
In \S\ref{sec:ode}, we demonstrate a simple technique of the ODE method and illustrate the difficulty of proving the geometric convergence of stochastic algorithms via the ODE method.
In \S\ref{sec:linconv}, we provide a sufficient condition for $\potential(\theta_n)$ to converge globally and geometrically towards zero and derive the first-hitting-time bound from it.
We prove practical conditions for verifying the sufficient condition of our main theorem in \S\ref{sec:pra}, followed by the application of these practical conditions to \correct{the step-size adaptive evolution strategy and the proof of its geometric convergence in \S\ref{sec:appli}.}
The extension of the main theorem to cover the cases where the ODE associated to the stochastic algorithm has multiple local attractors is discussed in \S\ref{sec:llc}. We conclude this paper in \S\ref{sec:conc}.

\paragraph{Notations}

Let $\ind{A}$ be the indicator function, with value $1$ if event $A$ occurs, $0$ otherwise.
Let $\E_{\omega \sim P}$ denote the expectation taken over $\omega$ having probability measure $P$, i.e., $\E_{\omega \sim P}[g(\omega)] = \int g(\omega) P(\rmd \omega)$.
When the expectation is taken over all random variables appearing in the expression, we will write just $\E$ if no confusion will result. Let $\E_{n} = \E[\cdot\mid\sfield_{n}]$ denote the conditional expectation given a filtration $\{\sfield_n\}_{n\geq0}$.
For any event $A$, let $\Pr[A] = \E[\ind{A}]$ be the probability of $A$ occurring.
Let $\Pr_n[A] = \E_n[\ind{A}]$ be the conditional probability given filtration $\sfield_{n}$.

Let $\R$ be the set of real numbers.
Let $\R_{\geq 0}$ and $\R_{+}$ be the sets of nonnegative and positive real numbers, respectively. 
Let $\N$ and $\N_{+}$ be the sets of nonnegative and positive integers, respectively. 
For $a \in \R$ and $b \in \R$ satisfying $a < b$, the open, closed, left open, right open intervals are denoted $(a, b)$, $[a, b]$, $(a, b]$, and $[a, b)$, respectively.
For $a \in \N$ and $b \in \N$ satisfying $a \leq b$, $\llbracket a, b \rrbracket$ denotes the set of integers between $a$ and $b$, inclusive. 
Let $\norm{x}$ be the Euclidean norm for any real vector $x$, and $\norm{x}_Q = \norm{Q^{1/2} x}$ be the Mahalanobis norm given a positive definite symmetric matrix $Q$. 
For any real value $a$, $\abs{a}$ denotes the absolute value. 

\section{Example Algorithms}\label{sec:ex}

The methodology presented in this paper is motivated by the class of comparison-based or rank-based stochastic search algorithms.
While we have sketched an example of such algorithms in the introduction, we define more generally in this section the class of algorithms and derive the function $F$ in \eqref{eq:associated-ode}.

The objective is to minimize (without loss of generality) a function $f: \X \to \R$, where $\X$ is an arbitrary search space.
The algorithm updates $\theta_n \in \Theta$ parametrizing a family of probability distributions $(P_\theta)_{\theta \in \Theta}$ defined over $\X$.
More precisely, at each iteration $n \geq 0$, multiple candidate solutions $x_{n,1},\dots,x_{n,\lambda}$ ($\lambda \geq 2$) are generated from the probability distribution $P_{\theta_n}$.
Their objective values $f(x_{n,i})$ are evaluated and higher weight values are assigned to candidate solutions with smaller $f$-values.
The weight value assigned to each $x_{n,i}$ is defined as follows:
Let $w_1,\dots,w_\lambda$ be real, predefined constants that are nonincreasing, i.e., $w_1 \geq w_2 \geq \dots \geq w_\lambda$.
Define a function $W: \llbracket1,\lambda\rrbracket \times \X^{\lambda} \to \R$ as
\begin{equation}
  \begin{split}
    &\textstyle W(i; x_{1},\dots,x_\lambda) = \sum_{j = {k_i}}^{{k_i+l_i}} \frac{ w_{j+1} }{{l_i}+1}, \\
    &\textstyle \text{where }
    {k_i} = \sum_{j\neq i} \1{f(x_j) < f(x_i)},\ {l_i} = \sum_{j \neq i} \1{f(x_j) = f(x_i)}
    \enspace. 
  \end{split}
  \label{eq:weight-rank}
\end{equation}
Note that $\sum_{i=1}^\lambda W(i; x_1,\dots,x_\lambda) = \sum_{i=1}^{\lambda} w_i$ for any combination of $\{x_i\}_{i=1,\dots,\lambda}$.
Then $W(i; x_{n,1},\dots,x_{n,\lambda})$ is the value assigned to $x_{n,i}$. 
If all the candidate solutions have distinct objective values, we have $W(i; x_{n,1},\dots,x_{n,\lambda}) = w_{k+1}$, where $k$ is the number of points with better $f$-values than $x_{n,i}$ among $\{x_{n,j}\}_{j\in\llbracket 1, \lambda\rrbracket}$ as defined in \eqref{eq:weight-rank}.
(In other words, $w_{k}$ is assigned to the $k$th-best point.)
Using the set of pairs of candidate solutions and their weight values, the parameter $\theta_n$ of the probability distribution is updated.
Let $g: \X \times \Theta \to \R^{\dimtheta}$.
The direction to adjust the parameter is defined by the sum of the products of the weight values and $g(x_{n,i}; \theta_n)$, i.e.,
\begin{equation}
\textstyle \theta_{n+1} = \theta_n + \deltat \sum_{i=1}^{\lambda} W(i; x_{n,1},\dots,x_{n,\lambda}) g(x_{n,i}; \theta_n) \enspace.
\label{eq:algo-rank}
\end{equation}
The algorithm \eqref{eq:es-igo-intro} sketched in the introduction follows the previous update equation with $\theta_n = (\mm_n, \correct{\sigma_n}) \in \Theta = \R^d \times \R_+$ and $g: \R^d \times \Theta \to (\R^d \times \R)$, defined as
\begin{equation}\label{eq:adaptES}
\textstyle g\left(x; \theta=(\mm,\correct{\sigma})\right) = \left(\Ratio (x - \mm), \
                \correct{ \frac{1}{2d\sigma} \big( \norm{x - \mm}^2 - d \sigma\big) }  \right) \enspace.
\end{equation}%
We see that the update equation for $\theta_n$ in \eqref{eq:algo-rank} depends on $f$ only through the weight functions, which themselves depend on $f$ only through the ranking of the candidate solutions (see \eqref{eq:weight-rank}); this explains the name ``comparison'' or ``rank-based'' algorithm.
Hence, such algorithms are invariant under rank-preserving transformations.
In other words, consider a function $h: \R \to \R$ that is strictly monotonically increasing: we obtain the same update \eqref{eq:algo-rank} for the algorithm optimizing $f$ or $h \circ f$. 

For this class of algorithms we can explicitly write the function $F(\theta) = \E\del{_n}{}[F_n \mid \theta_n = \theta]$ using the following proposition, the proof of which is included in \ref{apdx:prop:exp}.
\correct{
\begin{proposition}
\label{prop:exp}
  Given $\theta_n \in \Theta$, let $X_{n,1},\dots,X_{n,\lambda}$ be i.i.d.\ following $P_{\theta_n}$ and define $W$ as in \eqref{eq:weight-rank}.
  Let $\qftnle{n}(s)$ and $\qftneq{n}(s)$ be the probabilities of $\{f(x) < s\}$ and $\{f(x) = s\}$, respectively, when $x \sim P_{\theta_n}$.
  For any $P_{\theta_n}$-integrable function $g: \X \times \Theta \to \R^{\dimtheta}$, the random vector $Y_n = \sum_{i=1}^{\lambda} W(i; X_{n,1},\dots,X_{n,\lambda}) g(X_{n,i}; \theta_n)$ is $P_{\theta_n}$-integrable.
  Moreover, there exists a function $F$ such that the conditional expectation $F(\theta_n) = \E_n[Y_n] = \E[Y_n \mid \theta_n]$ is 
\begin{equation*}
F(\theta_n) = \int_{\X} u(\qftnle{n}(f(x)), \qftneq{n}(f(x))) g(x; \theta_n) P_{\theta_n}(\rmd x) \enspace.
\end{equation*}
Here,  $u: [0, 1]^2 \to \R$ is the function defined by 
\begin{equation}
u(p, q) = \sum_{k=0}^{\lambda-1} \sum_{l=0}^{\lambda-k-1} \Biggl( \sum_{j = k}^{k+l} \frac{ \ww_{j+1} }{ l + 1 } \Biggr) \lambda P_{T}(\lambda-1, k, l, p, q)
\enspace,
\label{eq:exp:weight}
\end{equation}
where $P_{T}(n, k, l, p, q) = \frac{n!}{k!l!(n-k-l)!}p^kq^l(1 - p - q)^{n-k-l}$ is the trinomial probability mass function. 
\end{proposition}
}
\begin{remark}
\label{rem:prop:exp}
In Proposition~\ref{prop:exp}, if $\qftneq{n}(f(x)) = 0$ holds $P_{\theta_n}$-almost-surely, we can simplify the result as
\begin{equation}
F(\theta_n) = \int_{\X} u (\qftnle{n}(f(x))) g(x; \theta_n) P_{\theta_n}(\rmd x)
\enspace,
\label{eq:meanfield-s}
\end{equation}
where $u: [0, 1] \to \R$ is defined by 
\begin{equation}
u(p) = \sum_{k=0}^{\lambda-1} \ww_{k+1} 
\lambda P_{B}(\lambda-1, k, p)
\enspace,
\label{eq:binom}
\end{equation}
and $P_{B}(n, k, p) = \frac{n!}{k!(n-k)!}p^k(1 - p)^{n-k}$ is the binomial probability mass function. 
\end{remark}


\subsection{Step-size Adaptive Evolution Strategy}\label{sec:es}

Step-size adaptive evolution strategies (SSA-ES) \cite{Rechenberg1973,Schwefel1975,Hansen:2013wf} are rank-based optimization algorithms in a continuous domain, that is $\X = \R^d$.
We consider an example corresponding to a simplification of the state-of-the-art CMA-ES algorithm. (See \cite{Hansen:2013wf} for a recent overview of different standard methods.)
This is the method that was already sketched in the introduction.
At each iteration, the algorithm samples candidate solutions $x_{n,1}, \dots, x_{n,\lambda}$ from a multivariate Gaussian distribution $\mathcal{N}(\mm_n, \correct{\sigma_n^2} \eye)$ parameterized by $\theta_n = (\mm_n, \correct{\sigma_n})$, where $\mm_n \in \R^d$ represents the mean vector of the Gaussian distribution and \correct{$\sigma_n >0$ its overall standard deviation}.
The parameter space is then $\Theta = \R^d \times \R_{+}$. 
The parameter update follows \eqref{eq:es-igo-intro}, i.e., \eqref{eq:algo-rank} with $g$ in the form \eqref{eq:adaptES}, and $\{\theta_n\}_{n\geq 0}$ never leaves $\Theta$ for $\deltat < \correct{2}/\sum_{i=1}^{\lambda}w_i$ if $w_i \geq 0$ for all $i = 1,\dots,\lambda$.

\subsection{Population-Based Incremental Learning}\label{sec:pbil}

The population-based incremental learning (PBIL) algorithm \cite{Baluja1995icml} is a probability-model-based search algorithm for optimization of a function defined on a binary space, $f: \{0, 1\}^d \to \R$.
At each iteration, PBIL samples candidate solutions $x_{n,1}, \dots, x_{n,\lambda}$ from a multivariate Bernoulli distribution parameterized by $\theta \in \Theta = (0, 1)^d$, where the $i$th element of $\theta$, denoted by $[\theta]_i$, represents the probability of the $i$th component of $x$, denoted by $[x]_i$, being $1$. 
The probability parameter $\theta$ is then updated by the formula
\begin{equation}
\textstyle \theta_{n+1} = \theta_{n} + \deltat \sum_{i=1}^{\lambda} W(i;x_{n,1},\dots,x_{n,\lambda}) (x_{n,i} - \theta_n) \enspace.
\label{eq:pbil}
\end{equation}
It is easy to see that this algorithm is of the form \eqref{eq:algo-rank} and that the sequence $\{\theta_n\}_{n \geq 0}$ never leaves $\Theta$ for $\deltat < 1/\sum_{i=1}^{\lambda}w_i$ if $w_i \geq 0$ for all $i = 1,\dots,\lambda$.

\subsection{Information Geometric Optimization}
\label{sec:igo}

The information geometric optimization (IGO) \cite{Ollivier2017jmlr} is a generic framework of probability-model-based search algorithms for black-box optimization of $f: \X \to \R$ in an arbitrary search space $\X$.
It takes a parametric family of probability distributions $P_\theta$ on $\X$ as an input to the framework and provides a procedure to update the distribution parameters $\theta \in \Theta \subseteq \mathbb{R}^\dimtheta$.
\correct{
  IGO fits in the form \eqref{eq:algo-rank} with $g(x; \theta) = \tilde\nabla_\theta \ln (p_\theta(x))$, where $p_\theta(x)$ is the probability density of $P_\theta$ \wrt the reference measure $\mathrm{d}x$ and $\tilde\nabla_\theta \ln (p_\theta(x))$ is the so-called \emph{natural gradient} \cite{Amari1998nc} of the log-likelihood $\ln(p_\theta(x))$. The natural gradient is defined as the product $\mathcal{I}(\theta)^{-1}\nabla_\theta \ln(p_\theta(x))$ of the inverse of the Fisher information matrix $\mathcal{I}(\theta) = \int \nabla_\theta \ln (p_\theta(x)) \nabla_\theta \ln (p_\theta(x))^\T P_\theta(\mathrm{d}x)$ and the vanilla gradient $\nabla \ln(p_\theta(x)) = [\partial \ln(p_\theta(x))/\partial \theta_1,\dots,\partial \ln(p_\theta(x))/\partial \theta_\dimtheta]$.
  Both algorithms in Section~\ref{sec:es} (if $\Ratio = 1$) and Section~\ref{sec:pbil} are instantiations of IGO (with the Gaussian and Bernoulli distributions as $P_\theta$, respectively). 
}

\subsection{Policy Gradient with Parameter Exploration}
\label{sec:pgpe}

Consider the control task of an agent that has a parametrized policy, i.e., a control law, $u(s; w)$, where $s \in \mathcal{S}$ is the observed state of the agent, $w \in \mathcal{W}$ is the policy parameter, and $u:\mathcal{S} \to \mathcal{A}$ maps the observed state to an action.
The agent interacts with the environment; its next state is determined by the unknown transition probability kernel, $p(s' \mid s, a)$, where $s \in \mathcal{S}$ and $a = u(s; w)$ are the current state and the current action, and $s'$ is the next state.
Given the history, $h = (s_0, a_0, s_1, a_1, \dots, s_{T-1}, a_{T-1}, s_T)$, of states and actions during $T$ steps, the agent receives a reward, $r(h) \in \mathbb{R}$.
The objective of the learning task is to optimize $w$ so as to maximize the expected reward $f(w) = \mathbb{E}[r(h)]$.
The control law is often modeled by a multi-layer perceptron, choosing $w \in \mathcal{W} = \mathbb{R}^d$ to be the vector consisting of link weights between neurons.
Since the transition probability kernel, which can be deterministic, is unknown to the agent, the function $f: \mathbb{R}^d \to \mathbb{R}$ to be maximized forms a black-box objective.

Policy Gradient with Parameter Exploration (PGPE) \cite{Sehnke2010nn} takes a probability model $p_\theta(w)$ on $\mathcal{W}$ parameterized by $\theta \in \Theta \subseteq \mathbb{R}^\dimtheta$ and uses the expectation of $f(w)$ over $p_\theta(w)$ to define the expected objective $J(\theta) = \int f(w) p_\theta(w) \mathrm{d}w$.
PGPE then takes the vanilla gradient $\nabla J$ at $\theta_n$.
The vanilla gradient of $J(\theta)$ can be approximated by the Monte Carlo method as
\begin{align*}
  \nabla J(\theta)
  &= \textstyle \int f(w) \nabla_\theta \ln(p_\theta(w)) p_\theta(w) \mathrm{d}w \\
  &\approx \textstyle \frac{1}{\lambda} \sum_{i=1}^{\lambda} f(w_i) \nabla_\theta \ln(p_\theta(w_{i})) \\ 
  &\approx \textstyle \frac{1}{\lambda} \sum_{i=1}^{\lambda} \left(\frac{1}{N_h} \sum_{j=1}^{N_h} r(h_{i,j}) \right) \nabla_\theta \ln(p_\theta(w_{i})) 
  =: \widehat{\nabla J(\theta)}\enspace,
\end{align*}
where $\lambda$ is the number of samples, $w_i$ for $i = 1, \dots, \lambda$ are independent samples drawn from $p_\theta$, $N_h$ is the number of histories observed for each $w_i$, and $h_{i,j}$ are the histories generated by the agent with policy $u(\cdot; w_i)$.
If the transition in the environment is deterministic, $f(w_i) = r(h_{i,j})$ for any $j = 1, \dots, N_h$; hence, there is no need to set $N_h > 1$. The parameter update reads
\begin{equation*}
\theta_{n+1} = \theta_{n} + A_{n} \widehat{\nabla J(\theta_{n})} \enspace,
\end{equation*}
where $n = 1, \dots$ is the iteration counter, and $A_{n} \in \mathbb{R}^{\dimtheta \times \dimtheta}$ is a generalized step-size.

In references \cite{Sehnke2010nn,Zhao2012nn}, the independent Gaussian model is employed, where the parameter $\theta$ encodes the mean vector $m \in \R^{d}$ and the standard deviation in each coordinate $\sigma_i > 0$ for $i = 1, \dots, d$.
The generalized step-size is then the diagonal matrix such that $[A]_{i, i} = [A]_{d + i, d + i} = \alpha \sigma_{i}^2$.
We find that it is proportional to the inverse of the Fisher information matrix $\mathcal{I}(\theta)$ of $p_\theta$.
Indeed, $\mathcal{I}(\theta)$ for this model is a diagonal matrix with $i$th and $d+i$th diagonal elements $1/\sigma_i^2$ and $2/\sigma_i^2$, respectively.
Ignoring the factor of $2$, the parameter update reads
\begin{equation*}
\theta_{n+1} = \theta_{n} + \alpha \mathcal{I}^{-1}(\theta_{n}) \widehat{\nabla J(\theta_{n})} \enspace.
\end{equation*}
This method is equivalent to the IGO algorithm, except that the nonlinear transformation $f(w_i) \mapsto W(i; w_1,\dots,w_\lambda)$ in \eqref{eq:weight-rank} is not performed.

\subsection{Other Example\new{s}\del{ Algorithms}}\label{sec:otherexample}

The rank-based search algorithms of the form \eqref{eq:algo-rank} include probabilistic-model-based algorithms, such as cross-entropy optimization algorithms \cite{Boer2005aor} and estimation-of-distribution algorithms \cite{Lozano2006book}.
Though we are primarily interested in rank-based algorithms, our approach is not limited to them.
It also includes randomized derivative-free methods such as natural-evolution strategies (NES) \cite{wierstra2014jmlr} or adaptive simultaneous perturbation stochastic approximation (SPSA) \cite{spall2000adaptive}.
The search space $\X$ is not necessarily a continuous domain; it can be discrete as long as the probability distribution is parametrized by a continuous parameter vector as in the PBIL in \S\ref{sec:pbil}.

\section{Ordinary Differential Equation Methods}
\label{sec:ode}

ODE methods are widely used for proving the convergence of recursive stochastic algorithms of the form \eqref{eq:algo-original}.
We rewrite the algorithms as:
\begin{equation}
	\theta_{n+1} = \theta_{n} + \deltat (F(\theta_n) + M_n) \enspace, \label{eq:algo}
\end{equation}
where $M_n = F_n - F(\theta_n)$ is a martingale difference sequence.
Throughout the paper, we assume that $\theta_n$ never leaves its domain $\Theta$ for sufficiently small $\deltat$.
Given certain regularity conditions, it can be proven that the stochastic sequence $\{ \theta_n\}_{n \geq 0}$ converges towards the attractor set of the associated ODE \eqref{eq:associated-ode}: see, e.g., Chapter~2 of \cite{Borkar2008book}.
Here we demonstrate a simple ODE method that relates the behavior of a recursive stochastic algorithm \eqref{eq:algo} and its associated ODE \eqref{eq:associated-ode}.
To illustrate the basic principle of the approach, we pose a relatively strong assumption that allows us to obtain an easy proof, namely $\E_{n}[\norm{F_{n}}^{2}] \leq K^{2}$;
i.e., that the second moment of $F_n$ is bounded by a constant over the domain of $\theta$.
This assumption may not always be satisfied, especially when the domain of $\theta$ is unbounded.

Let $\phi: \R_{\geq 0} \times \Theta$ be the flow derived from $F$, i.e., $\phi(t; \theta_0)$ is an extended solution to the initial value problem $\frac{\rmd \theta}{\rmd t} = F(\theta)$ with $\theta = \theta_0$ at $t = 0$.
More precisely, $\phi(\cdot; \theta_0)$ is an absolutely continuous function defined on $\R_{\geq 0}$ that satisfies
\begin{equation}\label{eq:integral-form}
\phi(t; \theta_0) = \theta_0 + \int_{0}^{t} F(\phi(\tau; \theta_0)) \rmd \tau \enspace.
\end{equation} 
This definition of the solution is extended in the sense that it may have non-differentiable points in a set of zero measure.
The existence of a solution of the ODE $\frac{\rmd \theta}{\rmd t} = F(\theta)$ of the form \eqref{eq:integral-form} is verifiable without knowing the explicit solution by using, for example, Carath{\'e}odory's existence theorem.

Let $\deltat_n$ be a deterministic and possibly time-dependent step-size and define 
\begin{align}
t_{n,k} &= \sum_{i=n}^{n+k-1}\deltat_{i}
&\text{and} &&
\epsilon_{n,k} &= \sum_{i=n}^{n+k-1}\deltat_{i}^2
\enspace.
\end{align}
Given $\theta_n$ and $N \geq 1$, we consider $\theta_{n+N}$ as the approximation stemming from \eqref{eq:algo} of the solution $\phi(t_{n,N}; \theta_n)$ to the ODE at time $t = t_{n,N}$ with the initial condition $\flow(0; \theta_n) = \theta_n$. 
It follows from \eqref{eq:algo} that $\theta_{n+N} = \theta_n + \sum_{i=n}^{n+N-1} \deltat_{i} (F(\theta_{i}) + M_i)$ and from \eqref{eq:integral-form} that if the solution exists then it satisfies $\phi(t_{n,N}; \theta_n) = \theta_n + \int_{0}^{t_{n,N}} F(\phi(\tau; \theta_n)) \rmd \tau$ for any $n\in \N$ and $N \in \N_{+}$, hence the difference between $\theta_{n+N}$ and $\phi(t_{n,N}; \theta_n)$ can be expressed as
\begin{multline}
	\textstyle
    \theta_{n+N} - \phi(t_{n,N}; \theta_n) 
    = 
    \sum_{i=0}^{N-1}\deltat_{n+i} (F(\theta_{n+i}) - F(\phi(t_{n,i}; \theta_n))) \\
    \textstyle
    + \sum_{i=0}^{N-1}\int_{t_{n,i}}^{t_{n,i+1}} (F(\phi(t_{n,i}; \theta_n)) - F(\phi(\tau; \theta_n))) \rmd \tau +\sum_{i=0}^{N-1}\deltat_{n+i} M_{n+i} \enspace.
  \label{eq:ode-diff}
\end{multline}

Assume that $F$ is Lipschitz continuous with $L$ as the Lipschitz constant.
The second term on the right-hand side (RHS) of \eqref{eq:ode-diff} is the sum of the errors due to the time-discretization of the ODE solution.
It is not difficult to imagine that we can obtain an $\bigO(\epsilon_{n,N})$ bound for the second term because of the Lipschitz continuity.
The third term is the sum of the martingale difference noises.
By the uncorrelation of martingale difference sequences, we can obtain an $\bigO(\epsilon_{n,N}^{1/2})$ bound.
The first term is bounded as $\norm{F(\theta_{n+i}) - F(\phi(t_{n,i}; \theta_n))} \leq L \norm{\theta_{n+i} - \phi(t_{n,i}; \theta_n)}$ because of the Lipschitz continuity of $F$.
Then, we can apply the discrete Gronwall inequality \cite{Clark1987279} to obtain the following basic result, which can be derived by following a standard argument such as that in Chapter~2 of \cite{Borkar2008book}.
Its proof is found in the proof of Theorem~\ref{lem:ode-igo}, which is an extension of the following theorem:
\begin{theorem}\label{lem:ode}
Consider an algorithm of the form \eqref{eq:algo} on the domain $\Theta = \R^\dimtheta$ with deterministic and possibly time-dependent step-size $\deltat = \deltat_n$. 
Assume that $\E_{n}[\norm{F_{n}}^{2}] \leq K^{2}$ for a finite constant $K > 0$, where $F(\theta)$ is Lipschitz continuous with the Lipschitz constant $L$, i.e., $\norm{F(\theta_{a}) - F(\theta_{b})} \leq L \norm{\theta_{a} - \theta_{b}}$. Let $t_{n,k} = \sum_{i=n}^{n+k-1}\deltat_{i}$ and $\epsilon_{n,k} = \sum_{i=n}^{n+k-1}\deltat_{i}^2$.
Then, for any $\theta \in \Theta$ there exists a unique solution $\flow(\cdot; \theta): \R_{\geq 0} \to \Theta$ satisfying \eqref{eq:integral-form} and for any $n\in \N$ and $N \in \N_{+}$, 
\begin{equation}
	\sup_{0 \leq k \leq N} \E_n[\norm{\theta_{n+k} - \phi(t_{n,k}; \theta_n)}]
	\leq ( (L K/2) \epsilon_{n,N} + K \epsilon_{n,N}^{1/2}) \exp(L t_{n,N})
        \enspace.
	\label{eq:ode-lemma}
\end{equation}
\end{theorem}%

If the step-size satisfies $\sum_{n} \deltat_{n} = \infty$ and $\sum_{n} \deltat_n^{2} < \infty$, then for any fixed $N \geq 1$, $\epsilon_{n,N} \to 0$ and $t_{n,N} \to 0$ as $n \to \infty$.
Hence, we obtain with \eqref{eq:ode-lemma} that
\begin{equation}
	\lim_{n\to\infty}\sup_{0 \leq k \leq N}\E_n[\textstyle\norm{\theta_{n+k} - \phi(t_{n,k}; \theta_n)}]
	= 0\enspace.
	\label{eq:ode-decrease}
\end{equation}
Roughly speaking, the limit behavior of $\{\theta_{n+k}\}_{0 \leq k \leq N}$ for $n \to \infty$ follows the continuous-time trajectory $\phi(\cdot; \theta_n)$ for any fixed $N$.
Therefore, if the associated mean ODE has a \emph{global attractor} $\theta^*$ in the sense that $\lim_{t\to\infty}\phi(t; \theta) = \theta^*$ for any $\theta \in \Theta$, the sequence $\theta_n$ will converge toward it.
The algorithms that we are interested in (e.g., rank-based search algorithms) typically use a constant step-size and converge, but
geometric convergence is not observed if the step-size is decreasing.

If the step-size is constant over the time index $n$, for any $n \geq 0$ and $N \geq 1$, $\epsilon_{n,N} \to 0$ and $t_{n,N} \to 0$ as $\deltat \to 0$.
Therefore, \eqref{eq:ode-lemma} immediately implies that
\begin{equation}
	\lim_{\deltat\to0}\sup_{0 \leq k \leq N}\E_n[\norm{\theta_{n+k} - \phi(k\deltat; \theta_n)}]
	= 0\enspace.
	\label{eq:ode-const}
\end{equation}
Roughly speaking, the stochastic sequence $\{\theta_{n+k}\}_{k \in \llbracket 0, N \rrbracket}$ can follow the continuous-time trajectory $\phi(t; \theta_n)$ for $t \in [0, N \deltat]$ with arbitrary precision by taking a sufficiently small $\deltat$.
This is not only true in the limit but holds for any $n$. However, since $t_{n,N} = N\deltat$ and $\epsilon_{n,N} = N\deltat^2$, one can see from the RHS of \eqref{eq:ode-lemma} that the upper-bound found for the term $\sup_{0 \leq k \leq N}\E_n[\norm{\theta_{n+k} - \phi(k\deltat; \theta_n)}]$ increases infinitely as $N \to \infty$ for any fixed $\deltat > 0$.
Therefore, we do not get the convergence of $\{\theta_{n}\}_{n \geq 0}$ from this argument.

Theorem~3 in Chapter~9 of \cite{Borkar2008book}, for example, deals with a recursive algorithm with a constant step-size by utilizing an approach similar to Theorem~\ref{lem:ode}, assuming that the associated ODE has a \emph{global exponential attractor} $\theta^*$, that is, one can pick $T > 0$ independently of $\theta \in \Theta$ such that $\norm{\flow(T; \theta) - \theta^*} \leq \frac12\norm{\theta - \theta^*}$. Then, letting $N = \lceil T / \deltat \rceil$, we have that $t_{n,N} \in [T, T + \deltat)$ and $\epsilon_{n,N} \in [T\deltat, (T+\deltat)\deltat)$, which indicates that for any $C > 0$ there exists $\bar\deltat > 0$ such that the RHS of \eqref{eq:ode-lemma} is no greater than $C$ for any $\deltat \in (0, \bar \deltat]$. Then we have
\begin{align*}
 \E[ \norm{\theta_{N (k+1)} - \theta^*} ] 
 &\leq \E[ \norm{\flow(N\deltat; \theta_{N k}) - \theta^*} + \norm{\theta_{N (k+1)} - \flow(N\deltat; \theta_{N k})} ] \\
 &\correct{=}{} \E[ \norm{\flow(N\deltat; \theta_{N k}) - \theta^*}] + \E[ \norm{\theta_{N (k+1)} - \flow(N\deltat; \theta_{N k})} ] \\
 &\leq (1/2) \E[\norm{\theta_{N k} - \theta^*}] + C \enspace.
\end{align*}
From the above inequality we obtain an upper bound $\limsup_{k \to \infty} \E[ \norm{\theta_{N k} - \theta^*} ] \leq 2 C$.
We can conclude that $\E[ \norm{\theta_{N \cdot k} - \theta^*} ]$ eventually becomes arbitrarily small for small enough $\deltat$. However, we do not derive the convergence of $\{\theta_{n}\}$ towards $\theta^*$. 

There have been a few applications of ODE-based methods to analyzing deterministic-optimization algorithms. For example,
Yin et al.~\cite{Yin:1995wj,Yin1995ec} have analyzed the ($1$, $\lambda$)-evolution strategy by the ODE method.
However, $\theta$ encoded only the mean vector $\mm$ in those studies; the overall standard deviation $\sigma$ was controlled by a function of the gradient of the objective function.
Hence the algorithm analyzed is significantly different from state-of-the-art evolution strategies that adapt $\mm$ and $\sigma$ simultaneously.
Moreover, since prior studies were based on the asymptotic relation between the parameter sequence and the ODE where they assume the decreasing learning rate or a constant but infinitesimal learning rate, they could not obtain geometric convergence (Theorem~5.2 of \cite{Yin1995ec}).
Recently, a probability-model-based search algorithm that can be regarded as a variation of CMA-ES has been proved globally convergent \cite{Zhou:2018tq}.
This algorithm is more practical than the one considered in \cite{Yin:1995wj,Yin1995ec}, yet the analysis relies on the ODE method with decreasing step-size.
Therefore, the convergence is not geometric.

Gerencs{\'e}r and V{\'a}g{\'o} \cite{gerencser2001mathematics} have developed an approach based on the ODE method to prove the geometric convergence of SPSA for noise free optimization with a constant step-size.
In their work, $\theta_n$ represents the current candidate solution $X_n$ in $\R^d$ and $F_n$ is an estimate of the gradient of function $f$ that is approximated by taking a pair of symmetric perturbations with perturbation strength $c_n$, i.e., $X_n \pm c_n \Delta_n$, where $\Delta_n$ is a random vector.
The geometric convergence of the SPSA with a constant $c_n = c$ on a convex quadratic function has been proven, but the approach cannot be generalized to a context where both $X_n$ and $c_n$ are adapted.
Note also that, to obtain a good practical algorithm, it is crucial to adapt $c_n$, in the same way that it is crucial to use line-search procedures to determine the step-size in gradient-based algorithms.
This adaptive scenario is not covered in \cite{gerencser2001mathematics}.

The fact that the algorithms we are interested in analyzing run with a constant step-size is not the only difficulty that must be overcome to prove the geometric convergence of rank-based search algorithms.
For example, the optimal parameter $\theta^*$ (i.e., the parameter, usually clear for a given context, to which convergence would be optimal) is typically located on the boundary of the domain $\Theta$ and thus not necessarily included in $\Theta$ itself.
In the case of SSA-ES, the optimal parameter is $\theta^* = (x^*,0)$.
Since the standard deviation of a probability distribution should be positive, $\theta^*$ is on the boundary $\partial \Theta$ and is not included in $\Theta$.
Convergence towards points on the boundary that are not necessarily in the parameter set is not standard in ODE analysis.

Another issue is that the optimal parameter $\theta^*$ may not be the unique equilibrium point of the associated ODE \eqref{eq:associated-ode} in the closure $\cl(\Theta)$.
In some cases, we might be able to extend the domain $\Theta$ to $\cl(\Theta)$ by extending $F$ by continuity.
However, for the algorithms we are interested in, we have $F(\theta) = 0$ for some $\theta \in \partial\Theta$.
For example, if $\sigma_n$ is zero in the SSA-ES described in \S\ref{sec:es}, we have $F_n = 0$, implying that $F(\theta) = 0$ for any $\theta \in \partial\Theta$ where $\vv$ is zero.
Therefore, in any \new{$\epsilon$-}neighborhood $\{\theta \in \cl(\Theta); \norm{\theta - \theta^*} < \epsilon\}$ of $\theta^*$, there exist infinitely many equilibrium points of the associated ODE \eqref{eq:associated-ode}.
This prevents us from choosing the most commonly used Lyapunov function $V(\theta) = \norm{\theta - \theta^*}^2$ to prove that the optimal point $\theta^*$ is a global attractor of the ODE \eqref{eq:associated-ode}.
Indeed, the time derivative of $V$ along the ODE solution, i.e., $\nabla V^\T F(\theta)$, is arbitrarily close to zero near equilibrium points.
This violates Lyapunov's stability criterion.

The last issue is that the geometric convergence of the parameter vector $\{\theta_n\}_{n\geq 0}$ using the Euclidean metric is generally not what we want.
For the rank-based optimization algorithm described in~\S\ref{sec:ex}, we typically want to show the geometric convergence of candidate solutions, $\{X_{n,i}\}_{i\in\llbracket 1, \lambda\rrbracket, n \geq 0}$, generated from $\{P_{\theta_n}\}_{n \geq 0}$ at each iteration $n$, or the geometric convergence of the sequence of the mean vectors or any other representative points of $\{P_{\theta_n}\}_{n \geq 0}$.
They may not be directly connected to the geometric convergence of $\{\theta_n\}_{n \geq 0}$ in a Euclidean sense since one can pick arbitrary parametrizations of the probability distribution when designing the algorithm. 

\section{Global Geometric Convergence via ODE Method}
\label{sec:linconv}

We present now our main results to prove the geometric convergence of algorithms in \S\ref{sec:ex}.
We seek a function $\potential: \Theta \to \R_{\geq 0}$ that satisfies conditions A1 and A2 presented in the introduction, which imply the geometric convergence of $\potential(\theta_n)$ to zero for an appropriate choice of the learning rate $\alpha$.
The function $\potential$ will typically upper-bound those interesting quantities the geometric convergence we want to prove---for example, if our objective is to show the geometric convergence of the expected distance between the candidate solution $X_n \sim P_{\theta_n}$ and the optimal point $x^*$ of $f: \R^d \to \R$, $\potential(\theta)$ should be chosen such that $\E_{x \sim P_{\theta}}[\norm{x - x^*}] \leq C \potential(\theta)$ for some $C > 0$.
The first condition \Aone{} that should be satisfied by $\potential$ is that $\potential$ decreases along the trajectories of the solutions of the ODE \eqref{eq:associated-ode}, like the Lyapunov functions used to prove stability of equilibrium of ODEs. This decrease should however be \emph{geometric}. 
The second condition \Atwo{} is to control the deviation between the stochastic trajectory and the solution of the ODE.

After showing that \Aone\ and \Atwo\ imply the geometric convergence of $\potential(\theta_n)$ in Theorem~\ref{thm:lc-igo}, we will derive a corollary on the first-hitting-time bound and discuss the choice of $\potential$ for the two example algorithms that were presented in \S\ref{sec:ex}.

\subsection{Theorem: Sufficient Conditions for Geometric Convergence}

In the following theorem, we prove that the conditions \Aone\ and \Atwo\ imply the global geometric convergence of the expectation of $\potential(\theta_n)$ for a small enough $\alpha$ belonging to a set $\Lambda$ characterized below.
\begin{theorem}
\label{thm:lc-igo}
Let $\{\theta_n\}_{n \geq 0}$ be a sequence defined by \eqref{eq:algo} with $\theta_0$ given deterministically.
Let $\flow: \R_{\geq 0} \times \Theta \to \Theta$ be an absolutely continuous function satisfying \eqref{eq:integral-form}\del{ and $\flow(0; \theta) = \theta$}.
Suppose that there is a function $\potential: \Theta \to \R_{\geq 0}$ satisfying the following two conditions:
\begin{description}
\item[\Aone] There exists $\del{\exists} \DeltaAone: \R_{\geq 0} \to \R_{\geq 0}$ nonincreasing such that $\DeltaAone(t) \downarrow 0$ as $t \uparrow \infty$, and for any $\theta \in \Theta$ and any $t \in \R_{\geq 0}$
\begin{equation}\label{eq:Aone}
\potential(\flow(t; \theta)) \leq \DeltaAone(t) \potential(\theta) \enspace;
\end{equation}
\item[\Atwo] There exists $\del{\exists} \DeltaAtwo: \R_{+} \times \R_{+} \to \R_{+}$ nondecreasing \wrt each argument such that $\DeltaAtwo(\deltat, T) \downarrow 0$ as $\deltat \downarrow 0$ for any fixed $T$, and for any $N \in \N_{+}$ and $\theta_{0} \in \Theta$
\begin{equation}
\E_0[\potential(\theta_{N})] \leq \potential(\flow(N \deltat; \theta_{0}))  + \DeltaAtwo(\deltat, N\deltat) \potential(\theta_{0}) \enspace .
\end{equation}
\end{description}%
Then, there exists an $\bar{\deltat} > 0$ such that, for any $\alpha \in (0, \bar{\deltat}]$, there exist $\CR < 1$ and $C_{\deltat} > 0$ such that for all $ \theta_{0} \in \Theta$ and for all $n \in \N_{+}$, the stochastic algorithm \eqref{eq:algo} satisfies the inequality
\begin{equation}
\E[\potential(\theta_{n})] \leq \CR^{n} C_{\deltat} \potential(\theta_{0})\enspace.
\label{eq:lc-anytime}
\end{equation}
Moreover, 
\begin{equation}
\limsup_{n\to\infty} \frac{1}{n} \ln \E[\potential(\theta_{n})]\leq \ln \CR\enspace.
\label{eq:lc-limit}
\end{equation}
\end{theorem}

\begin{proof}

  We will first refine the statement as follows: 
  Let $\CR = \inf_{N \geq 1}(\DeltaAone(N\deltat) + \DeltaAtwo(\deltat, N \deltat))^{1/N}$ and define $\deltatset = \{ \deltat > 0: \CR < 1\}$.
  Then, $\deltatset$ is nonempty and there exists $\bar{\deltat} > 0$ such that $(0, \deltat] \correct{\subseteq} \deltatset$. 
  For any $\deltat \in \deltatset$ there exists at least one finite integer $N \in \N_{+}$ such that $(\DeltaAone(N\deltat) + \DeltaAtwo(\deltat, N\deltat ))^{1/N} = \CR<1$.
  Let the smallest such $N$ be denoted by $\optN$, and let $\remCR = 1$ if $\optN = 1$ and $\remCR = \max_{N \in \llbracket 1, \optN-1\rrbracket}(\DeltaAone(N\deltat) + \DeltaAtwo(\deltat, N \deltat))$ otherwise.
  Then for any $\deltat \in \deltatset$, \eqref{eq:lc-anytime} and \eqref{eq:lc-limit} hold with $C_{\deltat} = \remCR/ \CR^{\optN-1}$.

First, we show that $\deltatset$ is nonempty.
According to \Aone, we can take $T$ such that $\DeltaAone(T) < 1/2$.
Then, since according to \Atwo, $\DeltaAtwo(T / N, T) \downarrow 0$ as $N \uparrow \infty$, we can choose $N$ such that $\DeltaAtwo(T / N, T) < 1/2$ for the fixed $T$.
For such an $N$ and $\deltat = T / N$, we have $\DeltaAone(N\deltat) + \DeltaAtwo(\deltat, N\deltat) < 1$.
As $\CR \leq (\DeltaAone(N\deltat) + \DeltaAtwo(\deltat, N\deltat))^{1/N} < 1$ for such an $\deltat$, it follows that $\deltat \in \deltatset$.
Therefore, $\deltatset$ is nonempty.

Next, we show that there exists $\bar{\deltat} > 0$ such that $(0, \bar{\deltat}] \correct{\subseteq} \deltatset$.
Let $T$ be such that $\DeltaAone(T) < 1/2$, and let $\bar N$ be the smallest $N \in \N$ such that $\DeltaAtwo( 2 T / N, 2 T) < 1/2$ and $\bar \deltat = 2 T / \bar N$.
Then, for any $\deltat \in [\bar \deltat / 2, \bar \deltat]$, we have that $\deltat \bar N \geq T$ and that $\DeltaAone(\deltat \bar N) \leq \DeltaAone(T) < 1/2$, and $\DeltaAtwo( \deltat, \deltat \bar N) \leq \DeltaAtwo( 2 T / \bar N, 2 T) < 1/2$.
As $\CR \leq [\DeltaAone(\deltat \bar N) + \DeltaAtwo( \deltat, \deltat \bar N)]^{1 / \bar N} < 1$, we have $[\bar \deltat / 2, \bar \deltat] \subseteq \deltatset$.
If we let $\bar N = 2^{i+1} T / \bar{\deltat}$, the same argument holds for $\deltat \in [\bar{\deltat}/2^{i+1}, \bar{\deltat}/2^{i}]$.
Therefore, $[\bar{\deltat}/2^{i+1}, \bar{\deltat}/2^{i}] \subseteq \deltatset$ for all $i \in \N$.
Hence, $\cup_{n=0}^{\infty} [\bar \deltat / 2^{i+1}, \bar \deltat / 2^{i}] = (0, \bar \deltat] \subseteq \deltatset$.

Now we prove that for any $\deltat \in \deltatset$ there exists an $N \in \N_{+}$ such that $\CR = (\DeltaAone(N\deltat) + \DeltaAtwo(\deltat, N \deltat))^{1/N}$.
If we assume that there does not exist an $N$ that minimizes $(\DeltaAone(N\deltat) + \DeltaAtwo(\deltat, N\deltat ))^{1/N}$, then $\liminf_{N}(\DeltaAone(N\deltat) + \DeltaAtwo(\deltat, N\deltat ))^{1/N} = \CR < 1$ must hold.
However, since $\DeltaAone(N\deltat)$ is nonnegative and $\DeltaAtwo(\deltat, N\deltat) \geq \DeltaAtwo(\deltat, \deltat) > 0$, we have a contradiction:  $\liminf_{N}(\DeltaAone(N\deltat) + \DeltaAtwo(\deltat, N\deltat ))^{1/N} \geq \liminf_{N} \DeltaAtwo(\deltat, \deltat)^{1/N} = 1$.
Hence, there exists at least one $N \in \N_{+}$ that minimizes $(\DeltaAone(N\deltat) + \DeltaAtwo(\deltat, N \deltat))^{1/N}$.

Since an algorithm of the form~\eqref{eq:algo} is time-homogeneous, the assumption \Atwo{} implies $\E_n[\potential(\theta_{n+N})] \leq \potential(\flow(N \deltat; \theta_{n}))  + \DeltaAtwo(\deltat, N\deltat) \potential(\theta_{n})$ for any $n \in \N$.
From \Aone\ and \Atwo, we have for $n \in \N$
\begin{align}
\label{eq:n-step-dec}
\E_{n}[\potential(\theta_{n+\optN})] 
  &\leq \potential(\flow(\optN \deltat; \theta_{n}))  + \DeltaAtwo(\deltat, \optN\deltat) \potential(\theta_{n})  \\
  &\leq \DeltaAone(\optN\deltat) \potential(\theta_{n})  + \DeltaAtwo(\deltat, \optN\deltat) \potential(\theta_{n}) 
  = \CR^{\optN} \potential(\theta_{n}) \enspace.\notag
\end{align}
Then, for any $n \in \N$ and $k \in \N$,
\begin{equation}
  \begin{split}
\E[\potential(\theta_{n+\optN k})] 
&= \E[ \E_{n+ \optN(k-1)}[\potential(\theta_{n+\optN k})] ]
\leq \CR^{\optN} \E[ \potential(\theta_{n+\optN(k-1)})] \\
&= \CR^{\optN} \E[ \E_{n+ \optN(k-2)}[\potential(\theta_{n+\optN (k-1)})] ]
\leq \CR^{2 \optN} \E[ \potential(\theta_{n+\optN(k-2)})] \\
&\leq \cdots \leq \CR^{k \optN} \E[ \potential(\theta_{n})] \enspace.
\end{split}\label{eq:n-step-dec-exp}
\end{equation}
If $\optN = 1$, this immediately leads to \eqref{eq:lc-anytime}, where $C_\deltat = \remCR /\CR^{\optN-1} = 1$.
On the other hand, if $\optN > 1$, for $n \in \llbracket 1, \optN-1\rrbracket$, we have $\E[\potential(\theta_{n})] \leq \remCR \potential(\theta_{0})$, since $\remCR \geq (\DeltaAone(n\deltat) + \DeltaAtwo(\deltat, n \deltat))$.
Combining this with \eqref{eq:n-step-dec-exp}, we obtain
\begin{align*}
  \E[\potential(\theta_{n+\optN k})] \leq \CR^{k \optN} \remCR \potential(\theta_{0}) &= \CR^{n + k \optN} (\remCR / \CR^{n})\potential(\theta_{0})
  \leq \CR^{n + k \optN} (\remCR / \CR^{\optN-1}) \potential(\theta_{0})
\enspace.
\label{eq:lc-igo:exp}
\end{align*}
Rewriting $n + \optN k$ as $n$, we have $\E[\potential(\theta_{n})] \leq \CR^{n} C_\deltat \potential(\theta_{0})$ for any $n \in \N$, where $C_\deltat = \remCR / \CR^{\optN-1} > 0$.
This completes the proof of \eqref{eq:lc-anytime}.

Taking the natural logarithm, dividing by $n$, and taking the limit of the supremum, we finally have $\limsup_{n} \frac{1}{n} \ln(\E[\potential(\theta_{n})] / \potential(\theta_{0}))\leq \ln \CR$.
This completes the proof of \eqref{eq:lc-limit}.
\end{proof}

Suppose that we can choose $\potential$ such that \Aone{} and \Atwo{} hold.
Then, an algorithm of the form~\eqref{eq:algo} exhibits geometric convergence of $\{\potential(\theta_n)\}_n$ with a rate of convergence upper-bounded by $\ln \CR$, if the learning rate $\deltat$ is sufficiently small.
The upper bound of the rate of convergence, $\ln \CR$, depends on $\deltat$.
As $\deltat$ becomes smaller, the sequences $\{\potential(\theta_{n+N})\}_n$ and $\{\potential(\flow(N \deltat, \theta_{n}))\}_n$ become closer for a fixed $N$ because of \Atwo, but $\potential(\flow(N \deltat, \theta_{n})) / \potential(\theta_{n})$ becomes larger because of \Aone.
This may lead to a larger $\ln \CR$.
In a similar way, we lower-bound the expected decrease of $\potential(\theta_{n})$ in the next theorem: 
\begin{theorem}
  \label{thm:lc-igo-lower}
  Let $\{\theta_n\}_{n \geq 0}$ and $\flow$ be as in Theorem~\ref{thm:lc-igo}. Suppose that there is a function $\potential: \Theta \to \R_{\geq 0}$ satisfying the following two conditions:
  \begin{description}
  \item[\Bone] There exists $\DeltaBone: \R_{\geq 0} \to \R_{\geq 0}$ nonincreasing such that $\DeltaBone(t) \downarrow 0$ as $t \uparrow \infty$ and $\DeltaBone(t) \uparrow c$ as $t \downarrow 0$ for some $c \in (0, 1]$, and for any $\theta \in \Theta$ and any $t \in \R_{\geq 0}$
    \begin{equation}
      \potential(\flow(t; \theta)) \geq \DeltaBone(t) \potential(\theta) \enspace;
    \end{equation}
  \item[\Btwo] There exists $\DeltaBtwo: \R_{+} \times \R_{+} \to \R_{+}$ nondecreasing \wrt each argument such that $\DeltaBtwo(\deltat, T) \downarrow 0$ as $\deltat \downarrow 0$ for any fixed $T$, and for any $N \in \N_{+}$ and $\theta_{0} \in \Theta$
    \begin{equation}
      \E_0[\potential(\theta_{N})] \geq \potential(\flow(N \deltat; \theta_{0}))  - \DeltaBtwo(\deltat, N\deltat) \potential(\theta_{0}) \enspace .
    \end{equation}
  \end{description}%

    Then, there exists an $\bar{\deltat} > 0$ such that, for any $\alpha \in (0, \bar{\deltat}]$, there exist $\CR > 0$ and $C_{\deltat} > 0$ such that for all $ \theta_{0} \in \Theta$ and for all $n \in \N_{+}$, the stochastic algorithm \eqref{eq:algo} satisfies the inequality
    \begin{equation}
      \E[\potential(\theta_{n})] \geq \CR^{n} C_{\deltat} \potential(\theta_{0})\enspace.
      \label{eq:lc-anytime-lower}
    \end{equation}
  Moreover, 
  \begin{equation}
    \liminf_{n\to\infty} \frac{1}{n} \ln \E[\potential(\theta_{n})]\geq \ln \CR\enspace.
    \label{eq:lc-limit-lower}
  \end{equation}
\end{theorem}

\begin{proof}
  We prove the refined statement:
    For a given $\deltat$, let $J_{\deltat}$ be the maximum integer such that $\DeltaBone(N\deltat) > \DeltaBtwo(\deltat, N \deltat)$ for all $N \in \llbracket 1, J_{\deltat}\rrbracket$.
    Let $\CR = \max_{N \in \llbracket 1, J_{\deltat}\rrbracket}(\DeltaBone(N\deltat) - \DeltaBtwo(\deltat, N \deltat))^{1/N}$ and $\optN$ be the smallest integer realizing $\CR$.
    Define $\deltatset = \{ \deltat > 0: \CR > 0\}$.
    Then, $\deltatset$ is nonempty, and there exists $\bar{\deltat} > 0$ such that $(0, \bar\deltat] \correct{\subseteq} \deltatset$. 
    For $\deltat \in \deltatset$, let $\remCR = 1$ if $\optN = 1$ and $\remCR = \min_{N \in \llbracket 1, \optN-1\rrbracket}(\DeltaBone(N\deltat) - \DeltaBtwo(\deltat, N \deltat)) > 0$ otherwise.
    Then for any $\deltat \in \deltatset$, \eqref{eq:lc-anytime-lower} and \eqref{eq:lc-limit-lower} hold with $C_{\deltat} = \remCR > 0$. 
  
  The proof is similar to that for Theorem~\ref{thm:lc-igo}.
  Since $\DeltaBone(N\deltat)$ is decreasing to zero and $\DeltaBtwo(\deltat, N \deltat) \geq \DeltaBtwo(\deltat, \deltat) > 0$, it is easy to see that $J_{\deltat}$ exists as a finite number for any $\deltat > 0$.
  Take $T$ such that $\DeltaBone(T) > c/2$ and choose $N$ such that $\DeltaBtwo(T/N, T) < c/2$.
  Then, for $\alpha = T/N$, we have $\CR \geq \DeltaBone(N\deltat) - \DeltaBtwo(\deltat, N \deltat) > 0$.
  Therefore, $\deltatset$ is nonempty.

    Next, we show that there exists $\bar{\deltat} > 0$ such that $(0, \bar{\deltat}] \correct{\subseteq} \deltatset$.
    Let $T$ be such that $\DeltaBone(T) > c/2$ and let $\bar N$ be the smallest $N \in \N$ such that $\DeltaBtwo( T / N, T) < c/2$ and $\bar \deltat = T / \bar N$.
    Then, for any $\deltat \in (0, \bar \deltat]$, we have that $\deltat \bar N \leq T$ and that $\DeltaBone(\deltat \bar N) \geq \DeltaBone(T) > c/2$ and $\DeltaBtwo( \deltat, \deltat \bar N) \leq \DeltaBtwo( T / \bar{N}, T) < c/2$.
    Since $\CR \geq [\DeltaBone(\deltat \bar N) - \DeltaBtwo( \deltat, \deltat \bar N)]^{1 / \bar N} > 0$, we have $(0, \bar \deltat] \subseteq \deltatset$.

  Instead of \eqref{eq:n-step-dec}, we have from \Bone{} and \Btwo,
  \begin{align*}
    \E_{n}[\potential(\theta_{n+\optN})] 
    &\geq \DeltaBone(\optN\deltat) \potential(\theta_{n}) - \DeltaBtwo(\deltat, \optN\deltat) \potential(\theta_{n}) 
      = \CR^{\optN} \potential(\theta_{n}) \enspace.
  \end{align*}
  Then, for any $n \in \N$ and $k \in \N$,
  \begin{align*}
    \E[\potential(\theta_{n+\optN k})] 
    &\geq \cdots \geq \CR^{k \optN} \E[ \potential(\theta_{n})] \enspace.
  \end{align*}
  This leads to \eqref{eq:lc-anytime-lower} if $\optN = 1$ with $C_\deltat = \remCR = 1$. 
    For $\optN > 1$ and $n \in \llbracket 0 , \optN-1\rrbracket$, since $\remCR \leq \DeltaBone(n\deltat) - \DeltaBtwo(\deltat, n \deltat)$, we have $\E[\potential(\theta_{n})] \geq \remCR \potential(\theta_{0})$. Combining them, we obtain
  \begin{equation*}
    \E[\potential(\theta_{n+\optN k})] \geq \CR^{k \optN} \remCR \potential(\theta_{0}) \geq \CR^{n + k \optN} \remCR\potential(\theta_{0}) 
    \enspace,
    \label{eq:lc-igo:exp}
  \end{equation*}
  where we used $\CR \leq 1$. 
  Rewriting $n + \optN k$ as $n$, we have $\E[\potential(\theta_{n})] \geq \CR^{n} \remCR \potential(\theta_{0})$ for any $n \in \N$. 
  Since $\DeltaBone(\deltat n) - \DeltaBtwo( \deltat, \deltat n)$ is nonincreasing w.r.t.\ $n$, we have $\remCR \geq \DeltaBone(\deltat \optN) - \DeltaBtwo( \deltat, \deltat \optN) = \CR^{\optN} > 0$.
    Hence, $C_\deltat = \remCR > 0$.
  This completes the proof of \eqref{eq:lc-anytime-lower}.

  Taking the natural logarithm, dividing by $n$, and taking the limit of the infimum, we finally have $\liminf_{n} \frac{1}{n} \ln(\E[\potential(\theta_{n})] / \potential(\theta_{0}))\geq \ln \CR$.
  This completes the proof of \eqref{eq:lc-limit-lower}.
\end{proof}

\subsection{Consequences in Optimization: Geometric Convergence and Hitting-Time Bound}

Consider a rank-based stochastic algorithm of the form~\eqref{eq:algo-rank} minimizing a deterministic function $f: \X \to \R$.
Let $d: \X \times \X \to \R_{\geq0}$ be a distance function on $\X$ and $x^* \in \X$ be the well-defined optimum of the objective function, $x^* = \argmin_{x \in \X} f(x)$.
As a consequence of the geometric convergence of $\{\potential(\theta_{n})\}_n$, we can derive from \eqref{eq:lc-anytime} a bound on the rate of convergence of the sequence of random vectors $\{X_{n} \in \X\}_{n\geq0}$ drawn from $P_{\theta_n}$ provided the expected distance between $X_n$ and the optimum is upper-bounded by a constant times $\potential(\theta_n)$, as stated in the following corollary:
\begin{corollary}
If $\potential(\theta_n)$ satisfies the assumptions of Theorem~\ref{thm:lc-igo} and $\E_n[d(X_n, x^*)] \leq C \potential(\theta_n)$ for some constant $C > 0$, we obtain the any-time bound of the expected distance $\E[d(X_n, x^*)]$:
\begin{equation}
	\E[d(X_n, x^*)] \leq C C_{\deltat} \CR^{n}  \potential(\theta_{0})
	\label{eq:lc-dbound}
\end{equation}
for any $\deltat \in (0, \bar\deltat]$ and any $\theta_0 \in \Theta$, 
where $\bar\deltat$, $\CR$, $C_{\deltat}$ are as defined in Theorem~\ref{thm:lc-igo} and independent of $\theta_0$.
Moreover, 
\begin{equation}
	\limsup_{n} \frac1n \ln\E[d(X_n, x^*)] \leq \ln \CR \enspace.\label{eq:lc-CR}
\end{equation}
\end{corollary}%
That is, the asymptotic bound on the expected distance to the optimum is $\bigO(\CR^{n})$.
Similarly, if we can find $\potential$ satisfying the assumptions of Theorem~\ref{thm:lc-igo-lower} and $\E_n[d(X_n, x^*)] \geq C \potential(\theta_n)$ for some $C$, we can deduce a lower bound on the expected distance.
The asymptotic bound is then $\Omega(\CR^n)$, where $\CR$ is different from the one for the upper bound.
Combining them, we obtain the geometric decrease of the expected distance.

By contrast, in machine learning, stochastic gradient methods constitute an important class of algorithms, the convergence of which is analyzed with ODE methods or stochastic approximation techniques. 
In machine-learning problems, the bound on the generalization error achieved by those methods scales like $\bigO(n^{-\beta})$ for some constant $\beta > 0$ \cite{Bach:2011vy} and matches the asymptotic lower bound \cite{Agarwal2012}, where $\beta = 1/2$ for convex functions with $\ell_q$-Lipschitz continuity with $q \geq 1$, and $\beta = 1$ if in addition the function is $\ell_2$-strong convex.

The difference between these bounds shows that our approach using Theorem~\ref{thm:lc-igo} differs significantly from prior work based on conventional ODE methods. 
Note, however, that Le Roux et al.~\cite{Roux2012nips} analyzed a variant of the stochastic gradient method with a constant step-size in a finite batch setting and showed a training error bound scaling like $\bigO(\gamma^n)$ for the given batch, assuming strong convexity of the sample average and smoothness of each summand.
Their proof relies on finding a Lyapunov function $V:\Theta \to \R_{\geq 0}$ of a parameter $\theta$ and showing $\E_n[V(\theta_{n+1})] = c \cdot V(\theta_n)$ for some constant $c \in (0, 1)$ independent of the iteration count $n$.
The idea is similar to the one in our paper, but they do not go through the mean ODE.
In our setting, unlike theirs, it is rather non-trivial to find a function $\potential$ such that $\potential(\theta_{n})$ decreases at every iteration. This is due to the comparison-based and adaptive nature of the algorithms.
Going through the mean ODE allows us to show this decrease in $\optN$ iterations easily (see the proof of Theorem~\ref{thm:lc-igo}); we expect it would be tedious to show it by hand.

Another consequence of Theorem~\ref{thm:lc-igo} concerns the first-hitting-time bound, a commonly accepted runtime measure for algorithms in discrete domains.
Given $\epsilon > 0$, we define the first hitting time of $\potential(\theta_n)$ to $[0, \epsilon]$ by
\begin{equation}
\tau = \inf_{n} \{n \in \N \mid \potential(\theta_n) \leq \epsilon\} \enspace.
\label{eq:fht}
\end{equation}
If $\E_n[d(X_n, x^*)] \leq C \potential(\theta_n)$, the first hitting time of $\E_n[d(X_n, x^*)]$ to $[0, C \epsilon]$ is immediately upper-bounded. 
Combining Theorem~\ref{thm:lc-igo} and drift analysis \cite{Hajek1982}, we can deduce an upper bound of the first hitting time.

\begin{corollary}\label{cor:tau}
  Let $\bar\deltat$ and $\CR$ be as defined in Theorem~\ref{thm:lc-igo}. Then, for any $\deltat \in (0, \bar\deltat]$, there exist constants $\beta_{0}$, $\beta_{1}$, $\bar\beta_{0}$, $\bar\beta_{1} \in \R_+$ such that for any $\theta_0 \in \Theta$, $\epsilon \in (0, \potential(\theta_0)]$, and $\delta \in (0, 1]$, 
  \begin{align}
    &\E[ \tau ] < \frac{\beta_{1}}{\ln(1/\CR)} \left( \ln\left(\frac{\potential(\theta_0)}{\epsilon}\right) + \beta_{0} \right) \enspace,\label{eq:exptau}
    \\
    &\Pr\left[ \tau > \frac{\bar\beta_{1}}{\ln(1/\CR)} \left( \ln\left(\frac{\potential(\theta_0)}{\epsilon \delta}\right) + \bar\beta_{0} \right) \right] \correct{\leq}{} \delta  \enspace.\label{eq:protau}
  \end{align}
\end{corollary}
\begin{proof}
  Let $C_{\deltat}$ be the constant defined in Theorem~\ref{thm:lc-igo}.
  We prove the statement with
  $\beta_1 = (2/\ln(2)) \ln(2 C_{\deltat}/\CR)$,
  $\beta_0 = 1 + \ln(2) / 2$,
  $\bar\beta_1 = \ln(2 C_{\deltat}/\CR) / \ln(2)$, and
  $\bar\beta_0 = \ln(2)$. 
  Let $N = \lceil \ln(2 C_{\deltat}) / \ln(1/\CR) \rceil$ so that $\rho := C_{\deltat} \CR^{N} \correct{\leq}{} 1/2$.   
  Let $Y_k = \ln \potential(\theta_{N k})$ and $\mathcal{E}_k$ be its natural filtration.
  Then, it is easy to see that $\E[\exp(Y_{k+1} - Y_{k}) \mid \mathcal{E}_k] \leq \rho$.
  This satisfies the prerequisites of Theorem 2.3 of \cite{Hajek1982}.
  Therefore, for the first hitting time $\tau_Y$ of $Y_k$ to $(-\infty, \ln(\epsilon)]$,
  we obtain $\E[s^{\tau_Y}] \leq 1 + \frac{s-1}{1 - \rho s}\frac{\potential(\theta_0)}{\epsilon}$ for any $s \in (1, \rho^{-1})$ and
$\Pr[\tau_Y > k] \leq \frac{\potential(\theta_0)}{\epsilon} \rho^k$.
Since $Y_k \leq \ln(\epsilon)$ implies $\potential(\theta_{Nk}) \leq \epsilon$, we have
$\tau \leq N \tau_Y$, then $(s^{1/N})^{\E[\tau]} \leq \E[s^{\tau / N}] \leq \E[s^{\tau_Y}]$ and $\Pr[\tau > \bar\tau] \leq \Pr[\tau_Y > \lfloor \bar\tau / N \rfloor ]$.
  Let $s = \rho^{-1/2}$.
  By using the inequality $\ln(1 + x) \leq 1 + \ln(x)$ for any $x \geq 1$, we obtain
  $\E[\tau] \leq \frac{N \ln(\E[s^{\tau_Y}])}{\ln(s)} \leq \frac{2N}{\ln(1/\rho)} \left( 1 + \frac{\ln(1/\rho)}{2} + \ln\left(\frac{\potential(\theta_0)}{\epsilon} \right)\right)$.
  The right-most side (RMS) is further upper-bounded by the RHS of \eqref{eq:exptau} by using $N \correct{\leq}{} \ln(2 C_{\deltat}/\CR) / \ln(1/\CR)$ and $\rho \correct{\leq}{} 1/2$.
  On the other hand, we have $\Pr[\tau > \bar\tau] \leq (\potential(\theta_0) / \epsilon) \rho^{\lfloor \bar\tau / N \rfloor} \correct{\leq} (\potential(\theta_0) / \epsilon) (1/2)^{\bar\tau / N - 1}$.
  Solving $(\potential(\theta_0) / \epsilon) (1/2)^{\bar\tau / N - 1} \leq \delta$ for $\bar\tau$,
  we obtain $\bar\tau \geq N + N \ln(\potential(\theta_0) / (\epsilon\delta)) / \ln(2)$.
  Since $N \leq \ln(2 C_{\deltat}/\CR) / \ln(1/\CR) = \ln(2) \bar\beta_1 / \ln(1/\CR)$,
  we obtain $\Pr[\tau > \bar\tau ] \leq \delta$ for $\bar\tau = \frac{\bar\beta_{1}}{\ln(1/\CR)} \left( \ln\left(\frac{\potential(\theta_0)}{\epsilon \delta}\right) + \bar\beta_{0}\right)$. This proves \eqref{eq:protau}.
\end{proof}

\subsection{Connection with Lyapunov Functions and the Choice of $\potential$}\label{sec:comment}

We explain in this section the link between the function $\potential$ and a Lyapunov function to analyze the stability of equilibrium points of ODEs.
We then discuss how practically to choose $\potential$ and illustrate some possible choices with two examples. 

The function $\potential$ can be seen as a Lyapunov function for the stability analysis of an ODE at a stationary point.
Recall that a Lyapunov function $V: \R^\dimtheta \to \R_{\geq 0}$ used to investigate the asymptotic stability of an ODE $\frac{\rmd \theta}{\rmd t} = F(\theta)$ at a stationary point $\theta^*$ is required to be continuous and nonnegative in a neighborhood $U \subseteq \Theta$ of $\theta^*$.
In addition $V(\theta) = 0$ if and only if $\theta = \theta^*$, and there exists a time derivative $\dot V(\theta) = \new{\frac{\rmd V(\theta(t))}{\rmd t}  =} \nabla V(\theta)^\T F(\theta)$ along the solution which should be negative for $\theta \in U\setminus\{\theta^*\}$ ($V$ decreases along the trajectory $t \mapsto \theta(t)$).
Under those conditions, one can conclude that the solution of the ODE converges towards $\theta^*$ starting from any $\theta \in U$.
A first difference between $\potential$ and a Lyapunov function $V$ is that, since $\theta^*$ is typically on the boundary of $\Theta$ and $\theta^* \notin \Theta$, $\potential(\theta^*)$ or the limit $\lim_{\theta \to \theta^*} \potential(\theta)$ will not be well-defined.
As an example, consider $\potential(\theta) = \norm{\theta}\max(1, \abs{\theta_1}/\abs{\theta_2})$, where $\theta = (\theta_1, \theta_2) \in \R^2$ and the limit $\lim_{\theta \to \theta^*} \potential(\theta)$ is not well-defined.
In this case the geometric convergence of $\potential(\theta) \to 0$ not only implies the convergence of $\theta \to \theta^* = (0, 0)$, but also tells us how it approaches $\theta^*$, i.e., $\theta_2$ cannot converge faster than $\theta_1$ does.
Another difference with the typical Lyapunov functions used to prove the asymptotic stability of solutions of an ODE is that we require more than the negativity of the time derivative of $\potential$ along the solution of the ODE. The conditions are similar to those for the exponential stability analysis.
Also, the zeroes of $\potential$ are not necessarily unique, as we are interested in the geometric convergence of $\potential(\theta_n)$ itself, not the convergence of $\theta_n$ in $\Theta$.
Convergence towards a unique attractor $\stheta$ is achieved if we add the condition that $\potential(\theta) \to 0$ happens only when $\theta \to \stheta$.
Otherwise, the parameter will converge to a subset of the closure of $\Theta$.
Despite these differences, the function $\potential$ is qualitatively very similar to a Lyapunov function.

Condition \Aone{} (or \Athr{} in Theorem~\ref{thm:llc-igo} in the following section) is a condition on the speed of the convergence of $\potential$ along the solution of the ODE \eqref{eq:associated-ode}.
The requirement of \eqref{eq:Aone}
is stronger than assuming the global (or local) convergence of $\potential(\flow(t; \theta))$ to zero, since the function $\DeltaAone$ must be chosen independently of the initial parameter $\theta$.
The condition \eqref{eq:Aone} typically holds in the context we are interested in, because we often have geometric (i.e., exponential) convergence, that is $\potential(\flow(t; \theta)) \leq \potential(\theta) e^{- C t}$ for some $C > 0$, when the algorithm converges geometrically.
However, \eqref{eq:Aone} is slightly less restrictive or weaker than exponential convergence, which allows us to use a loose upper bound $\DeltaAone(t)$ for $\potential(\flow(t; \theta)) / \potential(\theta)$ such as $\DeltaAone(t) = 2 / (t + 1)$.
This is helpful when it is not trivial to derive a tight upper bound.

Assumption \Atwo{} (or \Afou) is a condition on the deviation of the stochastic process $\potential(\theta_{n+N})$ from the ODE solution $\potential(\flow(N \deltat; \theta_n))$, the counterpart of the conclusion of Theorem~\ref{lem:ode}.
It differs from Theorem~\ref{lem:ode} in that we compute the difference in $\potential$ instead of the difference in the Euclidean sense $\norm{\theta_{n+N} - \flow(N\deltat; \theta_n)}$, and that the upper bound is proportional to $\potential(\theta_n)$. 
The first point is because we want to prove the geometric convergence of $\potential(\theta_n)$ rather than the geometric convergence of the parameter sequence in the Euclidean sense $\norm{\theta_n - \theta^*}$.
The second point is because algorithms following \eqref{eq:algo-rank} typically adapt the variance of samples. The parameter $\theta$ may encode the variance of the distribution $P_{\theta}$. The variance of $F_n$ is controlled by $\theta$.
When $\theta_n$ is approaching $\theta^*$ where the variance parameter is zero, the variance of $F_n$, and hence the variance of $\potential(\theta_n)$, approaches zero, while $\potential(\theta_n) \to \potential(\theta^*) = 0$.
If we have an upper bound of the form $\E_{n}[\potential(\theta_{n+N})] \leq \potential(\flow(N \deltat; \theta_{n})) + \DeltaAtwo(\deltat, N\deltat)$, we cannot conclude that $\potential(\theta_n)$ converges geometrically.
However, we may obtain the convergence of $\potential(\theta_n)$ in a compact set including zero, as discussed after Theorem~\ref{lem:ode}.
Therefore, \Atwo{} is essential for deriving the geometric convergence of $\potential(\theta_n)$.

\section{Practical Conditions}\label{sec:pra}

In general, we cannot obtain the flow $\flow$ explicitly, but the conditions presented in Theorem~\ref{thm:lc-igo} and in Theorem~\ref{thm:lc-igo-lower} refer to $\flow$.
In this section, we develop practical conditions to verify \Aone{} and \Atwo{} (or \Athr{} and \Afou{} in Theorem~\ref{thm:llc-igo} in \S\ref{sec:llc}) without knowing the explicit solution $\flow$ of the ODE.
Those conditions are presented in Theorem~\ref{lem:prac-one} and Theorem~\ref{lem:ode-igo}.
We then show in \S\ref{sec:appli} how to apply those practical conditions in the step-size adaptive ES described in \S\ref{sec:comment}.

To prove \Aone{} or \Athr{} without knowing $\flow$ explicitly, we typically utilize Lyapunov's argument for exponential stability \cite{Khalil:2002wj}.
Assume that $\potential$ is differentiable.
If we can show that there exists a constant $C > 0$ such that for any $\theta \in \Theta$, 
\begin{equation}
D_{F(\theta)}\potential(\theta) =  (\nabla \potential(\theta))^\T F(\theta) \leq - C \potential(\theta) \enspace,
\label{eq:lyapunov}
\end{equation}
where $D_{F(\theta)} \potential(\theta)$ denotes the directional derivative of $\potential$ at $\theta$ in the direction of $F(\theta)$, then the flow will satisfy $\potential(\flow(t; \theta)) \leq \potential(\theta) \exp(- Ct)$, which satisfies \Aone{} with $\DeltaAone(t) = \exp(-C t)$.
Then, we do not need to know the explicit solution $\flow$ of the ODE.
An analogous argument holds for the lower bound in \Bone{}.
If $\potential$ has non-differentiable points, we may replace the LHS of the inequality \eqref{eq:lyapunov} with the upper and lower Dini directional derivatives of $\potential$ at $\theta$ in the direction of $F(\theta)$, defined as
\begin{align*}
  D^+_{F(\theta)} \potential(\theta) &= \limsup_{h \downarrow 0} [\potential(\theta + h F(\theta)) - \potential(\theta)]/h\\
  D^-_{F(\theta)} \potential(\theta) &= \liminf_{h \downarrow 0} [\potential(\theta + h F(\theta)) - \potential(\theta)]/h \enspace.
\end{align*}
This is stated formally in the following proposition, which is an extension of Lyapunov's well-known argument for exponential stability \cite{Khalil:2002wj}; the proof is included in \ref{apdx:prac-one}.
\begin{theorem}\label{lem:prac-one}
If $\potential$ is a locally Lipschitz function that for some $C > 0$
satisfies 
\begin{equation}\label{cond-dini}
D^{+}_{F(\theta)} \potential(\theta) \leq - C \potential(\theta)
\end{equation}
for all $\theta \in \Theta$, then the condition \Aone{} in Theorem~\ref{thm:lc-igo} holds with $\DeltaAone(t) = \exp(-C t)$. 

Similarly, if \eqref{cond-dini} holds for any $\theta \in U = \{\theta \in \Theta; \potential(\theta) < \supU\}$ for some $\supU > 0$, then condition \Athr{} in Theorem~\ref{thm:llc-igo} holds with these $\supU$, $U$, and $\DeltaAthr(t) = \exp(-C t)$.

Moreover, if there exists $C > 0$ such that 
\begin{equation}\label{cond-dini-lower}
D^{-}_{F(\theta)} \potential(\theta) \geq - C \potential(\theta)
\end{equation}
for any $\theta \in \Theta$, then the condition \Bone{} in Theorem~\ref{thm:lc-igo-lower} holds with $\DeltaBone(t) = \exp(-C t)$. 
\end{theorem}
If $\potential$ is differentiable, the constants for the lower and the upper bounds are given by
\begin{align*}
  C_\mathrm{upper} = \sup_{\theta \in \Theta} \nabla (\ln \potential(\theta))^\T F(\theta) \quad \text{and} \quad
  C_\mathrm{lower} = \inf_{\theta \in \Theta} \nabla (\ln \potential(\theta))^\T F(\theta) \enspace.
\end{align*}

To prove \Atwo{} or \Afou{} without knowing $\flow$ explicitly, we may use the following theorem and its corollary, the proofs of which are included in \ref{apdx:lem:ode-igo} and \ref{apdx:cor:ode-igo}.
For a positive definite symmetric matrix $Q$, $\norm{\theta}_{Q}^2 = \theta^\T Q \theta$ denotes the square Mahalanobis norm. 

\begin{theorem}\label{lem:ode-igo}
  Consider an algorithm of the form \eqref{eq:algo} with a deterministic and time-independent step-size $\deltat$.
  Given $\theta_{n} \in \Theta$ and $N \in \N_+$, let $\Theta_{\deltat,N}(\theta_n)$ and $\tilde \Theta_{N\deltat}(\theta_n)$ be subsets of $\Theta$ in which $\{\theta_{n+k}\}_{k=0}^{N}$ and $\flow(t; \theta_n)$ for $t \in [0, N\deltat]$, respectively, almost surely stay.
  Let $U \subseteq \Theta$. We assume that for any $\theta_n \in U$ \correct{and for any $N \in \N_{+}$},
  \begin{itemize}
  \item[\asm{P1}] there exists a positive definite symmetric matrix $Q_{\deltat,N}(\theta_n)$ satisfying the inequality $\sup_{\norm{h} = 1} \abs{D^{+}_{h} \potential(\theta)}\norm{h}_{ Q_{\deltat,N}(\theta_n)^{-1} } \leq 1$ for all $\theta$ in the convex hull of $\Theta_{\deltat,N}(\theta_n) \cup \tilde \Theta_{N\deltat}(\theta_n)$, where this condition
can be replaced with $\norm{\nabla\potential(\theta)}_{ Q_{\deltat,N}(\theta_n)^{-1} } \leq 1$ if $\potential$ is differentiable;

    \item[\asm{P2}] 
        there exists a function $L_{\deltat,N}(\theta, \theta_n) > 0$ such that $\norm{F(\theta) - F(\theta')}_{Q_{\deltat,N}(\theta_n)} \leq L_{\deltat,N}(\theta, \theta_n) \norm{\theta - \theta'}_{Q_{\deltat,N}(\theta_n)}$ for any $\theta \in \Theta_{\deltat,N}(\theta_n)$ and $\theta'  \in \tilde \Theta_{N\deltat}(\theta_n)$, and
        there exists a function $\Delta_{L}$ that is independent of $\theta_n \in U$, nondecreasing \wrt both arguments, and satisfies $\E_n[ \prod_{k=0}^{N-1} (1 + \deltat L_{\deltat,N}(\theta_{n+k}, \theta_n))^2 ]^{1/2} \leq \Delta_{L}(\deltat, N\deltat)$;

  \item[\asm{P3}] there exists a function $\Delta_{\tilde{L}}$ that is independent of $\theta_n \in U$, nondecreasing \wrt both arguments, and satisfies
    $\norm{F(\theta) - F(\theta')}_{Q_{\deltat,N}(\theta_n)} \leq \Delta_{\tilde{L}}(\deltat, N\deltat) \norm{\theta - \theta'}_{Q_{\deltat,N}(\theta_n)}$ for any $\theta, \theta' \in \tilde \Theta_{N\deltat}(\theta_n)$;
    
  \item[\asm{P4}] there is a function $R$ satisfying $\E[\norm{F_i}_{Q_{\deltat,N}(\theta_n)}^2 \mid \theta_i = \theta]^{1/2} \leq R(\theta)$ for any $\theta \in \Theta_{\deltat,N}(\theta_n) \cup \tilde \Theta_{N\deltat}(\theta_n)$ and $R(\theta_a)^2 \leq K_1^2\norm{\theta_a - \theta_b}_{Q_{\deltat,N}(\theta_n)}^2 + K_2^2 R(\theta_b)^2$ for any $\theta_a, \theta_b \in \Theta_{\deltat,N}(\theta_n) \cup \tilde \Theta_{N\deltat}(\theta_n)$, where $K_1, K_2 \geq 0$ are constants independent of $\theta_n \in U$.
  \end{itemize}
Then, $\flow: [0, N\deltat] \to \tilde \Theta_{N\deltat}(\theta_n)$ is a unique solution satisfying \eqref{eq:integral-form} and 
\begin{equation}
\sup_{0 \leq k\leq N}\E_{n}[\abs*{\potential(\theta_{n+k}) - \potential(\flow(k \deltat;\theta_n))}] \leq (C_1 + C_2) \Delta_{L}(\deltat, N\deltat) \enspace,
\label{eq:ode-igo-bound}
\end{equation}
where 
\begin{align*}
 C_1 &= \textstyle \sum_{i=0}^{N-1} \Delta_{\tilde{L}}(\deltat, N\deltat) \int_{i \deltat}^{(i+1)\deltat} ((i+1)\deltat - t) R(\flow(t; \theta_n)) \rmd t \ , \\
  C_2 &= \textstyle (1 + (K_1 N \deltat)^2 \exp(2 K_1 N \deltat))^{1/2} K_2 (N \deltat^2)^{1/2} R(\theta_n)
 \enspace. 
\end{align*}%
\end{theorem}

\begin{corollary}\label{cor:ode-igo}
  Let $\Theta_{\deltat,N}(\theta_n)$, $\tilde \Theta_{N\deltat}(\theta_n)$, $K_1$, $K_2$, $\Delta_{L}$ and $\Delta_{\tilde{L}}$ be the same as in Theorem~\ref{lem:ode-igo}, and assume \asm{P1}, \asm{P2}, \asm{P3} and \asm{P4}.
  In addition, assume that \correct{for any $\theta_n \in \Theta$ and any $N \in \N_{+}$}
\begin{itemize}
\item[\asm{P5}] there exists a function $\Delta_{R}$ that is independent of $\theta_n \in U$, is nondecreasing \wrt each argument, and satisfies $\sup_{\theta \in \tilde \Theta_{N\deltat}(\theta_n)} R(\theta) \leq \Delta_R(\deltat, N\deltat)\potential(\theta_n)$;
\end{itemize}
then, the condition \Afou{} of Theorem~\ref{thm:llc-igo} is satisfied by
\begin{equation}
  \DeltaAfou(x, y) = [(xy / 2)\Delta_{\tilde{L}}(x, y) + K_2 (xy + K_1^2 x y^3 \exp(2 K_1 y))^{1/2} ] \Delta_{R}(x, y) \Delta_{L}(x, y)
  \enspace. 
\label{eq:ode-igo-atwo}
\end{equation}
If $U = \Theta$, condition \Atwo{} of Theorem~\ref{thm:lc-igo} as well as condition \Btwo{} of Theorem~\ref{thm:lc-igo-lower} are satisfied by $\DeltaAtwo = \DeltaBtwo = \DeltaAfou$ defined above.
\end{corollary}

Theorem~\ref{lem:ode-igo} introduces the subsets $\Theta_{\deltat,N}(\theta_n) \subseteq \Theta$ and $\tilde \Theta_{N\deltat}(\theta_n) \subseteq \Theta$.
The motivation for introducing these sets is to limit the situations one has to deal with to satisfy the conditions \asm{P1}--\asm{P4}.
It is often the case that one can bound the parameter set that $\flow(t; \theta_n)$ can reach in a finite time, whereas one may not be able to obtain a bounded set for $\Theta_{\deltat, N}(\theta_{n})$ due to its stochastic nature.
Comparing to Theorem~\ref{lem:ode}, \asm{P1} is to treat the progress in $\potential$ rather than the Euclidean distance.
The conditions \asm{P2} and \asm{P3} are weaker than the corresponding assumptions in Theorem~\ref{lem:ode}; they will be replaced by the global Lipschitz continuity $\norm{F(\theta') - F(\theta)} \leq L \norm{\theta' - \theta}$.
The function $R$ in \asm{P4} is replaced by a constant $K$ in Theorem~\ref{lem:ode}.
The condition \asm{P5} in Corollary~\ref{cor:ode-igo} is the only additional assumption posed to obtain \emph{geometric} convergence (condition \asm{A2} of Theorem~\ref{thm:lc-igo} and condition~\asm{A4} of Theorem~\ref{thm:llc-igo}).

\section{Application to Step-Size Adaptive Evolution Strategy}\label{sec:appli}


  We apply the practical conditions derived in Theorem~\ref{lem:prac-one} and Theorem~\ref{lem:ode-igo} to the SSA-ES described in \S\ref{sec:es}.

  We consider algorithm \eqref{eq:es-igo-intro} solving a spherical function $f(x) = \nu(\norm{x - \xstar})$ on $\X = \R^d$, where $\nu: \R \to \R$ is an arbitrary increasing function.
  For this $f$, we have $\qftneq{n}(s) = 0$ for any $s \in \R$ and any $\theta_n \in \Theta$. Thanks to Remark~\ref{rem:prop:exp}, we can write the mean field $F(\theta)$ as \eqref{eq:meanfield-s} with $u$ defined in \eqref{eq:binom}. For the sake of simplicity, we let $\Ratio = 1$. We let $\xstar = 0$ without loss of generality for the derivation.
  
  We assume the following conditions on $w_1, \dots, w_\lambda$.
  \asm{W1}: $w_i \geq 0$ for all $i = 1,\dots, \lambda$ and $\sum_{i=1}^{\lambda} w_i = 1$;
  \asm{W2}: $w_1 \geq \cdots \geq w_\lambda \geq 0$ and $w_i > w_{i+1}$ at least for one $i$;
  \asm{W3}: $\int u(\Phi(z)) (z^2 - 1) \psi(z) \mathrm{d}z = \bar{C}_\sigma > 0$, where $\psi$ and $\Phi$ are the probability density function and the cumulative density function of the standard normal distribution, respectively.
  \asm{W2} implies that $u$ is strictly decreasing because
  \begin{equation}
    \frac{\mathrm{d} u(p)}{\mathrm{d}p} = \lambda (\lambda - 1) \sum_{k=1}^{\lambda - 1} (-w_{k} + w_{k+1}) P_B(\lambda - 2, k-1, p) \enspace,
  \end{equation}
  all the terms are non-positive, and there is at least one negative term.

 We choose $\potential(\theta) = V(\theta)^{1/2} \cdot \max(\eta,  \xi(\theta))^{1/2}$, where
 $\eta > 0$ is some constant specified later, $V(\theta) = \norm{\mm - x^*}^2 + d \sigma^2$, and $\xi(\theta) = \norm{\mm - x^*} / \sigma$. 
 This choice is justified as follows:
 The expected distance between the optimal solution $\xstar$ and a sample $X_n \sim P_{\theta_n}$ is bounded from above, because $\E_{X_n \sim P_{\theta_n}}[\norm{X_n - \xstar}] \leq (\norm{\mm_n - \xstar}^2 + d \sigma_n^2)^{1/2} = V(\theta_n)^{1/2} \leq \potential(\theta_n) / \eta^{1/2}$.
 Therefore, the geometric convergence of $\potential(\theta_n)$ leads to the (expected) geometric convergence of $d(X_n, x^*)$ to $0$.
 Moreover, $\potential(\theta) \to + \infty$ as $\theta \to (\bar{x}, 0)$ for any $\bar{x} \neq \xstar$.
 This leads to geometric decrease of $\potential(\theta)$ at any $\theta$ along the flow $\flow$, because $\sigma$ increases geometrically if $\xi(\theta) \geq \eta$ (as proved later).

 First, we prove \asm{A1} in Theorem~\ref{thm:lc-igo} by using Theorem~\ref{lem:prac-one}.
 According to Lemma 2.9 of \cite{giorgi1992dini} (the formula for the Dini derivatives of a maximum of continuous functions), we have
\begin{align}
  D^{+}_{h} \potential(\theta) =
  \begin{cases}
    \inner{\nabla V(\theta)^{1/2}}{h} \cdot \eta^{1/2} & \xi(\theta) < \eta \\
    \inner{\nabla V(\theta)^{1/2}}{h} \cdot \eta^{1/2} + \max\{0, V(\theta)^{1/2} \cdot \inner{\nabla \xi(\theta)^{1/2}}{h} \} & \xi(\theta) = \eta \\
    \inner{\nabla V(\theta)^{1/2}}{h} \cdot \xi(\theta)^{1/2} + V(\theta)^{1/2} \cdot \inner{\nabla \xi(\theta)^{1/2}}{h} & \xi(\theta) > \eta     
  \end{cases}
\enspace,\label{eq:d_es_cases}
\end{align}
where $\nabla V(\theta)^{1/2} = (m, d\sigma) / V(\theta)^{1/2}$
and $\nabla \xi(\theta)^{1/2} = \frac12 \big(\frac{m}{\norm{m} \sigma}, -\frac{\norm{m}}{\sigma^2}\big) / \xi(\theta)^{1/2}$.
To proceed, we show the following lemmas, proved in \ref{apdx:lemma:apdx:j} and \ref{apdx:lemma:apdx:ratio}, respectively.
\begin{lemma}\label{lemma:apdx:j}
  Suppose that \asm{W2} is satisfied.
  Let $C_V = (3 + 12/d)^{-1/2}\cdot \sum_{k=1}^{\lambda} w_k \cdot (1 - 2k / (\lambda + 1))$.
  Then, for all $\theta \in \Theta$, 
  \begin{align*}
    \inner*{\nabla V(\theta) }{ F(\theta) } \leq - C_{V} \cdot \sigma \cdot V(\theta)^{1/2} \enspace.
    \end{align*}
\end{lemma}
\begin{lemma}\label{lemma:apdx:ratio}
  Suppose that \asm{W2} and \asm{W3} are satisfied.
  Then, for any $C_\sigma \in (0, \bar{C}_\sigma)$, there exists an $\eta > 0$ such that, for any $\theta$ satisfying $\xi(\theta) \geq \eta$,
  \begin{equation*}
    \inner*{ \nabla \xi(\theta) }{ F(\theta) } < - (C_\sigma / (2d)) \cdot \xi(\theta) \enspace.
  \end{equation*}
\end{lemma}
Let $C_\sigma$ and $\eta$ be the constants appearing in Lemma~\ref{lemma:apdx:ratio}.
If $\xi(\theta) \leq \eta$, we have $\sigma \geq V(\theta)^{1/2} / (\eta^2 + d)^{1/2}$. 
Using Lemma~\ref{lemma:apdx:j} and Lemma~\ref{lemma:apdx:ratio}, we have from \eqref{eq:d_es_cases}
\begin{align*}
  D^{+}_{F(\theta)}\potential(\theta)
  \leq
  \begin{cases}
    - \frac{C_V}{2 (\eta^2 + d)^{1/2}} \cdot \potential(\theta) & \xi(\theta) \leq \eta \\
    - \frac{C_\sigma}{4d} \cdot \potential(\theta) & \xi(\theta) > \eta 
  \end{cases}
\end{align*}
which satisfies the condition in Theorem~\ref{lem:prac-one} with $C = \min\left\{\frac{C_V}{2 (\eta^2 + d)^{1/2}}, \frac{C_\sigma}{4d}\right\}$. Hence, \asm{A1} of Theorem~\ref{thm:lc-igo} is satisfied with $\DeltaAone(t) = \exp(-C t)$. 

Next, we derive \asm{A2} by using Theorem~\ref{lem:ode-igo} and Corollary~\ref{cor:ode-igo}.
First, we establish the following lemma, which states that the expected squared norm of $F_n$ is upper-bounded. The proof is in \ref{apdx:lemma:es:variance}.
\begin{lemma}\label{lemma:es:variance}
  Let $C_{\ref{lemma:es:variance}} = (d + (2d)^{-1}) \cdot \lambda \cdot \sum_{i=1}^{\lambda}w_i^2$. 
  Then, $\E[\norm{F_n}^2\mid \theta_n = \theta] \leq C_{\ref{lemma:es:variance}} \cdot \sigma^2$ for all $\theta \in \Theta$.   
\end{lemma}
Using this bound, we derive the subsets $\Theta_{\alpha,N}(\theta_0)$ and $\tilde{\Theta}_{N\alpha}(\theta_0)$ in Lemma~\ref{lemma:apdx:ratio-bound}, proved in \ref{apdx:lemma:apdx:ratio-bound}.
\begin{lemma}\label{lemma:apdx:ratio-bound}
  Suppose that \asm{W1} is satisfied.
  Let $C_{\ref{lemma:es:variance}}$ be as defined in Lemma~\ref{lemma:es:variance}. 
  Then, for any $k \in \llbracket 0, N\rrbracket$, we have $\theta_{k} \in \Theta_{\alpha,N}(\theta_0)$, where
  \begin{multline*}
    \Theta_{\alpha, N}(\theta_0) := \{\theta \in \Theta:
    \sigma \geq \sigma_0 (1 - \alpha/2)^{N}
    \\
    \wedge
    \xi(\theta) \leq (1 - \alpha/2)^{-N} \xi(\theta_0) + ((1 - \alpha/2)^{-N} - 1) (d (2 - \alpha) / \alpha)^{1/2}
    \} \enspace.
  \end{multline*}
  Moreover, for any $t \in [0, N\alpha]$, we have $\flow(t; \theta_0) \in \tilde{\Theta}_{N\alpha}(\theta_0)$, where
  \begin{multline*}
    \tilde{\Theta}_{N\alpha}(\theta_0) := \{\theta \in \Theta:
    \sigma_0 \exp(- C_{\ref{lemma:es:variance}}^{1/2} \cdot N \alpha) \leq \sigma \leq \sigma_0 \exp( C_{\ref{lemma:es:variance}}^{1/2} \cdot N\alpha)
    \\
    \wedge
    \xi(\theta) \leq (1 + \xi(\theta_0)) \cdot \exp(C_{\ref{lemma:es:variance}}^{1/2}\cdot N \alpha) - 1
    \} \enspace.
  \end{multline*}
\end{lemma}
Let
\begin{multline*}
  \Delta_{\xi}(x, y) := \sup_{z \in (0, x]}\big\{
  \max\{ (1 - z/2)^{-y/z}, \exp(C_{\ref{lemma:es:variance}}^{1/2}\cdot y) \}
  \\
  + \max\{ ((1 - z/2)^{-y/z} - 1) (d (2 - z)/ z)^{1/2}, 
  \exp(C_{\ref{lemma:es:variance}}^{1/2}\cdot y) - 1 
  \} \cdot \eta^{-1} \big\}\enspace.
\end{multline*}
This is nondecreasing \wrt both arguments, and for any $\theta \in \text{conv}( \Theta_{\alpha,N}(\theta_0) \cup \tilde{\Theta}_{N\alpha}(\theta_0))$, we have $\max\{\eta, \xi(\theta)\} \leq \max\{\eta, \xi(\theta_0)\} \cdot \Delta_{\xi}(\alpha, N \alpha)$.
We obtain
\begin{align*}
  \frac{\potential(\theta) / \sigma }{\potential(\theta_0) / \sigma_0}
  &=
    \frac{(\xi(\theta)^2 + d)^{1/2} \max\{\eta, \xi(\theta)\}^{1/2}}{ (\xi(\theta_0)^2 + d)^{1/2} \max\{\eta, \xi(\theta_0)\}^{1/2}  } \\
  &=
    \frac{ (\xi(\theta)^2 + d)^{1/2} }{ \max\{\eta, \xi(\theta)\} }
    \frac{ \max\{\eta, \xi(\theta_0)\} }{ (\xi(\theta_0)^2 + d)^{1/2} }    
    \frac{ \max\{\eta, \xi(\theta)\}^{3/2} }{ \max\{\eta, \xi(\theta_0)\}^{3/2} }    
  \\
  &\leq
    ( 1 + d / \eta^2 )^{1/2} \cdot \max\{\eta / d^{1/2}, 1\} \cdot \Delta_\xi(\alpha, N\alpha)^{3/2}
    \\
  &=
    ( 1 + \max\{d / \eta^2, \eta^2 / d\} )^{1/2} \cdot \Delta_\xi(\alpha, N\alpha)^{3/2}
    \enspace.
\end{align*}
Finally, we confirm \asm{P1}--\asm{P5} one by one.

(\asm{P1})
For any $\theta \in \Theta$, we have $\norm{\nabla V(\theta)^{1/2}} \leq d^{1/2}$ and $\norm{\nabla \xi(\theta)^{1/2}} = \frac12 (\xi(\theta)^{-1} + \xi(\theta))^{1/2} / \sigma $.
Then, from \eqref{eq:d_es_cases}, we have that for $\xi(\theta) < \eta$
\begin{align*}
  \abs{D^{+}_{h}\potential(\theta)} / \norm{h}
  \leq
  (d \eta)^{1/2} \leq (V(\theta_0) / \sigma_0^2)^{1/2} \cdot \max\{\eta, \xi(\theta_0)\}^{1/2}
  =
  \potential(\theta_0) / \sigma_0 \enspace,
\end{align*}
and for $\xi(\theta) \geq \eta$
\begin{align*}
  &\abs{D^{+}_{h}\potential(\theta)} / \norm{h}
    \\
  &\leq
    d^{1/2} \xi(\theta)^{1/2} + (1/2) V(\theta)^{1/2} \cdot (\eta^{-1} + \xi(\theta))^{1/2} / \sigma
  \\
  &\leq
    ( 1 + (1 + \eta^{-2})^{1/2}/ 2 ) \potential(\theta) / \sigma
  \\
  &\leq
    ( 1 + (1 + \eta^{-2})^{1/2}/ 2 ) \cdot ( 1 + \max\{d / \eta^2, \eta^2 / d\} )^{1/2} \cdot \Delta_\xi(\alpha, N\alpha)^{3/2} \cdot (\potential(\theta_0) / \sigma_0) \enspace.
\end{align*}
Let $\beta = ( 1 + (1 + \eta^{-2})^{1/2}/ 2 ) \cdot ( 1 + \max\{d / \eta^2, \eta^2 / d\} )^{1/2}$.
Then $\beta \geq 1$ and $\Delta_\xi(\alpha, N\alpha) \geq 1$. Therefore, we have $\abs{D^{+}_{h}\potential(\theta)} / \norm{h} \leq \beta \cdot \Delta_\xi(\alpha, N\alpha)^{3/2} \cdot (\potential(\theta_0) / \sigma_0)$ for all $\theta \in \text{conv}(\Theta_{\alpha,N}(\theta_0) \cup \tilde{\Theta}_{N\alpha}(\theta_0))$.
Let $Q_{\alpha,N}(\theta_0) = \big( \beta \cdot \Delta_\xi(\alpha, N\alpha)^{3/2} \cdot (\potential(\theta_0) / \sigma_0) \big)^2 \cdot I_\dimtheta$. 
Then, we obtain $\sup_{\norm{h} = 1} \abs{D^{+}_{h}\potential(\theta)} \cdot \norm{h}_{Q_{\alpha,N}(\theta_0)^{-1}} \leq 1$ for all $\theta \in \text{conv}(\Theta_{\alpha,N}(\theta_0) \cup \tilde{\Theta}_{N\alpha}(\theta_0))$. 

(\asm{P2} and \asm{P3}) The Lipschitz continuity of $F$ follows from Lemma~\ref{lemma:apdx:lip} (proved in \ref{apdx:lemma:apdx:lip}):
\begin{lemma}\label{lemma:apdx:lip}
  Let $L_{\ref{lemma:apdx:lip}} = (3d^2+3/2)^{1/2} \cdot \lambda (\lambda - 1) \max_{i=1,\dots,\lambda-1} \abs{w_{i+1} - w_{i}} + 2 \lambda \max_{i=1,\dots,\lambda}\abs{w_i}$. Then, 
$\norm{F(\theta_1) - F(\theta_2)} \leq L_{\ref{lemma:apdx:lip}} \norm{\theta_1 - \theta_2}$ for any $\theta_1 \in \Theta$ and any $\theta_2 \in \Theta$.
\end{lemma}
Hence, \asm{P2} is satisfied, with $L_{\alpha,N}(\theta, \theta_n) = L_{\ref{lemma:apdx:lip}}$ and $\Delta_{L}(x, y) = \exp(L_{\ref{lemma:apdx:lip}} y)$.
Moreover, \asm{P3} is satisfied, with $\Delta_{\bar{L}}(x, y) = L_{\ref{lemma:apdx:lip}}$.

(\asm{P4})
Let $R(\theta) = C_{\ref{lemma:es:variance}}^{1/2} \cdot \big( \beta \cdot \Delta_\xi(\alpha, N\alpha)^{3/2} \cdot (\potential(\theta_0) / \sigma_0) \big) \cdot \sigma$.
Because of Lemma~\ref{lemma:es:variance}, we have $\E_n[\norm{F_n}_{Q_{\alpha,N}(\theta_0)}^2]^{1/2} \leq R(\theta_n)$. 
Since $\sigma_a^2 \leq 2 ((\sigma_a - \sigma_b)^2 + \sigma_b^2)$, we have $R(\theta_a)^2 \leq 2 C_{\ref{lemma:es:variance}} \norm{\theta_a - \theta_b}_{Q_{\alpha,N}(\theta_0)}^2 + 2 R(\theta_b)^2$. Therefore, with $K_1 = (2 C_{\ref{lemma:es:variance}})^{1/2}$ and $K_2 = 2^{1/2}$, \asm{P4} is satisfied.

(\asm{P5})
Because, as shown in Lemma~\ref{lemma:apdx:ratio-bound}, $\sigma \leq \sigma_0 \cdot \exp(C_{\ref{lemma:es:variance}}^{1/2} \cdot N \alpha)$ for all $\theta \in \tilde{\Theta}_{N\alpha}(\theta_0)$, we have
\begin{align}
  \sup_{\theta \in \tilde{\Theta}_{N\alpha}(\theta_0)} R(\theta) 
\leq C_{\ref{lemma:es:variance}}^{1/2} \cdot \beta \cdot \Delta_\xi(\alpha, N\alpha)^{3/2} \cdot \exp(C_{\ref{lemma:es:variance}}^{1/2} \cdot N \alpha)\cdot \potential(\theta_0)  \enspace.
\end{align}
Let $\Delta_R(x, y) = C_{\ref{lemma:es:variance}}^{1/2} \cdot \beta \cdot \Delta_\xi(x, y)^{3/2} \cdot \exp(C_{\ref{lemma:es:variance}}^{1/2} \cdot y)$. This is nondecreasing \wrt both arguments. Moreover, we have $\sup_{\theta \in \tilde{\Theta}_{N\alpha}(\theta_0)} R(\theta) \leq \Delta_R(\alpha, N\alpha) \cdot \potential(\theta_0)$. Therefore, (P5) is satisfied with $\Delta_R(x, y)$. 
Now that we have confirmed (\asm{P1})--(\asm{P5}), \Atwo{} is satisfied with $\DeltaAtwo$ defined in \eqref{eq:ode-igo-atwo}.

In sum, we have shown the result formalized in the following proposition:
\begin{theorem}[Geometric Convergence of Step-Size Adaptive Evolution Strategy]\label{prop:lin-conv-es}
  Consider the SSA-ES \eqref{eq:es-igo-intro} solving a spherical function $f(x) = \nu(\norm{x - \xstar})$ on $\X = \R^d$, where $\nu: \R \to \R$ is an arbitrary increasing function.
  Suppose that $w_1,\dots, w_\lambda$ are selected so that \asm{W1}, \asm{W2} and \asm{W3} described above are satisfied.
  Let $\potential$ be as defined above. 
  Then, there exists an $\bar \deltat > 0$ such that $\limsup_{n\to\infty} \frac{1}{n} \ln( \E_0[\potential(\theta_n)] ) < 0$ for any $\theta_0 \in \Theta$ and for any $\deltat \in (0, \bar\deltat]$.
  That is, $\limsup_{n\to\infty} \frac{1}{n} \ln( \E_0[\norm{X_n - \xstar}] ) < 0$, where $X_n \sim P_{\theta_n}$. 
\end{theorem}

To the best of our knowledge, this is the first result showing the geometric convergence of the SSA-ES of form \eqref{eq:es-igo-intro} derived from the IGO principle. 
It has been shown in \cite{10.1007/978-3-642-32937-1_1,Akimoto2012ppsn} that the flow $\flow(t; \theta)$ converges towards $\theta^* = (x^*, 0)$, however, the speed of convergence has not been discussed and the convergence of the flow has not been related to the convergence of the recursive algorithm.
Related works \cite{10.1145/3205455.3205606,10.1145/3299904.3340303,10.1145/3449639.3459289} have shown the geometric convergence of (1+1)-ES---a variant of SSA-ES using elitism---by applying the drift analysis without going through the ODE method.
Such an analysis is difficult to apply to a strategy without elitism such as \eqref{eq:es-igo-intro} or CMA-ES.
Since CMA-ES is partly derived from the IGO principle, we expect that our approach is applicable to proving the geometric convergence of the CMA-ES as well.

\section{Local Geometric Convergence}\label{sec:llc}

Theorem~\ref{thm:lc-igo} provides a sufficient condition for the \emph{global} geometric convergence of $\potential(\theta_n)$.
The assumptions \Aone{} and \Atwo{} pose conditions over all $\theta \in \Theta$.
For \emph{local} geometric convergence of $\potential(\theta_n)$, we can relax these requirements by posing conditions on a subset $U \subseteq \Theta$, rather than the whole parameter set $\Theta$, as is stated in Theorem~\ref{thm:llc-igo} and Corollary~\ref{cor:llc-igo} below.

\begin{theorem}
\label{thm:llc-igo}
Instead of conditions \Aone{} and \Atwo{} in Theorem~\ref{thm:lc-igo}, assume that there exists an open subset $U \subseteq \Theta$ that satisfies 
the following: 
$[\theta : \potential(\theta) < \supU] \subseteq U$ for some $\supU > 0$, and 
\begin{description}
\item[\Athr] There exists $\DeltaAthr: \R_{\geq 0} \to \R_{\geq 0}$ nonincreasing such that $\DeltaAthr(t)\downarrow 0$ as $t \uparrow \infty$ and for any $\theta \in U$ and for any $t \geq 0$
\begin{equation}
\label{eq:Athr}
\potential(\flow(t; \theta)) \leq \DeltaAthr(t) \potential(\theta) \enspace.
\end{equation}
\item[\Afou] There exists $\DeltaAfou: \R_{+} \times \R_{+} \to \R_{+}$ nondecreasing \wrt each argument such that $\DeltaAfou(\deltat, T) \downarrow 0$ as $\deltat \downarrow 0$ for any fixed $T > 0$, and for any $N \in \N_{+}$ and $\theta_0 \in U$
\begin{equation}
\E_{0}[\potential(\theta_{N})] \leq \potential(\flow(N \deltat; \theta_{0}))  + \DeltaAfou(\deltat, N\deltat) \potential(\theta_{0})
\enspace.
\label{eq:Afou}
\end{equation}%
\end{description}

Let $\CR = \inf_{N \geq 1}(\DeltaAthr(N\deltat) + \DeltaAfou(\deltat, N \deltat))^{1/N}$ and $\deltatset = \{ \deltat > 0: \CR < 1\}$.
Then, $\deltatset$ is nonempty, and for any $\deltat \in \deltatset$ there exists at least one $N \geq 1$ such that $(\DeltaAthr(N\deltat) + \DeltaAfou(\deltat, N\deltat ))^{1/N} = \CR<1$.
Let the minimum of such $N$ be denoted by $\optN$. Let $\remCR = 1$ if $\optN = 1$ and $\remCR = \max_{N \in \llbracket 0, \optN-1\rrbracket}(\DeltaAthr(N\deltat) + \DeltaAfou(\deltat, N \deltat))$ otherwise.
Then, the following hold for any $\deltat \in \deltatset$ and $\theta_0 \in U$:
\begin{enumerate}
\item Let $\EventU_{k}$ be the event that $\{\theta_{i \optN}\}_{i \in \llbracket 0, k\rrbracket}$ stay in $U$, i.e., $\EventU_{k} = \bigcap_{i \in \llbracket 0, k\rrbracket} [ \theta_{i \optN} \in U]$.
  Then, for any $k \geq 0$, $\Pr[\EventU_{ k}] \geq 1 - (\potential(\theta_{0}) / \supU) (\CR^{\optN} - \CR^{\optN (k+1)})/ (1 - \CR^{\optN})$ and $\Pr[\EventU_{\infty}] \geq 1 - (\potential(\theta_{0}) / \supU) \CR^{\optN}/ (1 - \CR^{\optN})$ for $\EventU_{\infty} = \lim_{k\to\infty} \EventU_{k}$.

\item $\E[\potential(\theta_{n}) \1{\EventU_{\lfloor (n-1) / \optN \rfloor}}] \leq (\remCR / \CR^{\optN-1}) \CR^{n} \potential(\theta_{0})$ for any $n \in \N_{+}$. 

\item $\Pr[\potential(\theta_{n}) < \epsilon] \geq \Pr[\EventU_{\lfloor (n-1) / \optN \rfloor}] - (\potential(\theta_{0}) / \epsilon) \CR^{n} (\remCR/\CR^{\optN-1})$ for any $n \in \N_{+}$ and $\epsilon > 0$. 
\end{enumerate}
\end{theorem}

The second statement is the counterpart of the consequence of Theorem~\ref{thm:lc-igo}; however, $\potential(\theta_{n})$ is multiplied by $\1{\EventU_{\lfloor (n-1) / \optN \rfloor}}$.
This is because $\potential(\theta_n)$ may not converge toward zero and may leave the neighborhood $U$, which is intuitively the basin of attraction of the desired point $\theta^* \in \Theta$.
The first statement provides a lower bound on the probability that all $\theta_{i\optN}$ for $i \geq 0$ stay in $U$.
The statement reads that this probability can be arbitrarily close to $1$ if we take the initial point $\theta_0$ such that $\potential(\theta_0)$ is sufficiently small.
However, there is always a positive probability that $\theta_n$ leaves $U$.
The third statement is the most interesting result.
The probability of $\potential(\theta_n)$ being smaller than a given $\epsilon > 0$ will eventually be lower-bounded by the probability that all $\theta_{i\optN}$ stay in $U$, which is lower-bounded according to the first statement.
Roughly speaking, $\potential(\theta_n)$ converges towards $0$ as long as $\theta_{i\optN}$ stays in $U$ for all $i \geq 0$.
The asymptotic results are derived in Corollary~\ref{cor:llc-igo} after the proof of Theorem~\ref{thm:llc-igo}.

\begin{proof}
  In a way analogous to the proof of Theorem~\ref{thm:lc-igo}, we can show that $\deltatset$ is nonempty and there exists at least one $N \in \N_{+}$ that minimizes $(\DeltaAthr(N\deltat) + \DeltaAfou(\deltat, N \deltat))^{1/N}$.
  Moreover, by applying the same argument, we find that $\E_{n}[\potential(\theta_{n+\optN})] \leq \CR^{\optN} \potential(\theta_{n})$ given $\theta_{n} \in U$.
  In other words, $\E_{n}[\potential(\theta_{n+\optN})]\1{\theta_{n} \in U} \leq \CR^{\optN} \potential(\theta_{n}) \1{\theta_{n} \in U}$.

  First, we prove the second statement in the theorem.
  Noting that $\1{\EventU_{k}} = \1{\EventU_{k-1}} \1{\theta_{k \optN} \in U} =\prod_{i=0}^{k} \1{\theta_{i \optN} \in U}$, the inequality $\E_{n}[\potential(\theta_{n+\optN})]\1{\theta_{n} \in U} \leq \CR^{\optN} \potential(\theta_{n}) \1{\theta_{n} \in U}$ implies
\begin{equation}
\label{eq:llc-igo:exp}
\begin{split}
\E[\potential(\theta_{\optN k}) \1{\EventU_{k-1}}] 
&= \E[\potential(\theta_{\optN k}) \1{\EventU_{k-2}}\1{\theta_{\optN (k-1)} \in U}] \\
&= \E[ \E_{\optN (k-1)}[\potential(\theta_{\optN k})] \1{\theta_{\optN (k-1)} \in U} \1{\EventU_{k-2}}] \\
&\leq \CR^{\optN} \E[ \potential(\theta_{\optN (k-1)}) \1{\theta_{\optN (k-1)} \in U} \1{\EventU_{k-2}}] \\
&\leq \CR^{\optN} \E[ \potential(\theta_{\optN (k-1)}) \1{\EventU_{k-2}}] \leq \cdots \leq \CR^{\optN k} \potential(\theta_{0}) \enspace.
\end{split}
\end{equation}
On the other hand, since $\remCR \geq (\DeltaAthr(n\deltat) + \DeltaAfou(\deltat, n \deltat))$ for any $n \in \llbracket 0, \optN-1\rrbracket$, we have that for such $n$ and for any $k \geq 0$, $\E_{\optN k}[\potential(\theta_{n + \optN k})] \1{\EventU_{k}} \leq \remCR \potential(\theta_{\optN k}) \1{\EventU_{k}}$.
Combining this with \eqref{eq:llc-igo:exp}, we obtain 
\begin{multline}
  \E[\potential(\theta_{n + \optN k}) \1{\EventU_{k}}]
  \leq \remCR \E[\potential(\theta_{\optN k}) \1{\EventU_{k}}]
  \leq \remCR \E[\potential(\theta_{\optN k}) \1{\EventU_{k-1}}]  \\
  \leq \remCR \CR^{\optN k} \potential(\theta_{0}) 
= (\remCR / \CR^n) \CR^{n + \optN k} \potential(\theta_{0}) \leq (\remCR / \CR^{\optN - 1}) \CR^{n + \optN k} \potential(\theta_{0}) \enspace.
\label{eq:llc-igo:exp2}
\end{multline}
Note that $\remCR / \CR^{\optN - 1} \geq 1$. Rewriting $n + \optN k$ as $n$, from \eqref{eq:llc-igo:exp} and \eqref{eq:llc-igo:exp2} we have that $\E[\potential(\theta_{n}) \1{\EventU_{\lfloor(n-1)/\optN\rfloor}}] \leq (\remCR / \CR^{\optN - 1}) \CR^{n} \potential(\theta_{0})$. 

Next, we prove the third statement of the theorem.
To do so, we will find a lower bound on the probability of $\potential(\theta_{n+\optN k})$ being smaller than a given $\epsilon > 0$. For $n = 0$,
\begin{equation}
\begin{split}
\Pr[\potential(\theta_{\optN k}) < \epsilon]
&\geq \Pr[\potential(\theta_{\optN k})\1{\EventU_{k-1}} < \epsilon \1{\EventU_{k-1}}] \\
&= 1 - \Pr[\potential(\theta_{\optN k})\1{\EventU_{k-1}} \geq \epsilon \1{\EventU_{k-1}}] \\
&= 1 - \Pr[\potential(\theta_{\optN k})\1{\EventU_{k-1}} + \epsilon (1-\1{\EventU_{k-1}}) \geq \epsilon ] \\
&\geq 1 - \E[\potential(\theta_{\optN k})\1{\EventU_{k-1}}]/\epsilon - \E[\epsilon (1-\1{\EventU_{k-1}})] / \epsilon \\
&= \Pr[\EventU_{k-1}]  - \E[\potential(\theta_{\optN k})\1{\EventU_{k-1}}] / \epsilon \\
&\geq \Pr[\EventU_{k-1}]  - \CR^{\optN k} \potential(\theta_{0}) / \epsilon \enspace,
\end{split}
\label{eq:llc-igo:prob}
\end{equation}
and similarly for $n \in \llbracket 1, \optN-1 \rrbracket$,
\begin{equation}
\begin{split}
\Pr[\potential(\theta_{n+\optN k}) < \epsilon]
&\geq \Pr[\potential(\theta_{n+\optN k})\1{\EventU_{k}} < \epsilon \1{\EventU_{k}}] \\
&= 1 - \Pr[\potential(\theta_{n+\optN k})\1{\EventU_{k}} \geq \epsilon \1{\EventU_{k}}] \\
&= 1 - \Pr[\potential(\theta_{n+\optN k})\1{\EventU_{k}} + \epsilon (1-\1{\EventU_{k}}) \geq \epsilon ] \\
&\geq 1 - \E[\potential(\theta_{n+\optN k})\1{\EventU_{k}} + \epsilon (1-\1{\EventU_{k}})] / \epsilon \\
&= \Pr[\EventU_{k}]  - \E[\potential(\theta_{n+\optN k})\1{\EventU_{k}}] / \epsilon \\
&\geq \Pr[\EventU_{k}]  - (\remCR / \CR^{\optN - 1}) \CR^{n+\optN k} \potential(\theta_{0}) / \epsilon \enspace.
\end{split}
\label{eq:llc-igo:prob2}
\end{equation}
Here we have applied Markov's inequality and the inequalities \eqref{eq:llc-igo:exp} and \eqref{eq:llc-igo:exp2}.
Rewriting $n + \optN k$ as $n$, \eqref{eq:llc-igo:prob} and \eqref{eq:llc-igo:prob2} imply that 
\begin{equation*}
\Pr[\potential(\theta_{n}) < \epsilon] \geq \Pr[\EventU_{\lfloor (n-1)/\optN \rfloor}]  - (\remCR / \CR^{\optN - 1}) \CR^{n} \potential(\theta_{0}) / \epsilon 
\enspace.
\end{equation*}

Last, we prove the first statement.
The definition of $\supU$ guarantees that $\theta \in U$ if $\potential(\theta) < \supU$.
Then, $\Pr[\EventU_{k}] = \E[\1{\EventU_{k}}] = \E[\1{\theta_{\optN k} \in U}\1{\EventU_{k-1}}] \geq \E[\1{\potential(\theta_{\optN k}) < \supU}\1{\EventU_{k-1}}]$.
The inside of the expectation on the RMS is equal to $\1{\potential(\theta_{\optN k})\1{\EventU_{k-1}} < \supU\1{\EventU_{k-1}}}$ and its expectation is the probability $\Pr[\potential(\theta_{\optN k})\1{\EventU_{k-1}} < \supU\1{\EventU_{k-1}}]$.
As is seen in \eqref{eq:llc-igo:prob}, this probability is lower-bounded by $\Pr[\EventU_{k-1}] - \CR^{\optN k} \potential(\theta_{0}) / \supU$.
Therefore, we have $\Pr[\EventU_{k}] \geq \Pr[\EventU_{k-1}] - \CR^{\optN k} \potential(\theta_{0}) / \supU$.
Hence, noting that $\Pr[\EventU_{0}] = \Pr[\theta_{0} \in U] = 1$, we finally have that, for $k \in \N$,
\begin{equation}
\Pr[\EventU_{k}] \geq 1 - \frac{\potential(\theta_{0})}{\supU} \sum_{i=1}^{k} \CR^{\optN i}
= 1 - \frac{\potential(\theta_{0})}{\supU} \frac{\CR^{\optN} - \CR^{\optN (k+1)}}{1 - \CR^{\optN}}
\enspace.
\label{eq:llc-igo:probU}
\end{equation}
Since $\EventU_{k+1} \subseteq \EventU_{k}$ for any $k \geq 0$, from the monotone continuity of a measure from above, we have $\Pr[\EventU_{\infty}] = \Pr[\lim_{k\to \infty}\cap_{i=0}^{k} \EventU_{k}] = \lim_{k\to\infty}\Pr[\EventU_{k}]$.
Therefore, from \eqref{eq:llc-igo:probU} we have $\Pr[\EventU_{\infty}] \geq 1 - (\potential(\theta_{0}) / \supU) \CR^{\optN}/ (1 - \CR^{\optN})$.
This completes the proof.
\end{proof}

\begin{corollary}
\label{cor:llc-igo}
Suppose that the assumptions of Theorem~\ref{thm:llc-igo} hold.
Then, for any $\deltat \in \deltatset$ and for any $\theta_0 \in U$ the following hold:
\begin{enumerate}
\item $\liminf_{n} \Pr[\potential(\theta_{n}) < \epsilon_{n}] \geq \Pr[\EventU_{\infty}]$ for any sequence $\{\epsilon_n > 0\}_{n\geq0}$ satisfying $\limsup_n \CR^n / \epsilon_n = 0$.

\item $\liminf_{n} \frac{1}{n} \ln \potential(\theta_{n}) \leq \ln \CR$ with probability at least $\Pr[\EventU_{\infty}]$.
\end{enumerate}
\end{corollary}

\begin{proof}
The first statement is immediately obtained by substituting $\epsilon = \epsilon_n$ in the third statement of Theorem~\ref{thm:llc-igo} and taking the limit infimum,  
\begin{equation*}
\liminf_{n} \Pr[\potential(\theta_{n}) < \epsilon_{n}] 
  \geq \liminf_{n} \Pr[\EventU_{\lfloor n / \optN \rfloor}] - \limsup_{n} \potential(\theta_{0}) \frac{\CR^{n}}{\epsilon_{n}}\frac{\remCR}{\CR^{\optN-1}} 
  = \Pr[\EventU_{\infty}] 
  \enspace.
\end{equation*}

For the second statement, substituting $\epsilon = a_{n}\CR^{n}$ in the third statement of Theorem~\ref{thm:llc-igo}, where $a_{n}$ is an increasing sequence with $a_{n} \uparrow \infty$ and $\frac{1}{n} \ln a_{n} \downarrow 0$, yields $\Pr[\potential(\theta_{n}) < a_{n}\CR^{n}] \geq \Pr[\EventU_{\lfloor n / \optN \rfloor}] - \potential(\theta_{0}) (\remCR/\CR)^{\optN-1} / a_{n}$.
Let $A_{n} = [\potential(\theta_{n}) < a_{n}\CR^{n}]$.
By the reverse Fatou's lemma, we have $\Pr[\limsup_{n} A_{n}] \geq \limsup_{n} \Pr[A_{n}]$.
The RHS is bounded below by 
\begin{equation*}
\limsup_{n} \big(\Pr[\EventU_{\lfloor n / \optN \rfloor}] - \potential(\theta_{0}) (\remCR/\CR)^{\optN-1} / a_{n} \big)
= \limsup_{n} \Pr[\EventU_{\lfloor n / \optN \rfloor}] = \Pr[\EventU_{\infty}]
\enspace,
\end{equation*}
and the inside of the probability on the left-hand side (LHS) is 
$\limsup_{n} A_{n} 
= \cap_{k=1}^{\infty}\cup_{n \geq k} A_{n} 
= \cap_{k=1}^{\infty}\cup_{n \geq k} [\frac{1}{n}\ln\potential(\theta_{n}) - \frac{1}{n} \ln a_{n} < \ln\CR ] 
\subseteq [\liminf_{n} \frac{1}{n}\ln\potential(\theta_{n}) \leq \ln\CR ] $.
This inclusion implies that the probability of the event $\limsup_{n} A_{n}$ 
satisfies $\Pr[\limsup_{n} A_{n}] \leq \Pr[\liminf_{n}\frac{1}{n}\ln(\potential(\theta_{n}) / \potential(\theta_{0})) \leq \ln\CR]$. 
Hence, $\Pr[\liminf_{n}\frac{1}{n}\ln\potential(\theta_{n}) \leq \ln\CR] \geq \Pr[\EventU_{\infty}]$. 
This ends the proof.
\end{proof}

\subsection{Consequences: Geometric Convergence and Hitting-Time Bound}

Consider a\del{n}{} \new{rank-based stochastic} algorithm of the form~\eqref{eq:algo-rank} minimizing a deterministic function $f: \X \to \R$.
Let $d: \X \times \X \to \R_{\geq0}$ be a distance function on $\X$ and $x^* \in \X$ be the well-defined optimum of the objective function, $x^* = \argmin_{x \in \X} f(x)$. 

As consequences of Theorem~\ref{thm:llc-igo} and Corollary~\ref{cor:llc-igo}, we obtain the following convergence results of the sequence of the samples $\{X_n \in \X\}$, which is an immediate consequence of the second statement of Corollary~\ref{cor:llc-igo}.
\begin{corollary}
If we can choose $\potential(\theta_n)$ such that it satisfies the assumptions of Theorem~\ref{thm:llc-igo} and if there is a constant $C > 0$ such that $\E_n[d(X_n, x^*)] \leq C \potential(\theta_n)$ for any $\theta_n \in U$, then $\forall \deltat \in \deltatset$
\begin{equation*}
  \liminf_{n} \frac1n \ln \E_{\new{n}}[d(X_n, x^*)] \leq \ln \CR \quad \text{with probability at least } \Pr[\EventU_{\infty}]
\enspace.
\end{equation*}
where $\deltatset$, $\CR$ are defined in Theorem~\ref{thm:llc-igo}.
\end{corollary}%

The above corollary implies that there exists at least one subsequence $\{n_k\}_{k \in \N_{+}}$ such that $\lim_{k\to\infty} \frac{1}{n_{k}} \ln\E_{\new{n_k}}[d(X_{n_k}, x^*)] \leq \ln \CR$ with probability at least $\Pr[\EventU_{\infty}]$, where a lower bound for $\Pr[\EventU_{\infty}]$ is provided in the first statement of Theorem~\ref{thm:llc-igo}.
In optimization settings, it is often sufficient to find a subsequence of solutions that converges geometrically towards the optimum;
however, it is not likely to happen in practical algorithms that only some subsequences converge geometrically.

The last consequence concerns the first-hitting-time bound.
Recall that the first hitting time $\tau$ defined in~\eqref{eq:fht} is the random variable that is the number of iterations spent before $X_i$ visits the $\epsilon$-neighborhood $B_\epsilon(x^*) = [x \in \X: d(x, x^*) < \epsilon]$ for the first time.
\begin{corollary}
Given $\delta \in (0, 1)$ and $\epsilon > 0$, the first hitting time $\tau$ defined in~\eqref{eq:fht} is upper-bounded by
\begin{equation*}
\bar \tau = 1 + \left\lceil\frac{\ln(\potential(\theta_0)/\epsilon)}{\ln(\CR^{-1})} + \frac{\ln(\Pr[\EventU_{\infty}]^{-1})}{\ln(\CR^{-1})} + \frac{\ln((1 - (1-\delta)^{1/2})^{-2})}{\ln(\CR^{-1})} + \frac{\ln( C \remCR/\CR^{\optN-1})}{\ln(\CR^{-1})}\right\rceil
\end{equation*}%
with probability at least $(1 - \delta) \Pr[\EventU_{\infty}]$, i.e., $\Pr[\tau \leq \bar \tau] \geq (1 - \delta) \Pr[\EventU_{\infty}]$.
\end{corollary}%

\begin{proof}
The probability of $X_n \in B_\epsilon(x^*)$ is lower-bounded by 
\begin{align*}
 \MoveEqLeft[2] \Pr[d(X_n, x^*) < \epsilon] \\
  &\textstyle= \E[ \ind{d(X_n, x^*) < \epsilon} ] \\
  &\textstyle\geq \E[ \ind{d(X_n, x^*) < \beta \E_n[d(X_n, x^*)]}\ind{\beta \E_n[d(X_n, x^*)] < \epsilon} ] \\
  &\textstyle= \E[ \E_n[\ind{d(X_n, x^*) < \beta \E_n[d(X_n, x^*)]}]\ind{\beta \E_n[d(X_n, x^*)] < \epsilon} ] \\
  &\textstyle= \E[ \Pr_n[ d(X_n, x^*) < \beta \E_n[d(X_n, x^*)]]\ind{\beta \E_n[d(X_n, x^*)] < \epsilon} ] \\
  &\textstyle= \E[ (1 - \Pr_n[ d(X_n, x^*) \geq \beta \E_n[d(X_n, x^*)]])\ind{\beta \E_n[d(X_n, x^*)] < \epsilon} ] \\
  &\textstyle\geq \E[ (1 - 1/\beta)\ind{\beta \E_n[d(X_n, x^*)] < \epsilon} ] \\
  &\textstyle= (1 - 1 / \beta) \E[\ind{\beta \E_n[d(X_n, x^*)] < \epsilon} ] \\
  &\textstyle= (1 - 1 / \beta) \Pr[\beta \E_n[d(X_n, x^*)] < \epsilon] \\
  &\textstyle\geq (1 - 1 / \beta) \Pr[\potential(\theta_n) < \epsilon/(\beta C) ] 
\end{align*}
for any $\beta \geq 1$.
The RHS is further bounded from below by using the first and third statements of Theorem~\ref{thm:llc-igo}: 
\begin{align*}
\Pr[d(X_n, x^*) < \epsilon] 
	&\textstyle\geq \sup_{\beta \geq 1}(1 - 1 / \beta) (\Pr[\EventU_{\infty}] - (\beta C \potential(\theta_{0}) / \epsilon) \CR^{n} (\remCR/\CR^{\optN-1})) \\
	&=\textstyle \max(0, \Pr[\EventU_{\infty}]^{1/2} - [(C \potential(\theta_{0}) / \epsilon) (\remCR/\CR^{\optN-1}) \CR^{n}]^{1/2} )^2
	\enspace.
\end{align*}
By solving $\Pr[\EventU_{\infty}]^{1/2} - [(C \potential(\theta_{0}) / \epsilon) (\remCR/\CR^{\optN-1}) \CR^{n}]^{1/2}  
\geq (1 - \delta)^{1/2}\Pr[\EventU_{\infty}]^{1/2}$ w.r.t.~$n$ for a given $\delta \in (0, 1)$ and incrementing it by 1, we obtain $\bar \tau$. 
\end{proof}

\section{Conclusion}\label{sec:conc}

We propose a novel methodology for proving the geometric convergence of adaptive stochastic algorithms, especially comparison-based ones for deterministic optimization problems.
The methodology is based on the ODE method, which relates the stochastic algorithm to its associated ordinary differential equation.
The main theorem, Theorem~\ref{thm:lc-igo}, provides sufficient conditions for an algorithm to exhibit geometric convergence, with an upper bound on the convergence rate.
A lower bound on the convergence rate is derived under similar sufficient conditions.
Theorem~\ref{lem:prac-one} and Theorem~\ref{lem:ode-igo} provide practically verifiable sufficient conditions to obtain both upper and lower bounds on the geometric convergence rate.
The use of the practical conditions has been illustrated on \correct{a comparison-based stochastic search algorithm on continuous domains, namely the step-size adaptive evolution strategy derived from the information-geometric optimization framework. To the best of our knowledge, this is the first result proving its geometric convergence.}
We further extend the methodology to cover the case of local convergence, in which geometric convergence is observed only if the initial parameter of the algorithm is close enough to the target parameter value.


\appendix

\section{Proofs of the Main Results}
\subsection{Proof of Proposition~\ref{prop:exp}}\label{apdx:prop:exp}

The conditional expectation of $Y_n$ is
\begin{align}
\E_n[Y_n]
  &=\textstyle \E_n\big[\sum_{i=1}^{\lambda} W(i; X_{n,1}, \dots, X_{n,\lambda}) g(X_{n,i}; \theta_n) \big] \notag\\
  &=\textstyle \sum_{i=1}^{\lambda} \E_n\left[W(i; X_{n,1}, \dots, X_{n,\lambda}) g(X_{n,i}; \theta_n) \right] \notag\\
  &=\textstyle \sum_{i=1}^{\lambda} \E_{X_{n,i} \sim P_{\theta_n}}\left[\E_{\{X_{n,k} \sim P_{\theta_n}\}_{k\neq i}}\left[W(i; X_{n,1}, \dots, X_{n,\lambda})\right] g(X_{n,i}; \theta_n) \right] \notag\\
  &=\textstyle \lambda\E_{X_{n,1} \sim P_{\theta_n}}\left[\E_{\{X_{n,k} \sim P_{\theta_n}\}_{k\geq 2}}\left[W(1; X_{n,1}, \dots, X_{n,\lambda})\right] g(X_{n,1}; \theta_n) \right] 
  \label{eq:igo-exp}
  \enspace.
\end{align}
Since $\{ X_{n,i} \}$ are i.i.d.~from $P_{\theta_n}$, the probabilities of $\1{f(X_{n,j}) < f(X_{n,1})} = 1$ and $\1{f(X_{n,j}) = f(X_{n,1})} = 1$ for each $j \neq 1$ given $X_{n,1}$ are $\qftle(f(X_{n,1}))$ and $\qfteq(f(X_{n,1}))$, respectively.
Then, the joint probability of the sums $\sum_{j=2}^{\lambda} \1{f(X_{n,j}) < f(X_{n,1})}$ and $\sum_{j=2}^{\lambda} \1{f(X_{n,j}) = f(X_{n,1})}$ being $k \in \llbracket 0, \lambda - 1 \rrbracket$ and $l \in \llbracket 0, \lambda - k - 1\rrbracket$, respectively, is given by $\trinomial(\lambda-1, k,l, p, q)$ with $p = \qftnle{n}(f(X_{n,1}))$ and $q = \qftneq{n}(f(X_{n,1}))$.
Then, $\E[ W(1; X_{n,1}, \dots, X_{n,\lambda}) \mid X_{n,1}]$ can be written as the sum of the product of $\sum_{j=k}^{k+l} \ww_{j+1}/(l+1)$ and $\trinomial(k,l; \lambda-1, \qftnle{n}(f(X_{n,1})), \qftneq{n}(f(X_{n,1})))$ over $k$ and $l$, $0 \leq k \leq \lambda-1$, $0 \leq l \leq \lambda-k-1$, resulting in $u(\qftnle{n}(f(X_{n,1})), \qftneq{n}(f(X_{n,1}))) / \lambda$.
Substituting this into \eqref{eq:igo-exp}, we obtain the desired equality.\qed

\subsection{Proof of Theorem~\ref{lem:prac-one}}\label{apdx:prac-one}

By Theorem~4.3 in Appendix~1 of \cite{rouche1977stability}, it is known that for an absolutely continuous function $\flow(\cdot; \theta): \R \to \R^\dimtheta$ and a locally Lipschitz function $\potential: \R^\dimtheta \to \R$, the upper right Dini derivative of $\potential(\flow(t; \theta))$ \wrt $t$ (defined as $D^+ \potential(\flow(t; \theta)) = \limsup_{h \downarrow 0} [\potential(\flow(t+h; \theta)) - \potential(\flow(t; \theta))]/h$) is equal to the upper Dini directional derivative of $\potential$ at $\flow(t; \theta)$ in the direction of $F(\flow(t; \theta))$ almost everywhere in $t$.
Since $\flow(\cdot; \theta)$ is a solution of the ODE \eqref{eq:associated-ode}, it is by definition absolutely continuous \wrt $t$. Then, we have
\begin{equation*}
  D^{+} \ln \potential(\flow(t; \theta)) 
  = \frac{D^{+} \potential(\flow(t; \theta))}{\potential(\flow(t; \theta))} \\
  = \frac{D^{+}_{F(\flow(t; \theta))} \potential(\flow(t; \theta))}{\potential(\flow(t; \theta))} \\
  \leq - C \enspace.
\end{equation*}

Note that $D^{+}[- \ln \potential(\flow(t; \theta)) - C t] = - D^{+} \ln \potential(\flow(t; \theta)) - C \geq 0$.
We find by Theorem~2.1 in Appendix~1 of \cite{rouche1977stability} that $- \ln \potential(\flow(t; \theta)) - C t$ is nondecreasing.
This implies that $- \ln \potential(\flow(t; \theta)) + C t \geq - \ln \potential(\flow(0; \theta)) = - \ln \potential(\theta)$, resulting in
\begin{equation*}
  \ln \potential(\flow(t; \theta)) 
  - \ln \potential(\theta) 
  \leq - C t \enspace.
\end{equation*}
Taking the exponential of each side of this inequality, we obtain $\potential(\flow(t; \theta)) / \potential(\theta) \leq \exp(-C t)$.
This ends the proof of the first statement.

To prove the second statement, it suffices to show that $\flow(t; \theta)$ stays within $U$.
We prove this by contradiction.
Assume that for some $\theta \in U$, $\flow(t; \theta)$ leaves $U$ at $t = \tau$ for the first time.
This means that $\potential(\flow(\tau; \theta)) \geq \supU$.
However, since the time derivative of $\potential(\flow(t; \theta))$ is negative, it must hold that $\potential(\flow(t; \theta)) \leq \potential(\theta) < \supU$, which leads to a contradiction.
Hence, $\flow(t; \theta)$ stays within $U$ for any $\theta \in U$ and $t \geq 0$.
This ends the proof of the second statement.

The proof of the last statement is analogous to that of the first and is omitted. \qed

\subsection{Proof of Theorem~\ref{lem:ode-igo}}\label{apdx:lem:ode-igo}
Remember that $t_{n,N} = \sum_{k=0}^{N-1} \deltat_{n+k}$ and $\epsilon_{n,N} = \sum_{k=0}^{N-1} \deltat_{n+k}^2$; if $\deltat$ is time-independent, these are $N\deltat$ and $N \deltat^2$, respectively.
Though the theorem is stated for a time-independent step-size $\deltat$, we use the above notation to cover the proof of Theorem~\ref{lem:ode}.
For the sake of notational simplicity, we drop $(\theta_n)$ from $Q_{\deltat, N}$, $\Theta_{\deltat,N}$, and $\tilde \Theta_{N\deltat}$, and write $Q = Q_{\deltat, N}$.

First, we show the unique existence of $\flow(t; \theta_n)$ for any $t \in [0, t_{n,N}]$.
Since $F$ is a globally Lipschitz function on $\tilde \Theta_{N\deltat}(\theta_n)$ thanks to \asm{P3}, according to Corollary~2.6 in \cite{teschl2012ordinary} we have that the initial value problem \eqref{eq:associated-ode} admits a unique solution for $t \in [0, t_{n,N}]$. Therefore, $\flow(t; \theta_n)$ is uniquely determined.

Second, we show
\begin{equation}
  \E_{n}[\abs{\potential(\theta_{n+k}) - \potential(\flow(t_{n,k}; \theta_n))}] \leq \E_{n}[\norm{\theta_{n+k} - \flow(t_{n,k}; \theta_n)}_{Q}] \enspace.
  \label{eq:proof:ode-igo:main}
\end{equation}
If $\potential$ is differentiable, $D^{+}_{h} \potential(\theta) = \nabla \potential(\theta)^\T h$.
Multiplying both sides by $h$, we obtain $\correct{\abs{D^{+}_{h} \potential(\theta)}}{} \norm{h}_{Q^{-1}} = \norm{(h h^\T) \nabla \potential(\theta)}_{Q^{-1}} \leq \norm{\nabla \potential(\theta)}_{Q^{-1}}$ \correct{for any $h$ such that $\norm{h} = 1$}.
Therefore, $\norm{\nabla \potential(\theta)}_{Q^{-1}} \leq 1$ leads to $\sup_{\norm{h} = 1} \abs{D^+_{h}\potential(\theta)} \norm{h}_{Q^{-1}} \leq 1$. 
Now we derive \eqref{eq:proof:ode-igo:main} from the condition $\sup_{\norm{h} = 1} \abs{D^+_{h}\potential(\theta)} \norm{h}_{Q^{-1}} \leq 1$. 
For any $\theta_a$ and $\theta_b$ in the convex hull of $\Theta_{\deltat,N}(\theta_n) \cup \tilde \Theta_{N\deltat}(\theta_n)$, let $h = (\theta_a - \theta_b) / \norm{\theta_a - \theta_b}$.
Since $1 = \norm{h} \leq \norm{h}_{Q}\norm{h}_{Q^{-1}}$ for any positive-definite symmetric $Q$, we have 
\begin{align*}
  \abs{\potential(\theta_a) - \potential(\theta_b)}
  &\leq\textstyle
  \int_0^{\norm{\theta_a - \theta_b}} \abs{D^+_{h}\potential(\theta_b + t h)} \rmd t
  \\
  &\leq\textstyle
    \int_0^{\norm{\theta_a - \theta_b}} (\abs{D^+_{h}\potential(\theta_b + t h)} \norm{h}_{Q^{-1}}) \norm{h}_{Q} \rmd t
    \\
  &\leq \norm{\theta_a - \theta_b} \norm{h}_{Q}
  = \norm{\theta_a - \theta_b}_{Q} \enspace,    
\end{align*}
where we used \asm{P1} to obtain the last inequality. 
Substituting $\theta_{n+k}$ and $\flow(t_{n,k}; \theta_n)$ for $\theta_a$ and $\theta_b$ and taking the conditional expectation, we have \eqref{eq:proof:ode-igo:main}.

Third, we prove
\begin{multline}
  \textstyle
  \sup_{0 \leq k \leq N}\norm{\theta_{n+k} - \phi(t_{n,k}; \theta_n)}_{Q}
  \\
  \textstyle
  \leq (C_1 + \norm[\big]{\sum\deltat_{n+i} M_{n+i}}_{Q}) \prod (1 + \deltat_{n+i} L_{\deltat,N}(\theta_{n+i}, \theta_n))
  \enspace.\label{eq:gronwall}
\end{multline}
Applying the triangle inequality to the RHS of \eqref{eq:ode-diff}, we obtain
\begin{multline}
\textstyle
\norm{\theta_{n+N} - \phi(t_{n,N}; \theta_n)}_{Q}
\leq\textstyle\sum \deltat_{n+i} \norm[\big]{F(\theta_{n+i}) - F(\phi(t_{n,i}; \theta_n))}_{Q} 
\\
\textstyle + \sum\int_{t_{n,i}}^{t_{n,i+1}} \norm{F(\phi(t_{n,i}; \theta_n)) - F(\phi(\tau; \theta_n))}_{Q} \rmd \tau 
\textstyle+  \norm[\big]{\sum\deltat_{n+i} M_{n+i}}_{Q}
\enspace.
\label{eq:ode-decom}
\end{multline}
Here and in the rest of the proof, sum $\sum$ and product $\prod$ are taken over $i \in \llbracket 0, N-1\rrbracket$.
We will bound each term of the RHS of \eqref{eq:ode-decom}. 
For the first term, we have from \asm{P2} that $\norm{F(\theta_{n+i}) - F(\flow(t_{n,i}; \theta_n))}_{Q} \leq L_{\deltat,N}(\theta_{n+i}, \theta_n)
\norm{\theta_{n+i} - \flow(t_{n,i}; \theta_n)}_{Q}$.
For the second term, from the definition of $F(\theta)$ and \asm{P4}, we have $\norm{F(\theta)}_{Q} = \norm{\E[F_{n}\mid \theta_n = \theta]}_{Q} \leq \E[\norm{F_{n}} _{Q}\mid \theta_n = \theta] \leq \E[\norm{F_{n}}_{Q}^{2} \mid \theta_n=\theta]^{1/2} \leq R(\theta)$.
Moreover, $\norm{\flow(t_{n,i}; \theta_n) - \flow(\tau; \theta_n)}_{Q} = \norm{\int_{t_{n,i}}^{\tau} F(\flow(t; \theta_n)) \rmd t}_{Q} \leq \int_{t_{n,i}}^{\tau} \norm{F(\flow(t; \theta_n))}_{Q} \rmd t \leq \int_{t_{n,i}}^{\tau} R(\flow(t; \theta_n)) \rmd t$.
Using Cauchy's repeated integral formula and \asm{P3}, we find that
\begin{align*}
  \MoveEqLeft[2]\textstyle \int_{t_{n,i}}^{t_{n,i+1}} \norm{F(\flow(t_{n,i}; \theta_n)) - F(\flow(\tau; \theta_n))}_{Q} \rmd \tau
  \\
&\textstyle\leq \Delta_{\tilde{L}}(\deltat, t_{n,N}) \int_{t_{n,i}}^{t_{n,i+1}} \norm{\flow(t_{n,i}; \theta_n) - \flow(\tau; \theta_n)}_{Q} \rmd \tau
  \\
&\textstyle\leq \Delta_{\tilde{L}}(\deltat, t_{n,N})\int_{t_{n,i}}^{t_{n,i+1}}\!\!\int_{t_{n,i}}^{\tau} R(\flow(t; \theta_n)) \rmd t \rmd \tau
  \\
&\textstyle= \Delta_{\tilde{L}}(\deltat, t_{n,N})\int_{t_{n,i}}^{t_{n,i+1}} (t_{n,i+1} - \tau) R(\flow(\tau; \theta_n)) \rmd \tau \enspace,
\end{align*}
where $\Delta_{\tilde{L}}(\deltat, t_{n,N})$ should be understood as $\Delta_{\tilde{L}}(\{\deltat_{n+i}\}_{i=0}^{N-1}, t_{n,N})$ in case of time-dependent step-size.
Therefore, the second term on the RHS of \eqref{eq:ode-decom} is upper-bounded by $\Delta_{\tilde{L}}(\deltat, t_{n,N})\sum\int_{t_{n,i}}^{t_{n,i+1}} (t_{n,i+1} - \tau) R(\flow(\tau; \theta_n)) \rmd \tau = C_1$.
Combining these inequalities with \eqref{eq:ode-decom}, we have
\begin{multline*}
  \textstyle \norm{\theta_{n+N} - \phi(t_{n,N}; \theta_n)}_{Q}
  \\
\leq \textstyle \sum \deltat_{n+i} L_{\deltat,N}(\theta_{n+i}, \theta_n) \norm{\theta_{n+i} - \flow(t_{n,i}; \theta_n)}_{Q} + C_1 + \norm[\big]{\sum\deltat_{n+i} M_{n+i}}_{Q}\enspace.
\end{multline*}
Using the discrete Gronwall inequality \cite{Clark1987279}, we obtain \eqref{eq:gronwall}. 

Fourth, we derive
\begin{equation}
  \textstyle
  \E_n[\norm{\sum\deltat_{n+i} M_{n+i}}_Q^2]^{1/2} 
  \leq \epsilon_{n,N}^{1/2} (1 + K_1^2 t_{n,N}^2 \exp(2 K_1 \cdot t_{n,N}))^{1/2} K_2 R(\theta_n) = C_2 \enspace.
  \label{eq:martingale}
\end{equation}
Since $\{M_{n+i}\}_{i=0}^{N-1}$ are uncorrelated martingale differences, we have $\E[M_k^{\T} Q M_j] = \E[M_k^\T Q\E_k[M_j]] = 0$ for $j > k$, and $\E_i[\norm{M_i}_{Q}^2] \leq \E_i[\norm{F_i}_{Q}^2] \leq R(\theta_i)^2$ thanks to \asm{P4}.
Then, we find $\E_n[\norm{\sum\deltat_{n+i} M_{n+i}}_{Q}^2] = \sum\deltat_{n,i}^2 \E_n[\norm{M_{n+i}}_{Q}^2] \leq \sum\deltat_{n+i}^2 \E_n[R(\theta_{n+i})^2]$.
With \asm{P4}, we have $\E_n[R(\theta_{n+i})^2] \leq K_1^2 \E_n[\norm{\theta_{n+i} - \theta_n}_{Q}^2] + K_2^2 R(\theta_n)^2$.
The expectation in the first term on the RHS is bounded by
\begin{equation}
  \sup_{0\leq k \leq N}\E[\norm{\theta_{n+k} - \theta_n}_{Q}^{2}]^{1/2} \leq K_2 t_{n,N} R(\theta_n) \exp(K_1 t_{n,N}) \enspace.
  \label{eq:loose-bound}
\end{equation}
The derivation of this is as follows: 
By $\mathcal{L}_2$ norm inequality, we have
\begin{align*}
  \textstyle \E_n[\norm{\theta_{n+k} - \theta_n}_Q^{2}]^{1/2}
  = \E_n[\norm{\sum_{i=0}^{k-1}\deltat_{n+i} F_{n+i}}_Q^{2}]^{1/2}
  \leq \sum_{i=0}^{k-1}\deltat_{n+i} \E_n[\norm{F_{n+i}}_Q^2 ]^{1/2} \enspace.    
\end{align*}
With \asm{P4}, each summand is upper-bounded:
\begin{align*}
  \E_n[\norm{F_{n+i}}_Q^2 ]^{1/2} \leq \E_n[R(\theta_{n+i})^2]^{1/2}
  \leq K_1 \E_n[\norm{\theta_{n+i} - \theta_n}_Q^2]^{1/2} + K_2 R(\theta_n) \enspace.
\end{align*}
Hence, we obtain
\begin{equation*}
  \begin{split}
  \E_n[\norm{\theta_{n+k} - \theta_n}_Q^{2}]^{1/2}
  \textstyle\leq
  K_1 \sum_{i=0}^{k-1}\deltat_{n+i} \E_n[\norm{\theta_{n+i} - \theta_n}_Q^2]^{1/2}
  + K_2 t_{n,k} R(\theta_n) 
\enspace.
  \end{split}
\end{equation*}
By the discrete Gronwall inequality and the inequality $\prod (1 + K_1 \deltat_{n_i}) \leq \exp(K_1 t_{n,N})$, we obtain \eqref{eq:loose-bound}.
Therefore,
\begin{align*}
\textstyle \sum\deltat_{n+i}^2 \E_n[R(\theta_{n+i})^2]
  &\leq \epsilon_{n,N} (1 + K_1^2 t_{n,N}^2 \exp(2 K_1 t_{n,N})) K_2^2 R(\theta_n)^2 = C_2^2
\end{align*}

Finally, taking the expectation of \eqref{eq:gronwall} and applying the Schwartz and Minkowski inequalities, we obtain
\begin{align*}
  \MoveEqLeft[2]\textstyle
  \E_n[\sup_{0 \leq k \leq N}\norm{\theta_{n+k} - \phi(t_{n,k}; \theta_n)}_{Q}]
  \\
  &\textstyle
  \leq \E_n[(C_1 + \norm[\big]{\sum\deltat_{n+i} M_{n+i}}_{Q})^2]^{1/2} \E_n[\prod(1 + \deltat_{n+i} L_{\deltat,N}(\theta_{n+i}, \theta_n))^2 ]^{1/2}
  \\
  &\textstyle
    \leq (C_1 + \E_n[\norm[\big]{\sum\deltat_{n+i} M_{n+i}}_{Q}^2]^{1/2}) \E_n[\prod(1 + \deltat_{n+i} L_{\deltat,N}(\theta_{n+i}, \theta_n))^2 ]^{1/2}
    \\
  &\textstyle
    \leq (C_1 + C_2) \Delta_{L}(\alpha, t_{n,N})
  \enspace,
\end{align*}
Where, for the last inequality, we used \eqref{eq:martingale} and \asm{P2}.
Along with \eqref{eq:proof:ode-igo:main}, this completes the proof. 

For the proof of Theorem~\ref{lem:ode}, we can simply set $R(\theta) = K$, where $K > 0$ appears in Theorem~\ref{lem:ode}, and set $Q$ to be the identity matrix.
The Lipschitz continuity assumption in Theorem~\ref{lem:ode} implies \asm{P2} with $\Delta_{L}(x, y) = \exp(L t_{n,N})$ ($\geq \prod(1 + L \deltat_{n+i})$) and \asm{P3} with $\Delta_{\tilde L}(x, y) = L$.
Then, \asm{P4} is satisfied with $K_1 = 0$ and $K_2 = 1$ and we find $C_1 = L K \epsilon_{n,N} / 2$ and $C_2 = K \epsilon_{n,N}^{1/2}$.
Since the function $\potential$ does not come into play for Theorem~\ref{lem:ode}, \asm{P1} is unnecessary. \qed

\subsection{Proof of Corollary~\ref{cor:ode-igo}}\label{apdx:cor:ode-igo}

Since $\sup_{t \in [0, N \deltat]} R(\flow(t; \theta_n)) \leq \Delta_{R}(\deltat, N\deltat) \potential(\theta_n)$,  we obtain
\begin{align*}
 C_1 &\leq \textstyle (N \deltat^2 / 2)\Delta_{\tilde{L}}(\deltat, N\deltat) \Delta_{R}(\deltat, N\deltat) \potential(\theta_n) \enspace,\\
  C_2 &\leq \textstyle K_2 (N \deltat^2)^{1/2} (1 + K_1^2 (N \deltat)^2 \exp(2 K_1 N \deltat))^{1/2} \Delta_{R}(\deltat, N\deltat) \potential(\theta_n)
 \enspace,
\end{align*}%
where we used the formula $\int_{t_{n,i}}^{t_{n,i+1}} (t_{n,i+1} - t)\rmd t = \deltat^2 / 2$.
The RHS of \eqref{eq:ode-igo-bound} is upper-bounded by $\potential(\theta_n) \DeltaAfou(\deltat, N\deltat)$. 
It is easy to see that $\DeltaAfou$ is nondecreasing \wrt each argument and $\DeltaAfou(x, y) \downarrow 0$ as $x \downarrow 0$ for any fixed $y$. \qed

\section{Improved Chebyshev Sum Inequality}

It is often the case that we want to prove the strong positivity of the covariance of two random variables.
According to \cite[Chapter~1]{Thorisson:2000uo}, two non-negatively correlated random variables have a non-negative covariance.
This result is also known as the Chebyshev sum inequality \cite[Theorem~43, Theorem~236]{hardy1952inequalities}.
In the following theorem, we extend the Chebyshev sum inequality to derive a tighter bound and show the strict positivity of the covariance of two non-negatively correlated random variables.

\newcommand{\F}{\mathcal{F}}
\newcommand{\PM}{\mathrm{P}}
\begin{theorem}[Improved Chebyshev sum inequality]
\label{thm:che}
Let $(\X, \F, \PM)$ be a probability space, where $\X$ is the domain of a random variable $X$, $\F$ is a $\sigma$-algebra on $\X$, and $\PM$ is the probability measure on $\F$.
Let $f: \X \to \R$ and $g:\X \to \R$ be $\PM$-integrable functions with finite expected values that are almost surely non-negatively correlated, i.e., $(f(x) - f(y))(g(x) - g(y)) \geq 0$ for any $(x, y) \in \X^2$ except the subset $A \subset \X^2$ for which $(\PM\otimes\PM)(A) = 0$, where $(\PM\otimes\PM)$ is the product measure on $\X^2$. Let $f_1 = f$, $g_1 = g$, and $f_i(x) = \E_{y}[(f_{i-1}(x) - f_{i-1}(y))\ind{f_{i-1}(y) < f_{i-1}(x)}]$, $g_i(x) = \E_{y}[(g_{i-1}(x) - g_{i-1}(y))\ind{g_{i-1}(y) < g_{i-1}(x)}]$ for any integer $i \geq 2$.
Then, for any $K \geq 1$,
\begin{equation}
\E[f(X) g(X)] 
\geq \sum_{i=1}^{K} \E[f_i(X)] \E[g_i(X)]
\enspace,
\label{eq:che:main}
\end{equation}
where $\E[f_i(X)]$ and $\E[g_i(X)]$ can be written by using an i.i.d.~copy $Y$ of $X$ as $\E[f_i(X)] = \frac12 \E[\abs{f_{i-1}(X) - f_{i-1}(Y)}]$ and $\E[g_i(X)] = \frac12 \E[\abs{g_{i-1}(X) - g_{i-1}(Y)}]$ for $i \geq 2$.
\end{theorem}

\begin{proof}
We will show that, for any almost surely non-negatively correlated functions $f$ and $g$, 
\begin{equation}
\E[f(X)g(X)] \geq \E[f(X)]\E[g(X)] + \E[F(X) G(X)] \enspace,
\label{thm:che:0}
\end{equation}
where $F: x \mapsto \E_{Y}[(f(x) - f(Y))\ind{f(Y) < f(x)}]$ and $G: x \mapsto \E_{Y}[(g(x) - g(Y))\ind{g(Y) < g(x)}]$.
Moreover, $F$ and $G$ are almost surely non-negatively correlated.
Once this is proven, the theorem statement follows immediately by repeated application of \eqref{thm:che:0}, recalling that $f_i$ and $g_i$ are non-negative for $i \geq 2$. 

Since $f$ and $g$ are almost surely non-negatively correlated, we have $(f(x) - f(y))(g(x) - g(y)) = \abs{f(x) - f(y)}\abs{g(x) - g(y)}$ almost everywhere w.r.t.~the product measure $\PM \otimes \PM$. Taking the expectation of both sides, we obtain
\begin{equation}
\E[f(X)g(X)] = \E[f(X)]\E[g(X)] + \frac12\E[\abs{f(X) - f(Y)}\abs{g(X) - g(Y)}] \enspace.
\label{eq:che:1}
\end{equation}
Since the second term on the RHS of \eqref{eq:che:1} is nonnegative, the result of \cite[Chapter~1]{Thorisson:2000uo} and \cite[Theorem~43, Theorem~236]{hardy1952inequalities} follows. 

To obtain a tighter bound, we take a closer look at the second term of \eqref{eq:che:1} and find
\begin{align}
\MoveEqLeft[1]
\textstyle \frac12\E[\abs{f(X) - f(Y)}\abs{g(X) - g(Y)}]
\notag \\
&\textstyle = \begin{aligned}[t]
	&\textstyle \frac12\E[\abs{f(X) - f(Y)}\abs{g(X) - g(Y)}\ind{f(Y) < f(X) \wedge g(Y) < g(X)}]
	 \\
	&\textstyle + \frac12\E[\abs{f(X) - f(Y)}\abs{g(X) - g(Y)}\ind{f(Y) > f(X) \wedge g(Y) > g(X)}]
	\end{aligned}
	\notag	\\
&\textstyle = \E[(f(X) - f(Y))\ind{f(Y) < f(X)} (g(X) - g(Y))\ind{g(Y) < g(X)}]
	\enspace.
	\label{eq:che:2}
\end{align}
For the second equality, we used the fact that
the first and the second terms on the LHS are equivalent due to the symmetry.

Given $x$, let $\tilde{f}_x: y \mapsto (f(x) - f(y))\ind{f(y) < f(x)}$ and $\tilde{g}_x: y \mapsto (g(x) - g(y))\ind{g(y) < g(x)}$.
It is easy to see that if $f(y) = f(z)$ then $\tilde{f}_x(y) = \tilde{f}_x(z)$, and if $f(y) > f(z)$ then $\tilde{f}_x(y) \leq \tilde{f}_x(z)$.
Analogously, we have that if $g(y) = g(z)$ then $\tilde{g}_x(y) = \tilde{g}_x(z)$, and if $g(y) > g(z)$ then $\tilde{g}_x(y) \leq \tilde{g}_x(z)$.
This implies that if $(f(y) - f(z))(g(y) - g(z)) \geq 0$, then $(\tilde{f}_x(y) - \tilde{f}_x(z))(\tilde{g}_x(y) - \tilde{g}_x(z)) \geq 0$ holds.
Since $(f(y) - f(z))(g(y) - g(z)) \geq 0$ holds $\PM \otimes \PM$-almost everywhere, so does $(\tilde{f}_x(y) - \tilde{f}_x(z))(\tilde{g}_x(y) - \tilde{g}_x(z)) \geq 0$.
Therefore, $\tilde{f}_x$ and $\tilde{g}_x$ are almost surely non-negatively correlated for each $x \in \X$.

Applying \eqref{eq:che:1}, we obtain
\begin{align}
\MoveEqLeft[2]\E[(f(X) - f(Y))\ind{f(Y) < f(X)} (g(X) - g(Y))\ind{g(Y) < g(X)}] \notag \\
&= \E_{X}[ \E_{Y}[ (f(X) - f(Y))\ind{f(Y) < f(X)} (g(X) - g(Y))\ind{g(Y) < g(X)}]] \notag \\
&\geq \E_{X}[ \E_{Y}[(f(X) - f(Y))\ind{f(Y) < f(X)}] \E_{Y}[(g(X) - g(Y))\ind{g(Y) < g(X)}]] \notag \\
&= \E_{X}[ F(X) G(X)] 
\enspace.
\label{eq:che:3}
\end{align}

We must still prove that $F$ and $G$ are almost surely non-negatively correlated.
It is easy to see from the definition of $F$ that if $f(x) = f(y)$ then $F(x) = F(y)$, and if $f(x) > f(y)$ then $F(x) \geq F(y)$.
Analogously, if $g(x) = g(y)$ then $G(x) = G(y)$, and if $g(x) > g(y)$ then $G(x) \geq G(y)$.
This implies that if $(f(x) - f(y))(g(x) - g(y)) \geq 0$ then $(F(x) - F(y))(G(x) - G(y)) \geq 0$.
Since $(f(x) - f(y))(g(x) - g(y)) \geq 0$ holds $\PM \otimes \PM$-almost everywhere, so does $(F(x) - F(y))(G(x) - G(y)) \geq 0$.
Therefore, $F$ and $G$ are almost surely non-negatively correlated.
This concludes the proof.
\end{proof}

In the above theorem, \eqref{eq:che:main} with $K = 1$ implies the standard Chebyshev sum inequality. With $K = 2$, we have
\begin{equation}
\textstyle \E[f(X) g(X)] \geq \E[f(X)] \E[g(X)] + \frac14 \E[\abs{f(X) - f(Y)}] \E[\abs{g(X) - g(Y)}] \enspace,
\label{eq:che:n2}
\end{equation}
where $X$ and $Y$ are i.i.d.\ random variables. To evaluate the second term on the RHS in \eqref{eq:che:n2}, the next proposition is useful. 

\begin{proposition}
\label{prop:abs-mom}
Let $Y$ and $\tildeY$ be i.i.d.\ random variables on $\R$ with finite mean $\mu_{(1)} = \E[Y]$, variance $\mu_{(2)} = \E[(Y - \mu_{(1)})^2]$ and fourth centered moment $\mu_{(4)} = \E[(Y - \mu_{(1)})^4]$.
Then, the following inequality holds:
\begin{equation*}
\textstyle \E[\abs{Y - \tildeY}] \geq 2\left(\frac{\mu_{(2)}^{3}}{\mu_{(4)} + 3 \mu_{(2)}^2} \right)^{1/2}\enspace.
\end{equation*}
\end{proposition}

\begin{proof}
  We employ a lower-bound technique that is known as the fourth-moment method \cite{Berger:1999wx}.
  For any random variable $X$ that satisfies $0 \leq \E[X^2] \leq \E[X^4] < \infty$ and any $q > 0$, the following inequality holds:
\begin{equation*}
\textstyle \E[\abs{X}] \geq \frac{3^{3/2}}{2 q^{1/2}}\left( \E[X^{2}] - \frac{\E[X^{4}]}{q} \right) \enspace.
\end{equation*}
Letting $q = 3 \E[X^{4}] / \E[X^{2}]$, we have
\begin{equation}\label{eq:berger}
\textstyle \E[\abs{X}] \geq \frac{\E[X^{2}]^{3/2}}{\E[X^{4}]^{1/2}}\enspace.
\end{equation}

In our case, $X$ has the form $Y - \tildeY$, where $Y$ and $\tildeY$ are i.i.d.\ random variables.
Provided that $\E[Y^{4}] < \infty$, we have $\E[Y^{2}] \leq \E[Y^{4}]^{1/2} < \infty$ by Schwarz inequality. Since $Y$ and $\tildeY$ are i.i.d., we have that 
\begin{gather}
\E[\lvert Y - \tildeY \rvert^{2}] = 2 \mu_{(2)} 
\label{eq:2mom}\\
\E[\lvert Y - \tildeY \rvert^{4}] = 2 \mu_{(4)} + 6 \mu_{(2)}^{2} \enspace. \label{eq:4mom}
\end{gather}
From \eqref{eq:2mom}, \eqref{eq:4mom}, and \eqref{eq:berger}, we have the proposition statement.
\end{proof}

\section{Proofs of the Linear Convergence of Adaptive ES}\label{apdx:es}


\subsection{Proof of Lemma~\ref{lemma:apdx:j}}\label{apdx:lemma:apdx:j}

To prove Lemma~\ref{lemma:apdx:j}, we use the following result: 
\begin{proposition}
\label{prop:w-diff}
Let $x, y \sim P_{\theta_n}$ be i.i.d. Then, 
\begin{equation}
\E_{x, y \sim P_{\theta_n}}[\abs{u(\qftnle{n}(f(x))) - u(\qftnle{n}(f(y)))}] 
= 2 \sum_{k=1}^{\lambda} \ww_{k} \left(1 - \frac{2k}{\lambda+1}\right)
\enspace.
\label{eq:w-diff}
\end{equation}
\end{proposition}
\begin{proof}
  Since $\qftneq{n}(f(x)) = 0$ holds $P_{\theta_n}$-almost-surely, the distribution of $f(x)$ with $x \sim P_{\theta_{n}}$ is absolutely continuous \wrt the Lebesgue measure. Hence, $\qftnle{n}: \R \to [0, 1]$ is a continuous function. From the well-known fact \cite[Theorem~2.1]{Devroye1986} that the distribution of $P(X)$ with a random variable $X$ having a continuous cumulative density function $P$ is uniform within $[0, 1]$, we know that $\qftnle{n}(f(x))$ with $x \sim P_{\theta_{n}}$ is uniformly distributed in $[0, 1]$. Using this fact and the monotonicity of $u$ to remove the absolute value, the LHS of \eqref{eq:w-diff} can be written as
  $\E_{x, y \sim P_{\theta_n}}[\abs{u(\qftnle{n}(f(x))) - u(\qftnle{n}(f(y)))}]  = 2 \int_{0}^{1}\int_{x}^{1} (u(x) - u(y)) \rmd y \rmd x $. Note that $\int_{x}^{1} u(x) \rmd y = (1-x)u(x)$ and $\int_{0}^{1}\int_{x}^{1} u(y) \rmd y \rmd x = \int_{0}^{1} x u(x) \rmd x$. Using these equalities, we obtain $\E_{x, y \sim P_{\theta_n}}[\abs{u(\qftnle{n}(f(x))) - u(\qftnle{n}(f(y)))}] = 2 \int_{0}^{1} (1 - 2 x) u(x) \rmd x$. 
Substituting \eqref{eq:binom} and using the formula $\int_0^1 x^k (1 - x)^{n - k} \rmd x = (n+1)^{-1} \binom{n}{k}^{-1}$, we obtain \eqref{eq:w-diff}.
\end{proof}

Now we prove Lemma~\ref{lemma:apdx:j}.
Note that $F$ can be written as \eqref{eq:meanfield-s} with $g$ as defined in \eqref{eq:adaptES} and $\Ratio = 1$.
We have $\inner{\nabla V(\theta)}{g(x; \theta)} = \norm{x}^2 - \E[\norm{x}^2]$. 
Therefore,
\begin{align}
  \inner*{\nabla V(\theta) }{ F(\theta) } &= \int u(q^{<}_{\theta}(f(x))) (f(x) - \E[f(x)]) p_{\theta}(x) \mathrm{d}x \enspace. \label{eq:lamma:j:1}
\end{align}
Since $u$ is monotonically decreasing under \asm{W2} and $f(x) = \nu(\norm{x})$ with a decreasing function $\nu$, $u(q^{<}_{\theta}(f(x)))$ and $\norm{x}^2 - \E[\norm{x}^2]$ are negatively correlated.
Then, by Theorem~\ref{thm:che} with $K = 2$ and Propositions~\ref{prop:abs-mom} and~\ref{prop:w-diff}, we can upper-bound the RHS of the above equality by
\begin{align}
  - \left( \sum_{k=1}^{\lambda} w_k \left(1 - \frac{2k}{\lambda + 1}\right)\right) \left(\frac{\E[(\norm{x}^2 - \E[\norm{x}^2])^2]^{3}}{\E[(\norm{x}^2 - \E[\norm{x}^2])^4] + 3 \E[(\norm{x}^2 - \E[\norm{x}^2])^2]^2} \right)^{1/2} \enspace,
  \label{eq:apdx:j:2}
\end{align}
where we used the fact that $q^f_{\theta}(f(x))$ is uniformly distributed on $[0, 1]$ under $x \sim p(\cdot; \theta)$.
Here, $\norm{x}^2 / \sigma^2$ is non-centrally $\chi^2$ distributed with $d$ degrees of freedom and non-centrality parameter $\xi(\theta)^2$. The centered moments of this are well-known, so we have
\begin{align*}
  \E[(\norm{x}^2 - \E[\norm{x}^2])^2] &= \sigma^4 \cdot (2(d + 2 \xi(\theta)^2)) \\
  \E[(\norm{x}^2 - \E[\norm{x}^2])^4] &= \sigma^8 \cdot (12(d + 2\xi(\theta)^2)^2 + 48(d + 4\xi(\theta)^2))\enspace.
\end{align*}
Therefore,
\begin{multline*}
  \frac{\E[(\norm{x}^2 - \E[\norm{x}^2])^2]^{3}}{\E[(\norm{x}^2 - \E[\norm{x}^2])^4] + 3 \E[(\norm{x}^2 - \E[\norm{x}^2])^2]^2}
  \\
  = \frac{\sigma^4}{3} \frac{ (d + 2 \xi(\theta)^2)^3 }{ (d + 2\xi(\theta)^2)^2 + 2(d + 4\xi(\theta)^2) }
  > \frac{\sigma^4}{3} \frac{ d + 2 \xi(\theta)^2 }{ 1 + 4/d }
  \geq \frac{\sigma^4}{3} \frac{ d + \xi(\theta)^2 }{ 1 + 4/d }
  = \frac{\sigma^2}{3} \frac{ V(\theta) }{ 1 + 4/d }
  \enspace.
\end{multline*}
Substituting the RMS of the above inequality into \eqref{eq:apdx:j:2} completes the proof.
\qed

\subsection{Proof of Lemma~\ref{lemma:apdx:ratio}}\label{apdx:lemma:apdx:ratio}

For $x \sim P_\theta$, let $z = (x - m) / \sigma$ and $z_{m} = \inner{m/\norm{m}}{z}$.
  First, we show that for any $\epsilon > 0$ there exists an $\eta_\epsilon > 0$ such that for any $\theta \in \Theta$ satisfying $\xi(\theta) \geq \eta_\epsilon$,
  \begin{align}
    \E\big[\abs*{ q^{<}_{\theta}(f(x)) - \Phi\left( z_e \right)}^2\big]^{1/2} \leq \epsilon + 2^{-1}\cdot(d / \pi)^{1/2} \cdot \xi(\theta)^{-1} \enspace.\label{eq:lem:quantile:1}
  \end{align}
  Let 
  \begin{align*}
    Q^{<}_{\theta}(t) &:= \Pr_{y\sim P_\theta} \left[ \left( \norm{y}^2 - (\norm{m}^2 + d \sigma^2) \right) / ( 2 \sigma \norm{m}) < t \right] \enspace.
  \end{align*}
  Then, 
  \begin{align*}
    q^{<}_{\theta}(f(x)) &= Q^{<}_{\theta}\left(\frac{ \norm{x}^2 - (\norm{m}^2 + d \sigma^2)}{ 2\sigma \norm{m}} \right)
                           = Q^{<}_{\theta}\left(
                           z_m + \frac{\sigma (\norm{z}^2-d) }{2\norm{m}} 
                           \right)\enspace.
  \end{align*}
  Lemma~8 of \cite{Akimoto2020tcs} shows that
  $\lim_{\norm{m} / \sigma \to \infty} \sup_{t \in [0, 1]} \abs*{ Q^{<}_{\theta}\left( t \right) - \Phi\left( t \right) } = 0$.
  That is, for any $\epsilon > 0$, there exists $\eta_{\epsilon} > 0$ such that for any $\theta \in \Theta$ satisfying $\xi(\theta) \geq \eta_{\epsilon}$ and any $x = m + \sigma z$,
  \begin{align*}
    \abs*{ q^{<}_{\theta}(f(x)) - \Phi\left( z_m + \xi(\theta)^{-1} (\norm{z}^2 - d)/2 \right) } \leq \epsilon \enspace.
  \end{align*}
  Because of the Lipschitz continuity of $\Phi$ with constant $\frac{1}{(2\pi)^{1/2}}$, we have 
  \begin{align*}
   \abs*{ \Phi\left( z_m + \xi(\theta)^{-1} (\norm{z}^2 -d) /2 \right) - \Phi\left( z_m \right) } \leq \xi(\theta)^{-1} (\norm{z}^2-d) / (8\pi)^{1/2}  \enspace.
  \end{align*}
  Then, by using the Minkowski inequality as well as $\E[(\norm{z}^2-d)^2] = 2d$, we obtain \eqref{eq:lem:quantile:1}. 

  Next, we show that
  \begin{align}
    \inner*{\nabla \xi(\theta) }{ \sigma \int u\left(\Phi\left( z_m \right)\right) \begin{bmatrix} z \\ (\norm{z}^2 - d) / (2d) \end{bmatrix}  \psi(z) \mathrm{d}z } \leq -\frac{\bar{C}_\sigma}{2d} \xi(\theta) \enspace.\label{eq:lem:quantile:2}
  \end{align}
  The LHS is
  \begin{align*}
    \int u\left(\Phi\left( z_m \right)\right) z_m  \psi(z) \mathrm{d}z 
    - \xi(\theta) \int u\left(\Phi\left( z_m \right)\right) ((\norm{z}^2 - d) / (2d))  \psi(z) \mathrm{d}z \enspace.
  \end{align*}
  Because $u$ is decreasing under \asm{W2} and because $\Phi\left( z_m \right)$ and $z_m$ are positively correlated, the first term is upper-bounded by $\int u\left(\Phi\left( z_m \right)\right) \psi(z) \mathrm{d}z \int z_m \psi(z) \mathrm{d}z = 0$ in light of Theorem~\ref{thm:che} with $K = 1$. 
  The integral in the second term is
  \begin{align*}
    \MoveEqLeft[2]\textstyle\int u\left(\Phi\left( z_m \right)\right) (\norm{z}^2 - d)/ (2d) \psi(z) \mathrm{d}z \\
    &\textstyle =\int u\left(\Phi\left( z_m \right)\right) (z_m^2 - 1)/ (2d) \psi(z) \mathrm{d}z \\
    &\textstyle\quad + \int u\left(\Phi\left( z_m \right)\right) ((\norm{z}^2 - z_m^2) - (d-1))/(2d) \psi(z) \mathrm{d}z \\
    &=\bar{C}_\sigma / (2d)
      \enspace,
  \end{align*}
  where, for the last equality, we used the facts that $z_m$ and $\norm{z}^2 - z_m^2$ are independent and the expectation of the latter is $d - 1$.
  Hence, we obtain \eqref{eq:lem:quantile:2}.

  Finally, we obtain
  \begin{multline*}
    \inner*{ \nabla \xi(\theta) }{ F(\theta) } 
    =
      \inner*{\nabla \xi(\theta) }{ \sigma \int u\left(\Phi\left( z_m \right)\right) \begin{bmatrix} z \\ (\norm{z}^2 - d) / (2d) \end{bmatrix}  \psi(z) \mathrm{d}z }
    \\
    + \sigma \int \left( u\left(q^{<}_{\theta}( f(m + \sigma z) )\right)- u\left(\Phi\left( z_m \right)\right)\right) \inner*{\nabla \xi(\theta) }{ \begin{bmatrix} z \\ (\norm{z}^2 - d) / (2d) \end{bmatrix} } \psi(z) \mathrm{d}z  \enspace.
  \end{multline*}
  The first term is upper-bounded by $-(\bar{C}_\sigma / (2d)) \xi(\theta)$, as shown in \eqref{eq:lem:quantile:2}.
  By using the Schwarz inequality as well as \eqref{eq:lem:quantile:1} and the $L_u$-Lipschitz continuity of $u$, the second term is upper-bounded by
  \begin{align*}
   L_u \cdot \left(\epsilon + \left(\frac{d}{4\pi}\right)^{1/2} \xi(\theta)^{-1} \right)\left( \int \left(z_m - \frac{\norm{z}^2 - d}{2d} \xi(\theta) \right)^2 \psi(z) \mathrm{d}z \right)^{1/2} \enspace.
  \end{align*}
  By using the Minkowski inequality as well as $\E[z_m^2] = 1$ and $\E[(\norm{z}^2 - d)^2] = 2d$, we can further bound it by 
  \begin{multline*}
    L_u \cdot\left(\epsilon + \left(\frac{d}{4\pi}\right)^{1/2} \xi(\theta)^{-1} \right)\left( 1 + \frac{1}{(2d)^{1/2}} \xi(\theta) \right)
    \\= L_u \cdot\left[\frac{\epsilon}{(2d)^{1/2}} \xi(\theta) + \epsilon + \left(\frac{d}{8 \pi}\right)^{1/2} + \left(\frac{d}{4\pi}\right)^{1/2} \xi(\theta)^{-1}\right]
    \enspace.
  \end{multline*}
 Then, we obtain 
  \begin{align*}
    \frac{\inner*{ \nabla \xi(\theta) }{ F(\theta) }}{ \xi(\theta) } 
    \leq - \left( \frac{\bar{C}_\sigma}{ 2d } - \frac{L_u \cdot\epsilon}{(2d)^{1/2}} \right) + \left( \epsilon + \left(\frac{d}{8 \pi}\right)^{1/2}\right) \frac{L_u}{\xi(\theta)} + \left(\frac{d}{4\pi}\right)^{1/2} \frac{L_u}{\xi(\theta)^2} \enspace.
  \end{align*}
  Since $\epsilon > 0$ is arbitrary, for any $C_\sigma \in (0, \bar{C}_\sigma]$ we can take a sufficiently small $\epsilon > 0$ and a sufficiently large $\eta \geq \eta_{\epsilon}$ for the RHS of the above inequality to be upper-bounded by $- C_\sigma / 2d$ under $\xi(\theta) \geq \eta$. This completes the proof.
\qed

\subsection{Proof of Lemma~\ref{lemma:es:variance}}\label{apdx:lemma:es:variance}

For simplicity, we write $W_{i} = w(i; x_{n,1}, \dots, x_{n,\lambda})$ and $Z_{i} = (x_{n,i} - m_n) / \sigma_n$. Then, 
\begin{align*}
  \E_n[\norm{F_n}^2]
  &= \sigma_n^2 \E\left[\norm*{ \sum_{i=1}^{\lambda} W_{i} \begin{bmatrix} Z_{i} \\ \frac{\norm{Z_{i}}^2 - d}{2d}\end{bmatrix}  }^2 \right]
  \leq \sigma_n^2 \E\left[\left( \sum_{i=1}^{\lambda} \abs{W_{i}} \norm*{\begin{bmatrix} Z_{i} \\ \frac{\norm{Z_{i}}^2 - d}{2d}\end{bmatrix}  }\right)^2 \right]
  \\
  &\leq \sigma_n^2 \E\left[\sum_{i=1}^{\lambda} \abs{W_{i}}^2 \sum_{i=1}^{\lambda} \norm*{\begin{bmatrix} Z_{i} \\ \frac{\norm{Z_{i}}^2 - d}{2d}\end{bmatrix}  }^2 \right]
  \leq \sigma_n^2 \left(\sum_{i=1}^{\lambda} w_i^2\right)  \E\left[\sum_{i=1}^{\lambda} \norm*{\begin{bmatrix} Z_{i} \\ \frac{\norm{Z_{i}}^2 - d}{2d}\end{bmatrix}  }^2 \right] 
  \\
  &\leq \sigma_n^2 \lambda \left(\sum_{i=1}^{\lambda} w_i^2\right) \left(d + \frac{1}{2d}\right) \enspace.
\end{align*}
This completes the proof.\qed

\subsection{Proof of Lemma~\ref{lemma:apdx:ratio-bound}}\label{apdx:lemma:apdx:ratio-bound}
For simplicity, we write $W_i = w(i; x_{n,1}, \dots, x_{n,\lambda})$ and $Z_i = (x_{n,i} - m_n) / \sigma_n$.
Under \asm{W1}, we have $\sigma_{n+1} \geq \sigma_{n} (1 - \alpha / 2)$ for any $n$. Hence, $\sigma_{k} \geq \sigma_0 (1 - \alpha / 2)^{k} \geq (1 - \alpha / 2)^{N}$ for any $k \in \llbracket 0, N\rrbracket$. 
  Next, observe
  \begin{align*}
    \frac{\norm{m_{n+1} - m_n}}{\sigma_{n+1}}
    &= \frac{\norm*{\alpha  \sum_{i=1}^{\lambda} W_i Z_i }}{1 + \alpha \sum_{i=1}^{\lambda} W_i (\norm{Z_i}^2 - d) / (2d)}
    \\
    &\leq \frac{\alpha (\sum_{i=1}^{\lambda} W_i \norm*{ Z_i }^2)^{1/2}}{1 + \alpha \sum_{i=1}^{\lambda} W_i (\norm{Z_i}^2 - d) / (2d)}
      \leq \left(\frac{d \alpha}{2 - \alpha}\right)^{1/2}
      \enspace.
  \end{align*}
  Here, for the last inequality, we used that the LHS takes its maximal value (which is equal to the RHS) at $\sum_{i=1}^{n} W_i \norm*{ Z_i }^2 = (2 - \alpha) d / \alpha$. Then,
  \begin{multline*}
    \frac{\norm{m_{n+1}}}{\sigma_{n+1}}
    \leq  \frac{\norm{m_{n+1} - m_n}}{\sigma_{n+1}} + \frac{\norm{m_n}}{\sigma_{n+1}}
    \\
    \leq \frac{\norm{m_{n+1} - m_n}}{\sigma_{n+1}} + \left(1 - \frac{\alpha}{2}\right)^{-1}\frac{\norm{m_n}}{\sigma_{n}}
    \leq
    \left(\frac{d \alpha}{2 - \alpha}\right)^{1/2}
    + \left(1 - \frac{\alpha}{2}\right)^{-1}\frac{\norm{m_n}}{\sigma_{n}}
    \enspace.
  \end{multline*}
  Then, it is easy to see that for any $k \in \llbracket 0, N \rrbracket$,
  \begin{align*}
    \frac{\norm{m_{k}}}{\sigma_{k}}
    &\leq \left(1 - \frac{\alpha}{2}\right)^{-N} \frac{\norm{m_{0}}}{\sigma_{0}}
      + \frac{(1 - \alpha/2)^{-N} - 1}{(1 - \alpha/2)^{-1} - 1}  \left(\frac{d \alpha}{2 - \alpha}\right)^{1/2}  \enspace.
  \end{align*}
  This completes the proof for $\Theta_{\alpha,N}(\theta_0)$.

  Next, we supply the proof for $\tilde{\Theta}_{N\alpha}(\theta_0)$. 
  Note $\norm{F(\theta)} = \norm{\E[F_n \mid \theta_n = \theta]} \leq \E[\norm{F_n}^2 \mid \theta_n = \theta]^{1/2} \leq C_{\ref{lemma:es:variance}}^{1/2} \cdot \sigma$. Then, it is easy to see that $\sigma_0 \exp(- C_{\ref{lemma:es:variance}}^{1/2} \cdot t) \leq \flow_{\sigma}(t; \theta_0) \leq \sigma_0 \exp( C_{\ref{lemma:es:variance}}^{1/2} \cdot t)$ for any $t \geq 0$. Hence, $\sigma_0 \exp(- C_{\ref{lemma:es:variance}}^{1/2} \cdot N \alpha) \leq \flow_{\sigma}(t; \theta_0) \leq \sigma_0 \exp( C_{\ref{lemma:es:variance}}^{1/2} \cdot N\alpha)$ for all $t \in [0, N\alpha]$.
Observe that
\begin{align*}
  \frac{\mathrm{d}}{\mathrm{d} t} \xi(\flow(t; \theta_0)) \leq \norm{\nabla \xi(\flow(t; \theta_0))} \cdot \norm{F(\flow(t; \theta_0))}
  \leq C_{\ref{lemma:es:variance}}^{1/2} (1 + \xi(\flow(t; \theta_0)))\enspace,
\end{align*}
where we used $\norm{\nabla \xi(\theta)} = \norm*{( m / \norm{m}, -\norm{m} / \sigma)} / \sigma = (1 + \xi(\theta)^2)^{1/2} / \sigma \leq (1 + \xi(\theta)) / \sigma$ and $\norm{F(\theta)} \leq C_{\ref{lemma:es:variance}}^{1/2}\cdot \sigma$. Then, for any $t \in [0, N\alpha]$, we have $\xi(\flow(t; \theta_0)) \leq (1 + \xi(\theta_0)) \cdot \exp(C_{\ref{lemma:es:variance}}^{1/2}\cdot N \alpha) - 1$. This completes the proof for $\tilde{\Theta}_{N\alpha}(\theta_0)$.     
\qed

\subsection{Proof of Lemma~\ref{lemma:apdx:lip}}\label{apdx:lemma:apdx:lip}
  Note first that $g(x; \theta) = \tilde\nabla_\theta \ln p_\theta(x) = \mathcal{I}(\theta)^{-1} \nabla_\theta \ln p_\theta(x)$. 
  It is enough to show that the operator norm of $\nabla_{\theta} F(\theta)$ is bounded by $L_{\ref{lemma:apdx:lip}}$, where
\begin{align*}
  \nabla_{\theta} F(\theta)
  &\textstyle= \nabla_{\theta}\int u(q^{<}_{\theta}(f(x))) \tilde\nabla_\theta \ln(p_\theta(x)) p_\theta(x) \mathrm{d}x\\
  &\textstyle= \int \tilde\nabla_\theta \ln(p_\theta(x)) \nabla_{\theta} u(q^{<}_{\theta}(f(x)))^\T  p_\theta(x) \mathrm{d}x\\
  &\textstyle\quad + \int u(q^{<}_{\theta}(f(x))) \nabla_{\theta} (\tilde\nabla_\theta \ln(p_\theta(x))) p_\theta(x) \mathrm{d}x \\
  &\textstyle\quad + \int u(q^{<}_{\theta}(f(x))) \tilde\nabla_\theta \ln(p_\theta(x)) \nabla_{\theta} \ln(p_\theta(x))^\T p_\theta(x) \mathrm{d}x  \enspace.  
\end{align*}

For the first term, noting that
\begin{multline*}
  \norm{\nabla_{\theta}q^{<}_{\theta}(f(x))}
  \textstyle= \norm{\int_{f(y) \leq f(x)} \nabla_{\theta}\ln(p_\theta(y)) p_\theta(y) \mathrm{d}y}
    \textstyle\leq \int_{f(y) \leq f(x)} \norm{\nabla_{\theta}\ln(p_\theta(y))} p_\theta(y) \mathrm{d}y
  \\
  \textstyle\leq \int \norm*{\nabla_{\theta}\ln(p_\theta(y))} p_\theta(y) \mathrm{d}y
  \textstyle\leq \left(\int \norm*{\nabla_{\theta}\ln(p_\theta(y))}^2 p_\theta(y) \mathrm{d}y\right)^{1/2}
  \textstyle= \Tr(\mathcal{I}(\theta))^{1/2} \enspace,
\end{multline*}
we have
\begin{align*}
  &\textstyle\norm{\int \tilde\nabla_{\theta}\ln(p_\theta(y)) \nabla_{\theta} u(q^{<}_{\theta}(f(x)))^\T p_\theta(x) \mathrm{d}x}\\
  &\textstyle= \norm{\int u'(q^{<}_{\theta}(f(x))) \tilde\nabla_{\theta}\ln(p_\theta(y)) \nabla_{\theta}q^{<}_{\theta}(f(x))^\T p_\theta(x) \mathrm{d}x}\\
  &\textstyle\leq \max_{s \in [0, 1]}\abs{u'(s)} \cdot \int \norm{\nabla_{\theta}q^{<}_{\theta}(f(x))} \norm{\tilde\nabla_{\theta}\ln(p_\theta(y))} p_\theta(x) \mathrm{d}x\\
  &\textstyle\leq \max_{s \in [0, 1]}\abs{u'(s)} \cdot \Tr(\mathcal{I}(\theta))^{1/2} \cdot \int \norm{\tilde\nabla_{\theta}\ln(p_\theta(y))} p_\theta(x) \mathrm{d}x\\
  &\textstyle\leq \max_{s \in [0, 1]}\abs{u'(s)} \cdot \Tr(\mathcal{I}(\theta))^{1/2} \cdot \left(\int \norm{\tilde\nabla_{\theta}\ln(p_\theta(y))}^2 p_\theta(x) \mathrm{d}x \right)^{1/2}\\
  &\textstyle= \max_{s \in [0, 1]}\abs{u'(s)} \cdot \Tr(\mathcal{I}(\theta))^{1/2} \cdot \Tr(\mathcal{I}(\theta)^{-1})^{1/2} 
  \textstyle= \max_{s \in [0, 1]}\abs{u'(s)} \cdot (3d^2+3/2)^{1/2} \enspace,
\end{align*}
where we used $\mathcal{I}(\theta) = \diag(I_d, 2d) / \sigma^2$. 

For the second term,
\begin{align*}\textstyle
  \nabla_{\theta} (\tilde\nabla_\theta \ln(p_\theta(x))) = \begin{bmatrix} - I & 0 \\ - (x - m)/ (d \sigma) & - \norm{x - m}^2 / (2 d \sigma^2) - 1/2\end{bmatrix} \enspace,
\end{align*}
implying that $\nabla_{\theta} (\tilde\nabla_\theta \ln(p_\theta(x)))$ is negative definite.
Then,
\begin{align*}
  \MoveEqLeft[2]\textstyle\norm{\int u(q^{<}_{\theta}(f(x))) \nabla_{\theta} (\tilde\nabla_\theta \ln(p_\theta(x))) p_\theta(x) \mathrm{d}x}\\
  &\textstyle\leq \max_{s \in [0, 1]} \abs{u(s)} \cdot \norm{\int \nabla_{\theta} (\tilde\nabla_\theta \ln(p_\theta(x))) p_\theta(x) \mathrm{d}x}\\
  &\textstyle= \max_{s \in [0, 1]} \abs{u(s)} \cdot \norm*{\begin{bmatrix} - I_d & 0 \\ 0 & - 1 \end{bmatrix}}
  \textstyle= \max_{s \in [0, 1]} \abs{u(s)} \enspace.
\end{align*}

For the third term, using the fact that $\nabla_{\theta} \ln(p_\theta(x)) \nabla_{\theta}\ln(p_\theta(x))^\T$ is positive definite, we have
\begin{align*}
  &\textstyle\norm*{\int u(q^{<}_{\theta}(f(x))) \mathcal{I}(\theta)^{-1} \nabla_{\theta} \ln(p_\theta(x)) \nabla_{\theta}\ln(p_\theta(x))^\T p_\theta(x) \mathrm{d}x}\\
  &\textstyle\leq \max_{s \in [0, 1]} \abs{u(s)} \cdot \norm*{\int \mathcal{I}(\theta)^{-1} \nabla_{\theta} \ln(p_\theta(x)) \nabla_{\theta}\ln(p_\theta(x))^\T p_\theta(x) \mathrm{d}x}\\
  &\textstyle= \max_{s \in [0, 1]} \abs{u(s)} \cdot \norm*{\mathcal{I}(\theta)^{-1} \mathcal{I}(\theta)}
  \textstyle= \max_{s \in [0, 1]} \abs{u(s)} \cdot \norm*{I_\dimtheta}
  \textstyle= \max_{s \in [0, 1]} \abs{u(s)} \enspace.
\end{align*}
Noting that $\abs{u(s)} \leq \lambda \abs{w_i}$ and $\abs{u'(s)} \leq \lambda (\lambda - 1) \max_{i}\abs{- w_i + w_{i+1}}$, the proof is complete.
\qed

\section*{Acknowledgements}
The authors would like to thank the anonymous reviewers for their valuable comments and
suggestions.
This work is partially supported by JSPS KAKENHI Grant Number 19H04179.


\end{document}